\documentclass{article}
\usepackage{amsmath}
\usepackage{amssymb}
\usepackage{amsthm}
\usepackage{cite}
\usepackage[arrow,curve,matrix,arc,2cell]{xy}
\UseAllTwocells
\begin{document}
%%%%%%%%%%%%%%%%%%%%%%%%%%%%%%%%%%%%%%%%%%%%%%%%%%%%%%%%%%%%%%%%%%%%%%%%
%%%%%%%%%%%%%%%%%%%%%%%%%%     Macros      %%%%%%%%%%%%%%%%%%%%%%%%%%%%%
%%%%%%%%%%%%%%%%%%%%%%%%%%%%%%%%%%%%%%%%%%%%%%%%%%%%%%%%%%%%%%%%%%%%%%%%
\def\e#1\e{\begin{equation}#1\end{equation}}
\def\ea#1\ea{\begin{align}#1\end{align}}
\def\eq#1{{\rm(\ref{#1})}}
\theoremstyle{plain}% default
\newtheorem{thm}{Theorem}[section]
\newtheorem{lem}[thm]{Lemma}
\newtheorem{prop}[thm]{Proposition}
\newtheorem{cor}[thm]{Corollary}
\newtheorem{conj}[thm]{Conjecture}
\newtheorem{quest}[thm]{Question}
\theoremstyle{definition}
\newtheorem{dfn}[thm]{Definition}
\newtheorem{ex}[thm]{Example}
\newtheorem{rem}[thm]{Remark}
\newtheorem{conv}[thm]{Convention}
\newtheorem{property}[thm]{Property}
\numberwithin{figure}{section} \numberwithin{equation}{section}
\def\dim{\mathop{\rm dim}\nolimits}
\def\supp{\mathop{\rm supp}}
\def\codim{\mathop{\rm codim}\nolimits}
\def\depth{\mathop{\rm depth}\nolimits}
\def\vdim{\mathop{\rm vdim}\nolimits}
\def\End{\mathop{\rm End}}
\def\Sh{\mathop{\rm Sh}}
\def\Ker{\mathop{\rm Ker}}
\def\Coker{\mathop{\rm Coker}}
\def\Sim{\mathop{\rm Sim}}
\def\Sing{\mathop{\rm Sing}}
\def\GL{\mathop{\rm GL}}
\def\GR{\mathop{\rm GR}}
\def\Im{\mathop{\rm Im}}
\def\Stab{\mathop{\rm Stab}\nolimits}
\def\rank{\mathop{\rm rank}}
\def\Hom{\mathop{\rm Hom}\nolimits}
\def\Ho{\mathop{\rm Ho}\nolimits}
\def\Ext{\mathop{\rm Ext}\nolimits}
\def\bHom{\mathop{\bf Hom}\nolimits}
\def\Aut{\mathop{\rm Aut}}
\def\coh{\mathop{\rm coh}}
\def\qcoh{\mathop{\rm qcoh}}
\def\vect{\mathop{\rm vect}}
\def\vqcoh{\mathop{\rm vqcoh}}
\def\vvect{\mathop{\rm vvect}}
\def\id{{\mathop{\rm id}\nolimits}}
\def\top{{\kern.05em\rm top}}
\def\red{{\mathop{\rm red}}}
\def\Iso{{\rm Iso}}
\def\inc{{\rm inc}}
\def\Obj{{\rm Obj}}
\def\Spec{\mathop{\rm Spec}\nolimits}
\def\MSpec{\mathop{\rm MSpec}\nolimits}
\def\Ho{{\mathop{\rm Ho}}}
\def\Top{{\mathop{\bf Top}}}
\def\CRings{{\mathop{\bf C^{\bs\iy}Rings}}}
\def\CRingsfg{{\mathop{\bf C^{\bs\iy}Rings^{fg}}}}
\def\CRingsfp{{\mathop{\bf C^{\bs\iy}Rings^{fp}}}}
\def\CRingsfa{{\mathop{\bf C^{\bs\iy}Rings^{fa}}}}
\def\CRingsgo{{\mathop{\bf C^{\bs\iy}Rings^{go}}}}
\def\CRS{{\mathop{\bf C^{\bs\iy}RS}}}
\def\LCRS{{\mathop{\bf LC^{\bs\iy}RS}}}
\def\DLCRS{{\mathop{\bf DLC^{\bs\iy}RS}}}
\def\ACSch{{\mathop{\bf AC^{\bs\iy}Sch}}}
\def\ACSchfp{{\mathop{\bf AC^{\bs\iy}Sch^{fp}}}}
\def\ACSchfa{{\mathop{\bf AC^{\bs\iy}Sch^{fa}}}}
\def\ACSchgo{{\mathop{\bf AC^{\bs\iy}Sch^{go}}}}
\def\CSch{{\mathop{\bf C^{\bs\iy}Sch}}}
\def\CSchlfp{{\mathop{\bf C^{\bs\iy}Sch^{lfp}}}}
\def\CSchlf{{\mathop{\bf C^{\bs\iy}Sch^{lf}}}}
\def\CSchlfssc{{\mathop{\bf C^{\bs\iy}Sch^{lf}_{ssc}}}}
\def\CSchlg{{\mathop{\bf C^{\bs\iy}Sch^{lg}}}}
\def\hCSchlfssc{{\mathop{\bf \hat{C}^{\bs\iy}Sch^{lf}_{ssc}}}}
\def\bCSch{{\mathop{\bf{\bar C}^{\bs\iy}Sch}}}
\def\bCSchlfp{{\mathop{\bf{\bar C}^{\bs\iy}Sch^{lfp}}}}
\def\bCSchlf{{\mathop{\bf{\bar C}^{\bs\iy}Sch^{lf}}}}
\def\bCSchlg{{\mathop{\bf{\bar C}^{\bs\iy}Sch^{lg}}}}
\def\CSta{{\mathop{\bf C^{\bs\iy}Sta}}}
\def\CStalfp{{\mathop{\bf C^{\bs\iy}Sta^{lfp}}}}
\def\CStalf{{\mathop{\bf C^{\bs\iy}Sta^{lf}}}}
\def\CStalg{{\mathop{\bf C^{\bs\iy}Sta^{lg}}}}
\def\DMCSta{{\mathop{\bf DMC^{\bs\iy}Sta}}}
\def\DMCStalf{{\mathop{\bf DMC^{\bs\iy}Sta^{lf}}}}
\def\DMCStalfssc{{\mathop{\bf DMC^{\bs\iy}Sta^{lf}_{ssc}}}}
\def\DMCStalg{{\mathop{\bf DMC^{\bs\iy}Sta^{lg}}}}
\def\DMCStalfp{{\mathop{\bf DMC^{\bs\iy}Sta^{lfp}}}}
\def\hDMCStalfssc{{\mathop{\bf D\hat{M}C^{\bs\iy}Sta^{lf}_{ssc}}}}
\def\Presta{{\mathop{\bf Presta}\nolimits}}
\def\Sta{{\mathop{\bf Sta}\nolimits}}
\def\GSta{{\mathop{\bf GSta}\nolimits}}
\def\Man{{\mathop{\bf Man}}}
\def\Manb{{\mathop{\bf Man^b}}}
\def\Manc{{\mathop{\bf Man^c}}}
\def\dMan{{\mathop{\bf dMan}}}
\def\bdMan{{\mathop{\bf d\bar{M}an}}}
\def\dManb{{\mathop{\bf dMan^b}}}
\def\dManc{{\mathop{\bf dMan^c}}}
\def\dOrb{{\mathop{\bf dOrb}}}
\def\bdOrb{{\mathop{\bf d\bar{O}rb}}}
\def\dOrbb{{\mathop{\bf dOrb^b}}}
\def\dOrbc{{\mathop{\bf dOrb^c}}}
\def\dSpa{{\mathop{\bf dSpa}}}
\def\dSta{{\mathop{\bf dSta}}}
\def\dSpab{{\mathop{\bf dSpa^b}}}
\def\dStab{{\mathop{\bf dSta^b}}}
\def\dSpac{{\mathop{\bf dSpa^c}}}
\def\dStac{{\mathop{\bf dSta^c}}}
\def\hdSpa{{\mathop{\bf d\hat{S}pa}}}
\def\bdSpa{{\mathop{\bf d\bar{S}pa}}}
\def\bdSta{{\mathop{\bf d\bar{S}ta}}}
\def\hMan{{\mathop{\bf \hat{M}an}}}
\def\bMan{{\mathop{\bf \bar{M}an}}}
\def\bManb{{\mathop{\bf \bar{M}an^b}}}
\def\bManc{{\mathop{\bf \bar{M}an^c}}}
\def\Orb{{\mathop{\bf Orb}}}
\def\hOrb{{\mathop{\bf \hat{O}rb}}}
\def\bOrb{{\mathop{\bf \bar{O}rb}}}
\def\bOrbb{{\mathop{\bf \bar{O}rb^b}}}
\def\bOrbc{{\mathop{\bf \bar{O}rb^c}}}
\def\Orbb{{\mathop{\bf Orb^b}}}
\def\Orbc{{\mathop{\bf Orb^c}}}
\def\Euc{{\mathop{\bf Euc}}}
\def\Sets{{\mathop{\bf Sets}}}
\def\SSets{{\mathop{\bf SSets}}}
\def\bo{{\rm bo}}
\def\hom{{\rm hom}}
\def\dbo{{\rm dbo}}
\def\orb{{\rm orb}}
\def\obo{{\rm obo}}
\def\dobo{{\rm dobo}}
\def\ul{\underline}
\def\bs{\boldsymbol}
\def\ge{\geqslant}
\def\le{\leqslant\nobreak}
\def\O{{\mathcal O}}
\def\CS{{\mathcal C}_S}
\def\DS{{\mathcal D}_S}
\def\ES{{\mathcal E}_S}
\def\FS{{\mathcal F}_S}
\def\IS{{\mathcal I}_S}
\def\NS{{\mathcal N}_{{}\kern -.1em\sst\upS}^{\sst\uS}}
\def\OS{{{\mathcal O}_{\kern -.1em S}}}
\def\OSp{{{\mathcal O}^{\smash{\prime}}_{\kern -.1em S}}}
\def\CU{{\mathcal C}_U}
\def\DU{{\mathcal D}_U}
\def\EU{{\mathcal E}_U}
\def\FU{{\mathcal F}_U}
\def\IU{{\mathcal I}_U}
\def\NU{{\mathcal N}_{{}\kern -.1em\sst\upU}^{\sst\uU}}
\def\OU{{{\mathcal O}_{\kern -.1em U}}}
\def\OUp{{{\mathcal O}^{\smash{\prime}}_{\kern -.1em U}}}
\def\CV{{\mathcal C}_V}
\def\DV{{\mathcal D}_V}
\def\EV{{\mathcal E}_V}
\def\FV{{\mathcal F}_V}
\def\IV{{\mathcal I}_V}
\def\NV{{\mathcal N}_{{}\kern -.1em\sst\upV}^{\sst\uV}}
\def\OV{{{\mathcal O}_{\kern -.1em V}}}
\def\OVp{{{\mathcal O}^{\smash{\prime}}_{\kern -.1em V}}}
\def\CW{{\mathcal C}_W}
\def\DW{{\mathcal D}_W}
\def\EW{{\mathcal E}_W}
\def\FW{{\mathcal F}_W}
\def\IW{{\mathcal I}_W}
\def\NW{{\mathcal N}_{{}\kern -.1em\sst\upW}^{\sst\uW}}
\def\NpW{{\mathcal N}_{{}\kern -.1em\sst\uptW}^{\sst\upW}}
\def\OW{{{\mathcal O}_{\kern -.1em W}}}
\def\OWp{{{\mathcal O}^{\smash{\prime}}_{\kern -.1em W}}}
\def\CX{{\mathcal C}_X}
\def\DX{{\mathcal D}_X}
\def\EX{{\mathcal E}_X}
\def\FX{{\mathcal F}_X}
\def\IX{{\mathcal I}_X}
\def\NX{{\mathcal N}_{{}\kern -.1em\sst\upX}^{\sst\uX}}
\def\NpX{{\mathcal N}_{{}\kern -.1em\sst\uptX}^{\sst\upX}}
\def\OX{{{\mathcal O}_{\kern -.1em X}}}
\def\OXp{{{\mathcal O}^{\smash{\prime}}_{\kern -.1em X}}}
\def\CY{{\mathcal C}_Y}
\def\DY{{\mathcal D}_Y}
\def\EY{{\mathcal E}_Y}
\def\FY{{\mathcal F}_Y}
\def\IY{{\mathcal I}_Y}
\def\NY{{\mathcal N}_{{}\kern -.1em\sst\upY}^{\sst\uY}}
\def\OY{{{\mathcal O}_{\kern -.1em Y}}}
\def\OYp{{{\mathcal O}^{\smash{\prime}}_{\kern -.1em Y}}}
\def\CZ{{\mathcal C}_Z}
\def\DZ{{\mathcal D}_Z}
\def\EZ{{\mathcal E}_Z}
\def\FZ{{\mathcal F}_Z}
\def\IZ{{\mathcal I}_Z}
\def\NZ{{\mathcal N}_{{}\kern -.1em\sst\upZ}^{\sst\uZ}}
\def\OZ{{{\mathcal O}_{\kern -.1em Z}}}
\def\OZp{{{\mathcal O}^{\smash{\prime}}_{\kern -.1em Z}}}
\def\CcS{{\mathcal C}_\cS}
\def\DcS{{\mathcal D}_\cS}
\def\EcS{{\mathcal E}_\cS}
\def\FcS{{\mathcal F}_\cS}
\def\IcS{{\mathcal I}_\cS}
\def\OcS{{{\mathcal O}_{\kern -.1em\cS}}}
\def\OcSp{{{\mathcal O}^{\smash{\prime}}_{\kern -.1em\cS}}}
\def\CcW{{\mathcal C}_\cW}
\def\DcW{{\mathcal D}_\cW}
\def\EcW{{\mathcal E}_\cW}
\def\FcW{{\mathcal F}_\cW}
\def\IcW{{\mathcal I}_\cW}
\def\OcW{{{\mathcal O}_{\kern -.1em\cW}}}
\def\OcWp{{{\mathcal O}^{\smash{\prime}}_{\kern -.1em\cW}}}
\def\CcX{{\mathcal C}_\cX}
\def\DcX{{\mathcal D}_\cX}
\def\EcX{{\mathcal E}_\cX}
\def\FcX{{\mathcal F}_\cX}
\def\IcX{{\mathcal I}_\cX}
\def\OcX{{{\mathcal O}_{\kern -.1em\cX}}}
\def\OcXp{{{\mathcal O}^{\smash{\prime}}_{\kern -.1em\cX}}}
\def\CcY{{\mathcal C}_\cY}
\def\DcY{{\mathcal D}_\cY}
\def\EcY{{\mathcal E}_\cY}
\def\FcY{{\mathcal F}_\cY}
\def\IcY{{\mathcal I}_\cY}
\def\OcY{{{\mathcal O}_{\kern -.1em\cY}}}
\def\OcYp{{{\mathcal O}^{\smash{\prime}}_{\kern -.1em\cY}}}
\def\CcZ{{\mathcal C}_\cZ}
\def\DcZ{{\mathcal D}_\cZ}
\def\EcZ{{\mathcal E}_\cZ}
\def\FcZ{{\mathcal F}_\cZ}
\def\IcZ{{\mathcal I}_\cZ}
\def\OcZ{{{\mathcal O}_{\kern -.1em\cZ}}}
\def\OcZp{{{\mathcal O}^{\smash{\prime}}_{\kern -.1em\cZ}}}
\def\im{\imath}
\def\jm{\jmath}
\def\bR{{\mathbin{\pmb{\mathbb R}}}}
\def\R{{\mathbin{\mathbb R}}}
\def\Z{{\mathbin{\mathbb Z}}}
\def\Q{{\mathbin{\mathbb Q}}}
\def\N{{\mathbin{\mathbb N}}}
\def\C{{\mathbin{\mathbb C}}}
\def\CP{{\mathbin{\mathbb{CP}}}}
\def\RP{{\mathbin{\mathbb{RP}}}}
\def\fC{{\mathbin{\mathfrak C}\kern.05em}}
\def\fD{{\mathbin{\mathfrak D}}}
\def\fE{{\mathbin{\mathfrak E}}}
\def\fF{{\mathbin{\mathfrak F}}}
\def\fn{{\mathbin{\mathfrak n}}}
\def\cA{{\mathbin{\cal A}}}
\def\cB{{\mathbin{\cal B}}}
\def\cC{{\mathbin{\cal C}}}
\def\cD{{\mathbin{\cal D}}}
\def\cE{{\mathbin{\cal E}}}
\def\cF{{\mathbin{\cal F}}}
\def\cG{{\mathbin{\cal G}}}
\def\cH{{\mathbin{\cal H}}}
\def\cI{{\mathbin{\cal I}}}
\def\cJ{{\mathbin{\cal J}}}
\def\cK{{\mathbin{\cal K}}}
\def\cL{{\mathbin{\cal L}}}
\def\cM{{\mathbin{\cal M}}}
\def\cN{{\mathbin{\cal N}}}
\def\cP{{\mathbin{\cal P}}}
\def\cQ{{\mathbin{\cal Q}}}
\def\cR{{\mathbin{\cal R}}}
\def\cS{{\mathbin{\cal S}\kern -0.1em}}
\def\cT{{\mathbin{\cal T}\kern -0.1em}}
\def\cW{{\mathbin{\cal W}}}
\def\oM{{\mathbin{\smash{\,\,\overline{\!\!\mathcal M\!}\,}}}}
\def\bcM{{\mathbin{\bs{\cal M}}}}
\def\fCmod{{\mathbin{{\mathfrak C}\text{\rm -mod}}}}
\def\fCmodco{{\mathbin{{\mathfrak C}\text{\rm -mod}^{\rm co}}}}
\def\fCmodfp{{\mathbin{{\mathfrak C}\text{\rm -mod}^{\rm fp}}}}
\def\fDmod{{\mathbin{{\mathfrak D}\text{\rm -mod}}}}
\def\fDmodco{{\mathbin{{\mathfrak D}\text{\rm -mod}^{\rm co}}}}
\def\fDmodfp{{\mathbin{{\mathfrak D}\text{\rm -mod}^{\rm fp}}}}
\def\OXmod{{\mathbin{\O_X\text{\rm -mod}}}}
\def\OYmod{{\mathbin{\O_Y\text{\rm -mod}}}}
\def\OZmod{{\mathbin{\O_Z\text{\rm -mod}}}}
\def\OcXmod{{\mathbin{{\O\kern -0.1em}_{\cal X}\text{\rm -mod}}}}
\def\OcYmod{{\mathbin{\O_{\cal Y}\text{\rm -mod}}}}
\def\OcZmod{{\mathbin{{\O\kern -0.1em}_{\cal Z}\text{\rm -mod}}}}
\def\m{{\mathfrak m}}
\def\ua{{\underline{a}}{}}
\def\ub{{\underline{b\kern -0.1em}\kern 0.1em}{}}
\def\uc{{\underline{c}}{}}
\def\ud{{\underline{d}}{}}
\def\ue{{\underline{e}}{}}
\def\uf{{\underline{f\!}\,}{}}
\def\ug{{\underline{g\!}\,}{}}
\def\uh{{\underline{h\!}\,}{}}
\def\ui{{\underline{i\kern -0.07em}\kern 0.07em}{}}
\def\uim{{\underline{\imath\kern -0.07em}\kern 0.07em}{}}
\def\uj{{\underline{j\kern -0.1em}\kern 0.1em}{}}
\def\uk{{\underline{k\kern -0.1em}\kern 0.1em}{}}
\def\um{{\underline{m\kern -0.1em}\kern 0.1em}{}}
\def\us{{\underline{s\kern -0.15em}\kern 0.15em}{}}
\def\ut{{\underline{t\kern -0.1em}\kern 0.1em}{}}
\def\uu{{\underline{u\kern -0.1em}\kern 0.1em}{}}
\def\uv{{\underline{v\!}\,}{}}
\def\uw{{\underline{w\!}\,}{}}
\def\ux{{\underline{x\!}\,}{}}
\def\uy{{\underline{y\!}\,}{}}
\def\uz{{\underline{z\!}\,}{}}
\def\uA{{{\underline{A\!}\,}}{}}
\def\uB{{{\underline{B\!}\,}}{}}
\def\uC{{{\underline{C\!}\,}}{}}
\def\uD{{{\underline{D\!}\,}}{}}
\def\uE{{{\underline{E\!}\,}}{}}
\def\uF{{{\underline{F\!}\,}}{}}
\def\uG{{{\underline{G\!}\,}}{}}
\def\uH{{{\underline{H\!}\,}}{}}
\def\uI{{{\underline{I\kern -0.15em}\kern 0.1em}}{}}
\def\uK{{{\underline{K\!}\,}}{}}
\def\uP{{{\underline{P\!}\,}}{}}
\def\uQ{{{\underline{Q\!}\,}}{}}
\def\uR{{{\underline{R\kern -0.1em}\kern 0.1em}}{}}
\def\uS{{{\underline{S\!}\,}}{}}
\def\uT{{{\underline{T\!}\,}}{}}
\def\uU{{{\underline{U\kern -0.25em}\kern 0.2em}{}}}
\def\utU{{{\underline{\ti U\kern -0.25em}\kern 0.2em}}{}}
\def\uV{{{\underline{V\kern -0.25em}\kern 0.2em}}{}}
\def\utV{{{\underline{\ti V\kern -0.25em}\kern 0.2em}}{}}
\def\utW{{{\underline{\ti W\!\!}\,\,}}{}}
\def\uW{{{\underline{W\!\!}\,\,}}{}}
\def\uX{{{\underline{X\!}\,}}{}}
\def\uY{{{\underline{Y\!\!}\,\,}}{}}
\def\uZ{{{\underline{Z\!}\,}}{}}
\def\tiuU{{{\underline{\ti U\kern -0.25em}\kern 0.2em}}{}}
\def\cU{{\cal U}}
\def\cV{{\cal V}}
\def\cW{{\cal W}}
\def\cX{{\cal X}}
\def\cY{{\cal Y}}
\def\cZ{{\cal Z}}
\def\bC{{\bs C}}
\def\bD{{\bs D}}
\def\bE{{\bs E}}
\def\bF{{\bs F}}
\def\bG{{\bs G}}
\def\bH{{\bs H}}
\def\bN{{\bs N}}
\def\bS{{\bs S}}
\def\bT{{\bs T}}
\def\bU{{\bs U}}
\def\bV{{\bs V}}
\def\bW{{\bs W}\kern -0.1em}
\def\bX{{\bs X}}
\def\bY{{\bs Y}\kern -0.1em}
\def\bZ{{\bs Z}}
\def\bpU{{\bs{\pd U}}}
\def\bpV{{\bs{\pd V}}}
\def\bpW{{\bs{\pd W}}}
\def\bpX{{\bs{\pd X}}}
\def\bpY{{\bs{\pd Y}}}
\def\bpZ{{\bs{\pd Z}}}
\def\rb{{\bs{\rm b}}}
\def\rc{{\bs{\rm c}}}
\def\rd{{\bs{\rm d}}}
\def\re{{\bs{\rm e}}}
\def\rf{{\bs{\rm f}}}
\def\rg{{\bs{\rm g}}}
\def\rh{{\bs{\rm h}}}
\def\ri{{\bs{\rm i}}}
\def\rj{{\bs{\rm j}}}
\def\rE{{\bs{\rm E}}}
\def\rS{{\bs{\rm S}}}
\def\rT{{\bs{\rm T}}}
\def\rU{{\bs{\rm U}}}
\def\rV{{\bs{\rm V}}}
\def\rW{{\bs{\rm W}}}
\def\rX{{\bs{\rm X}}}
\def\rY{{\bs{\rm Y}}}
\def\rZ{{\bs{\rm Z}}}
\def\ed{{\bs{\sf d}}}
\def\ee{{\bs{\sf e}}}
\def\ef{{\bs{\sf f}}}
\def\eg{{\bs{\sf g}}}
\def\eh{{\bs{\sf h}}}
\def\ei{{\bs{\sf i}}}
\def\eU{{\bs{\EuScript U}}}
\def\eV{{\bs{\EuScript V}}}
\def\eW{{\bs{\EuScript W}}}
\def\eX{{\bs{\EuScript X}}}
\def\eY{{\bs{\EuScript Y}}}
\def\eZ{{\bs{\EuScript Z}}}
\def\bcE{{\bs{\cal E}}}
\def\bcU{{\bs{\cal U}}}
\def\bcV{{\bs{\cal V}}}
\def\bcW{{\bs{\cal W}}}
\def\bcX{{\bs{\cal X}}}
\def\bcY{{\bs{\cal Y}}}
\def\bcZ{{\bs{\cal Z}}}
\def\cpU{{\pd{\cal U}}}
\def\cpV{{\pd{\cal V}}}
\def\cpW{{\pd{\cal W}}}
\def\cpX{{\pd{\cal X}}}
\def\cpY{{\pd{\cal Y}}}
\def\cpZ{{\pd{\cal Z}}}
\def\bcpZ{{\bs{\pd{\cal Z}}}}
\def\bcpU{{\bs{\pd{\cal U}}}}
\def\bcpV{{\bs{\pd{\cal V}}}}
\def\bcpW{{\bs{\pd{\cal W}}}}
\def\bcpX{{\bs{\pd{\cal X}}}}
\def\bcpY{{\bs{\pd{\cal Y}}}}
\def\bcpZ{{\bs{\pd{\cal Z}}}}
\def\upS{{\underline{\pd S\!}\,}}
\def\upT{{\underline{\pd T\!}\,}}
\def\upU{{\underline{\pd U\!}\,}}
\def\upV{{\underline{\pd V\!\!}\,\,}}
\def\upW{{\underline{\pd W\!\!}\kern .25em}}
\def\upX{{\underline{\pd X\!}\,}}
\def\upY{{\underline{\pd Y\!\!}\,\,}}
\def\upZ{{\underline{\pd Z\!}\,}}
\def\uptW{{\underline{\pd^2W\!\!}\,\,}}
\def\uptX{{\underline{\pd^2X\!}\,}}
\def\uptY{{\underline{\pd^2Y\!\!}\,\,}}
\def\uptZ{{\underline{\pd^2Z\!}\,}}
\def\uXi{{\underline{\Xi\!}\,}}
\def\uxi{{\underline{\xi\!}\,}{}}
\def\uchi{{\underline{\chi\!}\,}{}}
\def\upi{{\underline{\pi\!}\,}}
\def\uid{{\underline{\id\kern -0.1em}\kern 0.1em}}
\def\urho{{\underline{\rho\!}\,}}
\def\usi{{\underline{\si\!}\,}}
\def\uphi{{\underline{\phi\!}\,}}
\def\umu{{\underline{\smash{\mu}}\kern 0.1em}}
\def\uzi{{\underline{[0,\infty)\!}\,}}
\def\al{\alpha}
\def\be{\beta}
\def\ga{\gamma}
\def\de{\delta}
\def\io{\iota}
\def\ep{\epsilon}
\def\la{\lambda}
\def\ka{\kappa}
\def\th{\theta}
\def\ze{\zeta}
\def\up{\upsilon}
\def\vp{\varphi}
\def\si{\sigma}
\def\om{\omega}
\def\De{\Delta}
\def\La{\Lambda}
\def\Om{\Omega}
\def\Up{\Upsilon}
\def\Ga{\Gamma}
\def\Si{\Sigma}
\def\Th{\Theta}
\def\pd{\partial}
\def\ts{\textstyle}
\def\st{\scriptstyle}
\def\sst{\scriptscriptstyle}
\def\w{\wedge}
\def\sm{\setminus}
\def\lt{\ltimes}
\def\bu{\bullet}
\def\sh{\sharp}
\def\op{\oplus}
\def\od{\odot}
\def\op{\oplus}
\def\ot{\otimes}
\def\hot{{\mathop{\kern .1em\hat\otimes\kern .1em}\nolimits}}
\def\hbt{{\mathop{\kern .1em\hat\boxtimes\kern .1em}\nolimits}}
\def\ov{\overline}
\def\bigop{\bigoplus}
\def\bigot{\bigotimes}
\def\iy{\infty}
\def\es{\emptyset}
\def\ra{\rightarrow}
\def\rra{\rightrightarrows}
\def\Ra{\Rightarrow}
\def\ab{\allowbreak}
\def\longra{\longrightarrow}
\def\hookra{\hookrightarrow}
\def\dashra{\dashrightarrow}
\def\t{\times}
\def\ci{\circ}
\def\ti{\tilde}
\def\d{{\rm d}}
\def\ha{{\ts\frac{1}{2}}}
\def\md#1{\vert #1 \vert}
\def\bmd#1{\big\vert #1 \big\vert}
\def\ms#1{\vert #1 \vert^2}
\def\an#1{\langle #1 \rangle}
\def\ban#1{\bigl\langle #1 \bigr\rangle}
%%%%%%%%%%%%%%%%%%%%%%%%%%%%%%%%%%%%%%%%%%%%%%%%%%%%%%%%%%%%%%%%%%%%%%%%
%%%%%%%%%%%%%%%%%%%%%%%    Text of paper    %%%%%%%%%%%%%%%%%%%%%%%%%%%%
%%%%%%%%%%%%%%%%%%%%%%%%%%%%%%%%%%%%%%%%%%%%%%%%%%%%%%%%%%%%%%%%%%%%%%%%
\title{An introduction to d-manifolds \\ and derived differential geometry}
\author{Dominic Joyce}
\date{}
\maketitle

\begin{abstract} This is a survey of the author's book
\cite{Joyc5}. We introduce a 2-category $\dMan$ of {\it
d-manifolds}, new geometric objects which are `derived' smooth
manifolds, in the sense of the `derived algebraic geometry' of
To\"en and Lurie. Manifolds $\Man$ embed in $\dMan$ as a full
(2-)subcategory. There are also 2-categories $\dManb,\dManc$ of {\it
d-manifolds with boundary\/} and {\it with corners}, and orbifold
versions $\dOrb,\dOrbb,\dOrbc$, {\it d-orbifolds}.

Much of differential geometry extends very nicely to d-manifolds ---
immersions, submersions, submanifolds, transverse fibre products,
orientations, etc. Compact oriented d-manifolds have virtual
classes.

Many areas of symplectic geometry involve `counting' moduli spaces
$\oM_{g,m}(J,\be)$ of $J$-holomorphic curves to define invariants,
Floer homology theories, etc. Such $\oM_{g,m}(J,\be)$ are given the
structure of {\it Kuranishi spaces\/} in the work of Fukaya, Oh,
Ohta and Ono \cite{FOOO}, but there are problems with the theory.
The author believes the `correct' definition of Kuranishi spaces is
that they are d-orbifolds with corners. D-manifolds and d-orbifolds
will have applications in symplectic geometry, and elsewhere.

For brevity, this survey focusses on d-manifolds without boundary. A
longer and more detailed summary of the book is given
in~\cite{Joyc6}.
\end{abstract}

\baselineskip 11.9pt plus .1pt

\setcounter{tocdepth}{2} \tableofcontents

\section{Introduction}
\label{sd1}

This is a survey of \cite{Joyc5}, and describes the author's new
theory of `derived differential geometry'. The objects in this
theory are {\it d-manifolds}, `derived' versions of smooth
manifolds, which form a (strict) 2-category $\dMan$. There are also
2-categories of {\it d-manifolds with boundary\/} $\dManb$ and {\it
d-manifolds with corners\/} $\dManc$, and orbifold versions of all
these, {\it d-orbifolds\/} $\dOrb,\dOrbb,\dOrbc$.

Here `derived' is intended in the sense of {\it derived algebraic
geometry}. The original motivating idea for derived algebraic
geometry, as in Kontsevich \cite{Kont} for instance, was that
certain moduli schemes $\cM$ appearing in enumerative invariant
problems may be highly singular as schemes. However, it may be
natural to realize $\cM$ as a truncation of some `derived' moduli
space $\bcM$, a new kind of geometric object living in a higher
category. The geometric structure on $\bcM$ should encode the full
deformation theory of the moduli problem, the obstructions as well
as the deformations. It was hoped that $\bcM$ would be `smooth', and
so in some sense simpler than its truncation~$\cM$.

Early work in derived algebraic geometry focussed on {\it
dg-schemes}, as in Ciocan-Fontanine and Kapranov \cite{CiKa}. These
have largely been replaced by the {\it derived stacks\/} of To\"en
and Vezzosi \cite{Toen,ToVe}, and the {\it structured spaces\/} of
Lurie \cite{Luri1,Luri2}. {\it Derived differential geometry\/} aims
to generalize these ideas to differential geometry and smooth
manifolds. A brief note about it can be found in Lurie \cite[\S
4.5]{Luri2}; the ideas are worked out in detail by Lurie's student
David Spivak \cite{Spiv}, who defines an $\iy$-category (simplicial
category) of {\it derived manifolds}.

The author came to these questions from a different direction,
symplectic geometry. Many important areas in symplectic geometry
involve forming moduli spaces $\oM_{g,m}(J,\be)$ of $J$-holomorphic
curves in some symplectic manifold $(M,\om)$, possibly with boundary
in a Lagrangian $L$, and then `counting' these moduli spaces to get
`invariants' with interesting properties. Such areas include
Gromov--Witten invariants (open and closed), Lagrangian Floer
cohomology, Symplectic Field Theory, contact homology, and Fukaya
categories.

To do this `counting', one needs to put a suitable geometric
structure on $\oM_{g,m}(J,\be)$ --- something like the `derived'
moduli spaces $\bcM$ above --- and use this to define a `virtual
class' or `virtual chain' in $\Z,\Q$ or some homology theory. There
is no general agreement on what geometric structure to use ---
compared to the elegance of algebraic geometry, this area is
something of a mess. Two rival theories for geometric structures to
put on moduli spaces $\oM_{g,m}(J,\be)$ are the {\it Kuranishi
spaces\/} of Fukaya, Oh, Ohta and Ono \cite{FuOn,FOOO} and the {\it
polyfolds\/} of Hofer, Wysocki and Zehnder~\cite{HWZ1,HWZ2,HWZ3}.

The theory of Kuranishi spaces in \cite{FuOn,FOOO} does not go far
-- they define Kuranishi spaces, and construct virtual cycles upon
them, but they do not define morphisms between Kuranishi spaces, for
instance. The author tried to study Kuranishi spaces as geometric
spaces in their own right, but ran into problems, and became
convinced that a new definition of Kuranishi space was needed. Upon
reading Spivak's theory of derived manifolds \cite{Spiv}, it became
clear that some form of `derived differential geometry' was
required: {\it Kuranishi spaces in the sense of\/ {\rm\cite[\S
A]{FOOO}} ought to be defined to be `derived orbifolds with
corners'}.

The author tried to read Lurie \cite{Luri1,Luri2} and Spivak
\cite{Spiv} with a view to applications to Kuranishi spaces and
symplectic geometry, but ran into problems of a different kind: the
framework of \cite{Luri1,Luri2,Spiv} is formidably long, complex and
abstract, and proved too difficult for a humble trainee symplectic
geometer to understand, or use as a tool. So the author looked for a
way to simplify the theory, while retaining the information needed
for applications to symplectic geometry. The theory of d-manifolds
and d-orbifolds of \cite{Joyc5} is the result.

The essence of our simplification is this. Consider a `derived'
moduli space $\bcM$ of some objects $E$, e.g. vector bundles on some
$\C$-scheme $X$. One expects $\bcM$ to have a `cotangent complex'
${\mathbb L}_\bcM$, a complex in some derived category with
cohomology $h^i({\mathbb L}_\bcM)\vert_E\cong\Ext^{1-i}(E,E)^*$ for
$i\in\Z$. In general, ${\mathbb L}_\bcM$ can have nontrivial
cohomology in many negative degrees, and because of this such
objects $\bcM$ must form an $\iy$-category to properly describe
their geometry.

However, the moduli spaces relevant to enumerative invariant
problems are of a restricted kind: one considers only $\bcM$ such
that ${\mathbb L}_\bcM$ has nontrivial cohomology only in degrees
$-1,0$, where $h^0({\mathbb L}_\bcM)$ encodes the (dual of the)
deformations $\Ext^1(E,E)^*$, and $h^{-1}({\mathbb L}_\bcM)$ the
(dual of the) obstructions $\Ext^2(E,E)^*$. As in To\"en \cite[\S
4.4.3]{Toen}, such derived spaces are called {\it quasi-smooth}, and
this is a necessary condition on $\bcM$ for the construction of a
virtual fundamental class.

Our construction of d-manifolds replaces complexes in a derived
category $D^b\coh(\cM)$ with a 2-category of complexes in degrees
$-1,0$ only. For general $\bcM$ this loses a lot of information, but
for quasi-smooth $\bcM$, since ${\mathbb L}_\bcM$ is concentrated in
degrees $-1,0$, the important information is retained. In the
language of To\"en and Vezzosi \cite{Toen,ToVe}, this corresponds to
working with a subclass of derived schemes whose dg-algebras are of
a special kind: they are 2-step supercommutative dg-algebras
$A^{-1}\,{\buildrel\d\over \longra}\,A^0$ such that $\d(A^{-1})\cdot
A^{-1}=0$. Then $\d(A^{-1})$ is a square zero ideal in $A^0$, and
$A^{-1}$ is a module over~$H^0\bigl(A^{-1}\,{\buildrel
\d\over\longra}\,A^0\bigr)$.

The set up of \cite{Joyc5} is also long and complicated. But mostly
this complexity comes from other sources: working over
$C^\iy$-rings, and including manifolds with boundary, manifolds with
corners, and orbifolds. The 2-category style `derived geometry' of
\cite{Joyc5} really is far simpler than those
of~\cite{Luri1,Luri2,Spiv,Toen,ToVe}.

Following Spivak \cite{Spiv}, in order to be able to use the tools
of algebraic geometry --- schemes, stacks, quasicoherent sheaves ---
in differential geometry, our d-manifolds are built on the notions
of $C^\iy$-{\it ring\/} and $C^\iy$-{\it scheme\/} that were
invented in synthetic differential geometry, and developed further
by the author in \cite{Joyc3,Joyc4}. We survey the $C^\iy$-algebraic
geometry we need in \S\ref{sd2}. Section \ref{sd3} discusses the
2-category of {\it d-spaces\/} $\dSpa$, which are `derived'
$C^\iy$-schemes, and \S\ref{sd4} describes the 2-category of
d-manifolds $\dMan$, and our theory of `derived differential
geometry'. Appendix \ref{sdA} explains the basics of 2-categories.

For brevity, and to get the main ideas across as simply as possible,
this survey will concentrate on d-manifolds without boundary, apart
from a short section on d-manifolds with corners and d-orbifolds in
\S\ref{sd49}. A longer summary of the book, including much more
detail on d-manifolds with corners, d-orbifolds, and d-orbifolds
with corners, is given in~\cite{Joyc6}.
\medskip

\noindent{\it Acknowledgements.} I would like to thank Jacob Lurie
for helpful conversations.

\section{$C^\iy$-rings and $C^\iy$-schemes}
\label{sd2}

If $X$ is a manifold then the $\R$-algebra $C^\iy(X)$ of smooth
functions $c:X\ra\R$ is a $C^\iy$-{\it ring}. That is, for each
smooth function $f:\R^n\ra\R$ there is an $n$-fold operation
$\Phi_f:C^\iy(X)^n\ra C^\iy(X)$ acting by
$\Phi_f:c_1,\ldots,c_n\mapsto f(c_1,\ldots,c_n)$, and these
operations $\Phi_f$ satisfy many natural identities. Thus,
$C^\iy(X)$ actually has a far richer algebraic structure than the
obvious $\R$-algebra structure.

In \cite{Joyc3} (surveyed in \cite{Joyc4}) the author set out a
version of algebraic geometry in which rings or algebras are
replaced by $C^\iy$-rings, focussing on $C^\iy$-{\it schemes}, a
category of geometric objects which generalize manifolds, and whose
morphisms generalize smooth maps, {\it quasicoherent\/} and {\it
coherent sheaves\/} on $C^\iy$-schemes, and $C^\iy$-{\it stacks}, in
particular {\it Deligne--Mumford\/ $C^\iy$-stacks}, a 2-category of
geometric objects which generalize orbifolds. Much of the material
on $C^\iy$-schemes was already known in synthetic differential
geometry, see for instance Dubuc \cite{Dubu} and Moerdijk and
Reyes~\cite{MoRe}.

\subsection{$C^\iy$-rings}
\label{sd21}

\begin{dfn} A $C^\iy$-{\it ring\/} is a set $\fC$ together with
operations $\Phi_f:\fC^n\ra\fC$ for all $n\ge 0$ and smooth maps
$f:\R^n\ra\R$, where by convention when $n=0$ we define $\fC^0$ to
be the single point $\{\es\}$. These operations must satisfy the
following relations: suppose $m,n\ge 0$, and $f_i:\R^n\ra\R$ for
$i=1,\ldots,m$ and $g:\R^m\ra\R$ are smooth functions. Define a
smooth function $h:\R^n\ra\R$ by
\begin{equation*}
h(x_1,\ldots,x_n)=g\bigl(f_1(x_1,\ldots,x_n),\ldots,f_m(x_1
\ldots,x_n)\bigr),
\end{equation*}
for all $(x_1,\ldots,x_n)\in\R^n$. Then for all
$(c_1,\ldots,c_n)\in\fC^n$ we have
\begin{equation*}
\Phi_h(c_1,\ldots,c_n)=\Phi_g\bigl(\Phi_{f_1}(c_1,\ldots,c_n),
\ldots,\Phi_{f_m}(c_1,\ldots,c_n)\bigr).
\end{equation*}
We also require that for all $1\le j\le n$, defining
$\pi_j:\R^n\ra\R$ by $\pi_j:(x_1,\ldots,x_n)\mapsto x_j$, we have
$\Phi_{\pi_j}(c_1,\ldots,c_n)=c_j$ for
all~$(c_1,\ldots,c_n)\in\fC^n$.

Usually we refer to $\fC$ as the $C^\iy$-ring, leaving the
operations $\Phi_f$ implicit.

A {\it morphism\/} between $C^\iy$-rings $\bigl(\fC,(\Phi_f)_{
f:\R^n\ra\R\,\,C^\iy}\bigr)$, $\bigl({\mathfrak
D},(\Psi_f)_{f:\R^n\ra\R\,\,C^\iy}\bigr)$ is a map
$\phi:\fC\ra{\mathfrak D}$ such that $\Psi_f\bigl(\phi
(c_1),\ldots,\phi(c_n)\bigr)=\phi\ci\Phi_f(c_1,\ldots,c_n)$ for all
smooth $f:\R^n\ra\R$ and $c_1,\ldots,c_n\in\fC$. We will write
$\CRings$ for the category of $C^\iy$-rings.
\label{sd2def1}
\end{dfn}

Here is the motivating example:

\begin{ex} Let $X$ be a manifold. Write $C^\iy(X)$ for the set of
smooth functions $c:X\ra\R$. For $n\ge 0$ and $f:\R^n\ra\R$ smooth,
define $\Phi_f:C^\iy(X)^n\ra C^\iy(X)$ by
\e
\bigl(\Phi_f(c_1,\ldots,c_n)\bigr)(x)=f\bigl(c_1(x),\ldots,c_n(x)\bigr),
\label{sd2eq1}
\e
for all $c_1,\ldots,c_n\in C^\iy(X)$ and $x\in X$. It is easy to see
that $C^\iy(X)$ and the operations $\Phi_f$ form a~$C^\iy$-ring.

Now let $f:X\ra Y$ be a smooth map of manifolds. Then pullback
$f^*:C^\iy(Y)\ra C^\iy(X)$ mapping $f^*:c\mapsto c\ci f$ is a
morphism of $C^\iy$-rings. Furthermore (at least for $Y$ without
boundary), every $C^\iy$-ring morphism $\phi:C^\iy(Y)\ra C^\iy(X)$
is of the form $\phi=f^*$ for a unique smooth map~$f:X\ra Y$.

Write $\CRings^{\rm op}$ for the opposite category of $\CRings$,
with directions of morphisms reversed, and $\Man$ for the category
of manifolds without boundary. Then we have a full and faithful
functor $F_\Man^\CRings:\Man\ra\CRings^{\rm op}$ acting by
$F_\Man^\CRings(X)=C^\iy(X)$ on objects and $F_\Man^\CRings(f)=f^*$
on morphisms. This embeds $\Man$ as a full subcategory
of~$\CRings^{\rm op}$.
\label{sd2ex1}
\end{ex}

Note that $C^\iy$-rings are far more general than those coming from
manifolds. For example, if $X$ is any topological space we could
define a $C^\iy$-ring $C^0(X)$ to be the set of {\it continuous\/}
$c:X\ra\R$, with operations $\Phi_f$ defined as in \eq{sd2eq1}. For
$X$ a manifold with $\dim X>0$, the $C^\iy$-rings $C^\iy(X)$ and
$C^0(X)$ are different.

\begin{dfn} Let $\fC$ be a $C^\iy$-ring. Then we may give $\fC$ the
structure of a {\it commutative\/ $\R$-algebra}. Define addition
`$+$' on $\fC$ by $c+c'=\Phi_f(c,c')$ for $c,c'\in\fC$, where
$f:\R^2\ra\R$ is $f(x,y)=x+y$. Define multiplication `$\,\cdot\,$'
on $\fC$ by $c\cdot c'=\Phi_g(c,c')$, where $g:\R^2\ra\R$ is
$f(x,y)=xy$. Define scalar multiplication by $\la\in\R$ by $\la
c=\Phi_{\la'}(c)$, where $\la':\R\ra\R$ is $\la'(x)=\la x$. Define
elements $0,1\in\fC$ by $0=\Phi_{0'}(\es)$ and $1=\Phi_{1'}(\es)$,
where $0':\R^0\ra\R$ and $1':\R^0\ra\R$ are the maps $0':\es\mapsto
0$ and $1':\es\mapsto 1$. One can show using the relations on the
$\Phi_f$ that the axioms of a commutative $\R$-algebra are
satisfied. In Example \ref{sd2ex1}, this yields the obvious
$\R$-algebra structure on the smooth functions~$c:X\ra\R$.

An {\it ideal\/} $I$ in $\fC$ is an ideal $I\subset\fC$ in $\fC$
regarded as a commutative $\R$-algebra. Then we make the quotient
$\fC/I$ into a $C^\iy$-ring as follows. If $f:\R^n\ra\R$ is smooth,
define $\Phi_f^I:(\fC/I)^n\ra\fC/I$ by
\begin{equation*}
\bigl(\Phi_f^I(c_1+I,\ldots,c_n+I)\bigr)(x)=f\bigl(c_1(x),\ldots,
c_n(x)\bigr)+I.
\end{equation*}
Using Hadamard's Lemma, one can show that this is independent of the
choice of representatives $c_1,\ldots,c_n$. Then
$\bigl(\fC/I,(\Phi_f^I)_{ f:\R^n\ra\R\,\,C^\iy}\bigr)$ is a
$C^\iy$-ring.

A $C^\iy$-ring $\fC$ is called {\it finitely generated\/} if there
exist $c_1,\ldots,c_n$ in $\fC$ which generate $\fC$ over all
$C^\iy$-operations. That is, for each $c\in\fC$ there exists smooth
$f:\R^n\ra\R$ with $c=\Phi_f(c_1,\ldots,c_n)$. Given such
$\fC,c_1,\ldots,c_n$, define $\phi:C^\iy(\R^n)\ra\fC$ by
$\phi(f)=\Phi_f(c_1,\ldots,c_n)$ for smooth $f:\R^n\ra\R$, where
$C^\iy(\R^n)$ is as in Example \ref{sd2ex1} with $X=\R^n$. Then
$\phi$ is a surjective morphism of $C^\iy$-rings, so $I=\Ker\phi$ is
an ideal in $C^\iy(\R^n)$, and $\fC\cong C^\iy(\R^n)/I$ as a
$C^\iy$-ring. Thus, $\fC$ is finitely generated if and only if
$\fC\cong C^\iy(\R^n)/I$ for some $n\ge 0$ and some ideal $I$
in~$C^\iy(\R^n)$.
\label{sd2def2}
\end{dfn}

\subsection{$C^\iy$-schemes}
\label{sd22}

Next we summarize material in \cite[\S 4]{Joyc3} on $C^\iy$-schemes.

\begin{dfn} A {\it $C^\iy$-ringed space\/} $\uX=(X,\O_X)$ is a
topological space $X$ with a sheaf $\O_X$ of $C^\iy$-rings on $X$.

A {\it morphism\/} $\uf=(f,f^\sh):(X,\O_X)\ra (Y,\O_Y)$ of $C^\iy$
ringed spaces is a continuous map $f:X\ra Y$ and a morphism
$f^\sh:f^{-1}(\OY)\ra\OX$ of sheaves of $C^\iy$-rings on $X$, where
$f^{-1}(\OY)$ is the inverse image sheaf. There is another way to
write the data $f^\sh$: since direct image of sheaves $f_*$ is right
adjoint to inverse image $f^{-1}$, there is a natural bijection
\e
\Hom_X\bigl(f^{-1}(\O_Y),\O_X\bigr)\cong\Hom_Y\bigl(\O_Y,f_*(\O_X)\bigr).
\label{sd2eq2}
\e
Write $f_\sh:\O_Y\ra f_*(\O_X)$ for the morphism of sheaves of
$C^\iy$-rings on $Y$ corresponding to $f^\sh$ under \eq{sd2eq2}, so
that
\e
f^\sh:f^{-1}(\O_Y)\longra\O_X\quad \leftrightsquigarrow\quad
f_\sh:\O_Y\longra f_*(\O_X).
\label{sd2eq3}
\e
Depending on the application, either $f^\sh$ or $f_\sh$ may be more
useful. We choose to regard $f^\sh$ as primary and write morphisms
as $\uf=(f,f^\sh)$ rather than $(f,f_\sh)$, because we find it
convenient in \cite{Joyc5} to work uniformly using pullbacks, rather
than mixing pullbacks and pushforwards.

Write $\CRS$ for the category of $C^\iy$-ringed spaces. As in
\cite[Th.~8]{Dubu} there is a {\it spectrum functor\/}
$\Spec:\CRings^{\rm op}\ra\CRS$, defined explicitly in
\cite[Def.~4.12]{Joyc3}. A $C^\iy$-ringed space $\uX$ is called an
{\it affine\/ $C^\iy$-scheme\/} if it is isomorphic in $\CRS$ to
$\Spec\fC$ for some $C^\iy$-ring $\fC$. A $C^\iy$-ringed space
$\uX=(X,\OX)$ is called a $C^\iy$-{\it scheme\/} if $X$ can be
covered by open sets $U\subseteq X$ such that $(U,\O_X\vert_U)$ is
an affine $C^\iy$-scheme. Write $\CSch$ for the full subcategory of
$C^\iy$-schemes in~$\CRS$.

A $C^\iy$-scheme $\uX=(X,\OX)$ is called {\it locally fair} if $X$
can be covered by open $U\subseteq X$ with
$(U,\O_X\vert_U)\cong\Spec\fC$ for some finitely generated
$C^\iy$-ring $\fC$. Roughly speaking this means that $\uX$ is
locally finite-dimensional. Write $\CSchlf$ for the full subcategory
of locally fair $C^\iy$-schemes in~$\CSch$.

We call a $C^\iy$-scheme $\uX$ {\it separated, second countable,
compact}, or {\it paracompact}, if the underlying topological space
$X$ is Hausdorff, second countable, compact, or paracompact,
respectively.
\label{sd2def3}
\end{dfn}

We define a $C^\iy$-scheme $\uX$ for each manifold~$X$.

\begin{ex} Let $X$ be a manifold. Define a $C^\iy$-ringed space
$\uX=(X,\O_X)$ to have topological space $X$ and $\O_X(U)=C^\iy(U)$
for each open $U\subseteq X$, where $C^\iy(U)$ is the $C^\iy$-ring
of smooth maps $c:U\ra\R$, and if $V\subseteq U\subseteq X$ are open
define $\rho_{UV}:C^\iy(U)\ra C^\iy(V)$ by $\rho_{UV}:c\mapsto
c\vert_V$. Then $\uX=(X,\O_X)$ is a local $C^\iy$-ringed space. It
is canonically isomorphic to $\Spec C^\iy(X)$, and so is an affine
$C^\iy$-scheme. It is locally fair.

Define a functor $F_\Man^\CSch:\Man\ra\CSchlf\subset\CSch$ by
$F_\Man^\CSch=\Spec\ci F_\Man^\CRings$. Then $F_\Man^\CSch$ is full
and faithful, and embeds $\Man$ as a full subcategory of $\CSch$.
\label{sd2ex2}
\end{ex}

By \cite[Cor.~4.21 \& Th.~4.33]{Joyc3} we have:

\begin{thm} Fibre products and all finite limits exist in
the category\/ $\CSch$. The subcategory $\CSchlf$ is closed under
fibre products and all finite limits in\/ $\CSch$. The functor\/
$F_\Man^\CSch$ takes transverse fibre products in $\Man$ to fibre
products in\/~$\CSch$.
\label{sd2thm1}
\end{thm}

The proof of the existence of fibre products in $\CSch$ follows that
for fibre products of schemes in Hartshorne \cite[Th.~II.3.3]{Hart},
together with the existence of $C^\iy$-scheme products $\uX\t\uY$ of
affine $C^\iy$-schemes $\uX,\uY$. The latter follows from the
existence of coproducts $\fC\hat\ot\fD$ in $\CRings$ of
$C^\iy$-rings $\fC,\fD$. Here $\fC\hat\ot\fD$ may be thought of as a
`completed tensor product' of $\fC,\fD$. The actual tensor product
$\fC\ot_\R\fD$ is naturally an $\R$-algebra but not a $C^\iy$-ring,
with an inclusion of $\R$-algebras $\fC\ot_\R\fD\hookra
\fC\hat\ot\fD$, but $\fC\hat\ot\fD$ is often much larger than
$\fC\ot_\R\fD$. For free $C^\iy$-rings we have~$C^\iy(\R^m)\hat\ot
C^\iy(\R^n)\cong C^\iy(\R^{m+n})$.

In \cite[Def.~4.34 \& Prop.~4.35]{Joyc3} we discuss {\it partitions
of unity\/} on $C^\iy$-schemes.

\begin{dfn} Let $\uX=(X,\O_X)$ be a $C^\iy$-scheme. Consider a
formal sum $\sum_{a\in A}c_a$, where $A$ is an indexing set and
$c_a\in\O_X(X)$ for $a\in A$. We say $\sum_{a\in A}c_a$ is a {\it
locally finite sum on\/} $\uX$ if $X$ can be covered by open
$U\subseteq X$ such that for all but finitely many $a\in A$ we have
$\rho_{XU}(c_a)=0$ in~$\O_X(U)$.

By the sheaf axioms for $\O_X$, if $\sum_{a\in A}c_a$ is a locally
finite sum there exists a unique $c\in\O_X(X)$ such that for all
open $U\subseteq X$ with $\rho_{XU}(c_a)=0$ in $\O_X(U)$ for all but
finitely many $a\in A$, we have $\rho_{XU}(c)=\sum_{a\in
A}\rho_{XU}(c_a)$ in $\O_X(U)$, where the sum makes sense as there
are only finitely many nonzero terms. We call $c$ the {\it limit\/}
of $\sum_{a\in A}c_a$, written~$\sum_{a\in A}c_a=c$.

Let $c\in\O_X(X)$. Suppose $V_i\subseteq X$ is open and
$\rho_{XV_i}(c)=0\in\O_X(V_i)$ for $i\in I$, and let
$V=\bigcup_{i\in I}V_i$. Then $V\subseteq X$ is open, and
$\rho_{XV}(c)=0\in\O_X(V)$ as $\O_X$ is a sheaf. Thus taking the
union of all open $V\subseteq X$ with $\rho_{XV}(c)=0$ gives a
unique maximal open set $V_c\subseteq X$ such that
$\rho_{XV_c}(c)=0\in\O_X(V_c)$. Define the {\it support\/} $\supp c$
of $c$ to be $X\sm V_c$, so that $\supp c$ is closed in $X$. If
$U\subseteq X$ is open, we say that $c$ {\it is supported in\/} $U$
if~$\supp c\subseteq U$.

Let $\{U_a:a\in A\}$ be an open cover of $X$. A {\it partition of
unity on\/ $\uX$ subordinate to\/} $\{U_a:a\in A\}$ is
$\{\eta_a:a\in A\}$ with $\eta_a\in\O_X(X)$ supported on $U_a$ for
$a\in A$, such that $\sum_{a\in A}\eta_a$ is a locally finite sum on
$\uX$ with~$\sum_{a\in A}\eta_a=1$.
\label{sd2def4}
\end{dfn}

\begin{prop} Suppose $\uX$ is a separated, paracompact, locally fair
$C^\iy$-scheme, and\/ $\{\uU_a:a\in A\}$ an open cover of\/ $\uX$.
Then there exists a partition of unity $\{\eta_a:a\in A\}$ on $\uX$
subordinate to\/~$\{\uU_a:a\in A\}$.
\label{sd2prop}
\end{prop}

Here are some differences between ordinary schemes and
$C^\iy$-schemes:

\begin{rem}{\bf(i)} If $A$ is a ring or algebra, then points of the
corresponding scheme $\Spec A$ are prime ideals in $A$. However, if
$\fC$ is a $C^\iy$-ring then (by definition) points of $\Spec\fC$
are maximal ideals in $\fC$ with residue field $\R$, or
equivalently, $\R$-algebra morphisms $x:\fC\ra\R$. This has the
effect that if $X$ is a manifold then points of $\Spec C^\iy(X)$ are
just points of~$X$.
\smallskip

\noindent{\bf(ii)} In conventional algebraic geometry, affine
schemes are a restrictive class. Central examples such as $\CP^n$
are not affine, and affine schemes are not closed under open
subsets, so that $\C^2$ is affine but $\C^2\sm\{0\}$ is not. In
contrast, affine $C^\iy$-schemes are already general enough for many
purposes. For example:
\begin{itemize}
\setlength{\itemsep}{0pt}
\setlength{\parsep}{0pt}
\item All manifolds are fair affine $C^\iy$-schemes.
\item Open $C^\iy$-subschemes of fair affine $C^\iy$-schemes are
fair and affine.
\item Separated, second countable, locally fair $C^\iy$-schemes
are affine.
\end{itemize}
Affine $C^\iy$-schemes are always separated (Hausdorff), so we need
general $C^\iy$-schemes to include non-Hausdorff behaviour.
\smallskip

\noindent{\bf(iii)} In conventional algebraic geometry the Zariski
topology is too coarse for many purposes, so one has to introduce
the \'etale topology. In $C^\iy$-algebraic geometry there is no need
for this, as affine $C^\iy$-schemes are Hausdorff.
\smallskip

\noindent{\bf(iv)} Even very basic $C^\iy$-rings such as
$C^\iy(\R^n)$ for $n>0$ are not noetherian as $\R$-algebras. So
$C^\iy$-schemes should be compared to non-noetherian schemes in
conventional algebraic geometry.
\smallskip

\noindent{\bf(v)} The existence of partitions of unity, as in
Proposition \ref{sd2prop}, makes some things easier in
$C^\iy$-algebraic geometry than in conventional algebraic geometry.
For example, geometric objects can often be `glued together' over
the subsets of an open cover using partitions of unity, and if $\cE$
is a quasicoherent sheaf on a separated, paracompact, locally fair
$C^\iy$-scheme $\uX$ then $H^i(\cE)=0$ for~$i>0$.
\label{sd2rem1}
\end{rem}

\subsection{Modules over $C^\iy$-rings, and cotangent modules}
\label{sd23}

In \cite[\S 5]{Joyc3} we discuss modules over $C^\iy$-rings.

\begin{dfn} Let $\fC$ be a $C^\iy$-ring. A $\fC$-{\it module} $M$
is a module over $\fC$ regarded as a commutative $\R$-algebra as in
Definition \ref{sd2def2}. $\fC$-modules form an abelian category,
which we write as $\fCmod$. For example, $\fC$ is a $\fC$-module,
and more generally $\fC\ot_\R V$ is a $\fC$-module for any real
vector space $V$. Let $\phi:\fC\ra\fD$ be a morphism of
$C^\iy$-rings. If $M$ is a $\fC$-module then $\phi_*(M)=M\ot_\fC\fD$
is a $\fD$-module. This induces a functor~$\phi_*:\fCmod\ra\fDmod$.
\label{sd2def5}
\end{dfn}

\begin{ex} Let $X$ be a manifold, and $E\ra X$  a vector bundle.
Write $C^\iy(E)$ for the vector space of smooth sections $e$ of $E$.
Then $C^\iy(X)$ acts on $C^\iy(E)$ by multiplication, so $C^\iy(E)$
is a $C^\iy(X)$-module.
\label{sd2ex3}
\end{ex}

In \cite[\S 5.3]{Joyc3} we define the {\it cotangent module\/}
$\Om_\fC$ of a $C^\iy$-ring~$\fC$.

\begin{dfn} Let $\fC$ be a $C^\iy$-ring, and $M$ a
$\fC$-module. A $C^\iy$-{\it derivation} is an $\R$-linear map
$\d:\fC\ra M$ such that whenever $f:\R^n\ra\R$ is a smooth map and
$c_1,\ldots,c_n\in\fC$, we have
\begin{equation*}
\d\Phi_f(c_1,\ldots,c_n)=\ts\sum_{i=1}^n\Phi_{\frac{\pd f}{\pd
x_i}}(c_1,\ldots,c_n)\cdot\d c_i.
\end{equation*}
We call such a pair $M,\d$ a {\it cotangent module\/} for $\fC$ if
it has the universal property that for any $\fC$-module $M'$ and
$C^\iy$-derivation $\d':\fC\ra M'$, there exists a unique morphism
of $\fC$-modules $\phi:M\ra M'$ with~$\d'=\phi\ci\d$.

Define $\Om_\fC$ to be the quotient of the free $\fC$-module with
basis of symbols $\d c$ for $c\in\fC$ by the $\fC$-submodule spanned
by all expressions of the form
$\d\Phi_f(c_1,\ldots,c_n)-\sum_{i=1}^n \Phi_{\frac{\pd f}{\pd
x_i}}(c_1,\ldots,c_n)\cdot\d c_i$ for $f:\R^n\ra\R$ smooth and
$c_1,\ldots,c_n\in\fC$, and define $\d_\fC:\fC\ra \Om_\fC$ by
$\d_\fC:c\mapsto\d c$. Then $\Om_\fC,\d_\fC$ is a cotangent module
for $\fC$. Thus cotangent modules always exist, and are unique up to
unique isomorphism.

Let $\fC,\fD$ be $C^\iy$-rings with cotangent modules
$\Om_\fC,\d_\fC$, $\Om_\fD,\d_\fD$, and $\phi:\fC\ra\fD$ be a
morphism of $C^\iy$-rings. Then $\phi$ makes $\Om_\fD$ into a
$\fC$-module, and there is a unique morphism
$\Om_\phi:\Om_\fC\ra\Om_\fD$ in $\fC$-mod with
$\d_\fD\ci\phi=\Om_\phi\ci\d_\fC$. This induces a morphism
$(\Om_\phi)_*:\Om_\fC\ot_\fC\fD\ra\Om_\fD$ in $\fD$-mod
with~$(\Om_\phi)_*\ci (\d_\fC\ot\id_\fD)=\d_\fD$.
\label{sd2def6}
\end{dfn}

\begin{ex} Let $X$ be a manifold. Then the cotangent bundle $T^*X$
is a vector bundle over $X$, so as in Example \ref{sd2ex3} it yields
a $C^\iy(X)$-module $C^\iy(T^*X)$. The exterior derivative
$\d:C^\iy(X)\ra C^\iy(T^*X)$ is a $C^\iy$-derivation. These
$C^\iy(T^*X),\d$ have the universal property in Definition
\ref{sd2def6}, and so form a {\it cotangent module\/}
for~$C^\iy(X)$.

Now let $X,Y$ be manifolds, and $f:X\ra Y$ be smooth. Then
$f^*(TY),TX$ are vector bundles over $X$, and the derivative of $f$
is a vector bundle morphism $\d f:TX\ra f^*(TY)$. The dual of this
morphism is $\d f^*:f^*(T^*Y)\ra T^*X$. This induces a morphism of
$C^\iy(X)$-modules $(\d f^*)_*:C^\iy\bigl(f^*(T^*Y)\bigr)\ra
C^\iy(T^*X)$. This $(\d f^*)_*$ is identified with $(\Om_{f^*})_*$
in Definition \ref{sd2def6} under the natural isomorphism
$C^\iy\bigl(f^*(T^*Y)\bigr)\cong C^\iy(T^*Y)\ot_{C^\iy(Y)}C^\iy(X)$.
\label{sd2ex4}
\end{ex}

Definition \ref{sd2def6} abstracts the notion of cotangent bundle of
a manifold in a way that makes sense for any $C^\iy$-ring.

\subsection{Quasicoherent sheaves on $C^\iy$-schemes}
\label{sd24}

In \cite[\S 6]{Joyc3} we discuss sheaves of modules on
$C^\iy$-schemes.

\begin{dfn} Let $\uX=(X,\O_X)$ be a $C^\iy$-scheme. An $\O_X$-{\it
module\/} $\cE$ on $\uX$ assigns a module $\cE(U)$ over $\O_X(U)$
for each open set $U\subseteq X$, with $\O_X(U)$-action
$\mu_U:\O_X(U)\t\cE(U)\ra\cE(U)$, and a linear map
$\cE_{UV}:\cE(U)\ra\cE(V)$ for each inclusion of open sets
$V\subseteq U\subseteq X$, such that the following commutes:
\begin{equation*}
\xymatrix@R=10pt@C=60pt{ \O_X(U)\t \cE(U) \ar[d]^{\rho_{UV}\t
\cE_{UV}} \ar[r]_{\mu_U} & \cE(U) \ar[d]_{\cE_{UV}} \\
\O_X(V)\t \cE(V)\ar[r]^{\mu_V} & \cE(V),}
\end{equation*}
and all this data $\cE(U),\cE_{UV}$ satisfies the usual sheaf
axioms~\cite[\S II.1]{Hart} .

A {\it morphism of\/ $\O_X$-modules\/} $\phi:\cE\ra\cF$ assigns a
morphism of $\O_X(U)$-modules $\phi(U):\cE(U)\ra\cF(U)$ for each
open set $U\subseteq X$, such that $\phi(V)\ci\cE_{UV}=
\cF_{UV}\ci\phi(U)$ for each inclusion of open sets $V\subseteq
U\subseteq X$. Then $\O_X$-modules form an abelian category, which
we write as~$\OXmod$.

As in \cite[\S 6.2]{Joyc3}, the spectrum functor $\Spec:\CRings^{\rm
op}\ab\ra\CSch$ has a counterpart for modules: if $\fC$ is a
$C^\iy$-ring and $(X,\O_X)=\Spec\fC$ we can define a functor
$\MSpec:\fCmod\ra\OXmod$. If $\fC$ is a {\it fair\/} $C^\iy$-ring,
there is a full abelian subcategory $\fCmod^{\rm co}$ of {\it
complete\/} $\fC$-modules in $\fCmod$, such that
$\MSpec\vert_{\fCmod^{\rm co}}:\fCmod^{\rm co}\ra\OXmod$ is an
equivalence of categories, with quasi-inverse the global sections
functor $\Ga:\OXmod\ra\fCmod^{\rm co}$. Let $\uX=(X,\O_X)$ be a
$C^\iy$-scheme, and $\cE$ an $\O_X$-module. We call $\cE$ {\it
quasicoherent\/} if $\uX$ can be covered by open $\uU$ with
$\uU\cong\Spec\fC$ for some $C^\iy$-ring $\fC$, and under this
identification $\cE\vert_U\cong\MSpec M$ for some $\fC$-module $M$.
We call $\cE$ a {\it vector bundle of rank\/} $n\ge 0$ if $\uX$ may
be covered by open $\uU$ such that
$\cE\vert_\uU\cong\O_U\ot_\R\R^n$.

Write $\qcoh(\uX),\vect(\uX)$ for the full subcategories of
quasicoherent sheaves and vector bundles in $\OXmod$. Then
$\qcoh(\uX)$ is an abelian category. Since $\MSpec:\fCmod^{\rm
co}\ra\OXmod$ is an equivalence for $\fC$ fair and
$(X,\O_X)=\Spec\fC$, as in \cite[Cor.~6.11]{Joyc3} we see that if
$\uX$ is a locally fair $C^\iy$-scheme then every $\O_X$-module
$\cE$ on $\uX$ is quasicoherent, that is,~$\qcoh(\uX)=\OXmod$.
\label{sd2def7}
\end{dfn}

\begin{rem} If $\uX$ is a separated, paracompact, locally fair
$C^\iy$-scheme then vector bundles on $\uX$ are projective objects
in the abelian category~$\qcoh(\uX)$.
\label{sd2rem2}
\end{rem}

\begin{dfn} Let $\uf:\uX\ra\uY$ be a morphism of $C^\iy$-schemes,
and let $\cE$ be an $\OY$-module. Define the {\it pullback\/}
$\uf^*(\cE)$, an $\OX$-module, by $\uf^*(\cE)=f^{-1}(\cE)
\ot_{f^{-1}(\OY)}\OX$, where $f^{-1}(\cE),f^{-1}(\OY)$ are inverse
image sheaves, and the tensor product uses the morphism
$f^\sh:f^{-1}(\O_Y)\ra\O_X$. If $\phi:\cE\ra\cF$ is a morphism in
$\OYmod$ we have an induced morphism $\uf^*(\phi)=f^{-1}(\phi)
\ot\id_{\OX}:\uf^*(\cE)\ra\uf^*(\cF)$ in $\OXmod$. Then
$\uf^*:\OYmod\ra\OXmod$ is a right exact functor, which restricts to
a right exact functor~$\uf^*:\qcoh(\uY)\ra\qcoh(\uX)$.
\label{sd2def8}
\end{dfn}

\begin{rem} Pullbacks $\uf^*(\cE)$ are characterized by a universal
property, and so are unique up to canonical isomorphism, rather than
unique. Our definition of $\uf^*(\cE)$ is not functorial in $\uf$.
That is, if $\uf:\uX\ra\uY$, $\ug:\uY\ra\uZ$ are morphisms and
$\cE\in\OZmod$ then $(\ug\ci\uf)^*(\cE)$ and $\uf^*(\ug^*(\cE))$ are
canonically isomorphic in $\OXmod$, but may not be equal. In
\cite{Joyc5} we keep track of these canonical isomorphisms, writing
them as $I_{\uf,\ug}(\cE):(\ug\ci\uf)^*(\cE)\ra\uf^*(\ug^*(\cE))$.
However, in this survey, by an abuse of notation that is common in
the literature, we will for simplicity omit the isomorphisms
$I_{\uf,\ug}(\cE)$, and identify $(\ug\ci\uf)^*(\cE)$
with~$\uf^*(\ug^*(\cE))$.

Similarly, when $\uf$ is the identity $\uid_{\smash{\uX}}:\uX\ra\uX$
and $\cE\in\OXmod$, we may not have $\uid^*_{\smash{\uX}}(\cE)=\cE$,
but there is a canonical isomorphism
$\de_{\smash{\uX}}(\cE):\uid^*_{\smash{\uX}}(\cE)\ra\cE$, which we
keep track of in \cite{Joyc5}. But here, for simplicity, by an abuse
of notation we omit $\de_{\smash{\uX}}(\cE)$, and identify
$\uid^*_{\smash{\uX}}(\cE)$ with~$\cE$.
\label{sd2rem3}
\end{rem}

\begin{ex} Let $X$ be a manifold, and $\uX$ the associated
$C^\iy$-scheme from Example \ref{sd2ex2}, so that $\O_X(U)=C^\iy(U)$
for all open $U\subseteq X$. Let $E\ra X$ be a vector bundle. Define
an $\OX$-module $\cE$ on $\uX$ by $\cE(U)=C^\iy(E\vert_U)$, the
smooth sections of the vector bundle $E\vert_U\ra U$, and for open
$V\subseteq U\subseteq X$ define $\cE_{UV}:\cE(U)\ra\cE(V)$ by
$\cE_{UV}:e_U\mapsto e_U\vert_V$. Then $\cE\in\vect(\uX)$ is a
vector bundle on $\uX$, which we think of as a lift of $E$ from
manifolds to $C^\iy$-schemes.

Let $f:X\ra Y$ be a smooth map of manifolds, and $\uf:\uX\ra\uY$ the
corresponding morphism of $C^\iy$-schemes. Let $F\ra Y$ be a vector
bundle over $Y$, so that $f^*(F)\ra X$ is a vector bundle over $X$.
Let $\cF\in\vect(\uY)$ be the vector bundle over $\uY$ lifting $F$.
Then $\uf^*(\cF)$ is canonically isomorphic to the vector bundle
over $\uX$ lifting~$f^*(F)$.
\label{sd2ex5}
\end{ex}

We define {\it cotangent sheaves}, the sheaf version of cotangent
modules in~\S\ref{sd23}.

\begin{dfn} Let $\uX$ be a $C^\iy$-scheme. Define
${\cal P}T^*\uX$ to associate to each open $U\subseteq X$ the
cotangent module $\Om_{\O_X(U)}$, and to each inclusion of open sets
$V\subseteq U\subseteq X$ the morphism of $\O_X(U)$-modules
$\Om_{\rho_{UV}}:\Om_{\O_X(U)}\ra\Om_{\O_X(V)}$ associated to the
morphism of $C^\iy$-rings $\rho_{UV}:\O_X(U)\ra\O_X(V)$. Then ${\cal
P}T^*\uX$ is a {\it presheaf of\/ $\O_X$-modules on\/} $\uX$. Define
the {\it cotangent sheaf\/ $T^*\uX$ of\/} $\uX$ to be the
sheafification of ${\cal P}T^*\uX$, as an $\OX$-module.

Let $\uf:\uX\ra\uY$ be a morphism of $C^\iy$-schemes. Then by
Definition \ref{sd2def8}, $\uf^*\bigl(T^*\uY\bigr)=f^{-1}(T^*\uY)
\ot_{f^{-1}(\O_Y)}\O_X,$ where $T^*\uY$ is the sheafification of the
presheaf $V\mapsto\Om_{\O_Y(V)}$, and $f^{-1}(T^*\uY)$ the
sheafification of the presheaf $U\mapsto\lim_{V\supseteq f(U)}
(T^*\uY)(V)$, and $f^{-1}(\O_Y)$ the sheafification of the presheaf
$U\mapsto\lim_{V\supseteq f(U)}\O_Y(V)$. The three sheafifications
combine into one, so that $\uf^*\bigl(T^*\uY\bigr)$ is the
sheafification of the presheaf ${\cal P}(\uf^*(T^*\uY))$ acting by
\begin{equation*}
U\longmapsto{\cal P}(\uf^*(T^*\uY))(U)=
\ts\lim_{V\supseteq f(U)}\Om_{\O_Y(V)}\ot_{\O_Y(V)}\O_X(U).
\end{equation*}

Define a morphism of presheaves ${\cal P}\Om_\uf:{\cal
P}(\uf^*(T^*\uY))\ra{\cal P}T^*\uX$ on $X$ by
\begin{equation*}
({\cal P}\Om_\uf)(U)=\ts\lim_{V\supseteq f(U)}
(\Om_{\rho_{f^{-1}(V)\,U}\ci f_\sh(V)})_*,
\end{equation*}
where $(\Om_{\rho_{f^{-1}(V)\,U}\ci f_\sh(V)})_*:\Om_{\O_Y(V)}
\ot_{\O_Y(V)}\O_X(U)\ra\Om_{\O_X(U)}=({\cal P} T^*\uX)(U)$ is
constructed as in Definition \ref{sd2def6} from the $C^\iy$-ring
morphisms $f_\sh(V):\O_Y(V)\ra\O_X(f^{-1}(V))$ from $f_\sh:\O_Y\ra
f_*(\O_X)$ corresponding to $f^\sh$ in $\uf$ as in \eq{sd2eq3}, and
$\rho_{f^{-1}(V)\,U}:\O_X(f^{-1}(V))\ra\O_X(U)$ in $\O_X$. Define
$\Om_\uf:\uf^*\bigl(T^*\uY\bigr)\ra T^*\uX$ to be the induced
morphism of the associated sheaves.
\label{sd2def9}
\end{dfn}

\begin{ex} Let $X$ be a manifold, and $\uX$ the associated
$C^\iy$-scheme. Then $T^*\uX$ is a vector bundle on $\uX$, and is
canonically isomorphic to the lift to $C^\iy$-schemes from Example
\ref{sd2ex5} of the cotangent vector bundle $T^*X$ of~$X$.
\label{sd2ex6}
\end{ex}

Here \cite[Th.~6.17]{Joyc3} are some properties of cotangent
sheaves.

\begin{thm}{\bf(a)} Let\/ $\uf:\uX\ra\uY$ and\/ $\ug:\uY\ra\uZ$ be
morphisms of\/ $C^\iy$-schemes. Then
\begin{equation*}
\Om_{\ug\ci\uf}=\Om_\uf\ci \uf^*(\Om_\ug)
\end{equation*}
as morphisms $(\ug\ci\uf)^*(T^*\uZ)\ra T^*\uX$ in $\OXmod$. Here
$\Om_\ug:\ug^*(T^*\uZ)\ra T^*\uY$ is a morphism in $\OYmod,$ so
applying $\uf^*$ gives $\uf^*(\Om_\ug):(\ug\ci\uf)^*(T^*\uZ)=
\uf^*(\ug^*(T^*\uZ))\ra \uf^*(T^*\uY)$ in\/~$\OXmod$.
\smallskip

\noindent{\bf(b)} Suppose\/ $\uW,\uX,\uY,\uZ$ are locally fair\/
$C^\iy$-schemes with a Cartesian square
\begin{equation*}
\xymatrix@C=60pt@R=13pt{ \uW \ar[r]_\uf \ar[d]^\ue & \uY \ar[d]_\uh \\
\uX \ar[r]^\ug & \uZ}
\end{equation*}
in $\CSchlf,$ so that\/ $\uW=\uX\t_\uZ\uY$. Then the following is
exact in $\qcoh(\uW)\!:$
\begin{equation*}
\xymatrix@C=15pt{ (\ug\ci\ue)^*(T^*\uZ)
\ar[rrrr]^(0.45){\ue^*(\Om_\ug)\op -\uf^*(\Om_\uh)} &&&&
\ue^*(T^*\uX)\!\op\!\uf^*(T^*\uY) \ar[rr]^(0.63){\Om_\ue\op\Om_\uf}
&& T^*\uW \ar[r] & 0.}
\end{equation*}
\label{sd2thm2}
\end{thm}

\section{The 2-category of d-spaces}
\label{sd3}

We will now define the 2-category of {\it d-spaces\/} $\dSpa$,
following \cite[\S 2]{Joyc5}. D-spaces are `derived' versions of
$C^\iy$-schemes. In \S\ref{sd4} we will define the 2-category of
d-manifolds $\dMan$ as a 2-subcategory of $\dSpa$. For an
introduction to 2-categories, see Appendix~\ref{sdA}.

\subsection{The definition of d-spaces}
\label{sd31}

\begin{dfn} A {\it d-space\/} $\bX$ is a quintuple
$\bX=(\uX,\OXp,\EX,\im_X,\jm_X)$ such that $\uX=(X,\OX)$ is a
separated, second countable, locally fair $C^\iy$-scheme, and
$\OXp,\EX,\im_X,\jm_X$ fit into an exact sequence of sheaves on~$X$
\begin{equation*}
\smash{\xymatrix@C=25pt{ \EX \ar[rr]^(0.45){\jm_X} && \OXp
\ar[rr]^(0.55){\im_X} && \OX \ar[r] & 0,}}
\end{equation*}
satisfying the conditions:
\begin{itemize}
\setlength{\itemsep}{0pt}
\setlength{\parsep}{0pt}
\item[(a)] $\OXp$ is a sheaf of $C^\iy$-rings on $X$, with
$\uX'=(X,\OXp)$ a $C^\iy$-scheme.
\item[(b)] $\im_X:\OXp\ra\OX$ is a surjective morphism of sheaves of
$C^\iy$-rings on $X$. Its kernel $\IX$ is a sheaf of ideals in
$\OXp$, which should be a sheaf of square zero ideals. Here a
{\it square zero ideal\/} in a commutative $\R$-algebra $A$ is
an ideal $I$ with $i\cdot j=0$ for all $i,j\in I$. Then $\IX$ is
an $\OXp$-module, but as $\IX$ consists of square zero ideals
and $\im_X$ is surjective, the $\OXp$-action factors through an
$\OX$-action. Hence $\IX$ is an $\OX$-module, and thus a
quasicoherent sheaf on $\uX$, as $\uX$ is locally fair.
\item[(c)] $\EX$ is a quasicoherent sheaf on $\uX$, and
$\jm_X:\EX\ra\IX$ is a surjective morphism in~$\qcoh(\uX)$.
\end{itemize}
As $\uX$ is locally fair, the underlying topological space $X$ is
locally homeomorphic to a closed subset of $\R^n$, so it is {\it
locally compact}. But Hausdorff, second countable and locally
compact imply paracompact, and thus $\uX$ is {\it paracompact}.

The sheaf of $C^\iy$-rings $\OXp$ has a sheaf of cotangent modules
$\Om_{\OXp}$, which is an $\OXp$-module with exterior derivative
$\d:\OXp\ra\Om_{\OXp}$. Define $\FX=\Om_{\OXp}\ot_\OXp\OX$ to be the
associated $\OX$-module, a quasicoherent sheaf on $\uX$, and set
$\psi_X=\Om_{\im_X}\ot\id:\FX\ra T^*\uX$, a morphism in
$\qcoh(\uX)$. Define $\phi_X:\EX\ra\FX$ to be the composition of
morphisms of sheaves of abelian groups on~$X$:
\begin{equation*}
\smash{\xymatrix@C=7pt{ \EX  \ar[rr]^{\jm_X} && \IX
\ar[rr]^{\d\vert_{\IX}} && \Om_{\OXp} \ar@{=}[r]^(0.31)\sim &
\Om_{\OXp}\ot_\OXp\OXp \ar[rrr]^{\id\ot\im_X} &&&
\Om_{\OXp}\ot_\OXp\OX \ar@{=}[r] & \FX. }}
\end{equation*}
It turns out that $\phi_X$ is actually a morphism of $\OX$-modules,
and the following sequence is exact in~$\qcoh(\uX)\!:$
\begin{equation*}
\xymatrix@C=20pt{ \EX \ar[rr]^{\phi_X} && \FX \ar[rr]^{\psi_X} &&
T^*\uX \ar[r] & 0.}
\end{equation*}
The morphism $\phi_X:\EX\ra\FX$ will be called the {\it virtual
cotangent sheaf\/} of $\bX$, for reasons we explain in~\S\ref{sd43}.

Let $\bX,\bY$ be d-spaces. A 1-{\it morphism\/} $\bs f:\bX\ra\bY$ is
a triple $\bs f=(\uf,f',f'')$, where $\uf=(f,f^\sh):\uX\ra\uY$ is a
morphism of $C^\iy$-schemes, $f':f^{-1}(\OYp)\ra\OXp$ a morphism of
sheaves of $C^\iy$-rings on $X$, and $f'':\uf^*(\EY)\ra\EX$ a
morphism in $\qcoh(\uX)$, such that the following diagram of sheaves
on $X$ commutes:
\begin{equation*}
\xymatrix@C=11pt@R=1pt{
f^{-1}\!(\EY)\ot_{f^{-1}(\OY)}^{\id}\!f^{-1}\!(\OY) \ar@{=}[r]
\ar[dd]^(0.4){{}\,\id\ot f^\sh} & f^{-1}\!(\EY)
\ar[rr]_(0.45){\raisebox{-9pt}{$\scriptstyle f^{-1}(\jm_Y)$}}
 && f^{-1}\!(\OYp)
\ar[rr]_(0.52){\raisebox{-9pt}{$\scriptstyle f^{-1}(\im_Y)$}}
\ar[ddd]^{f'} && f^{-1}\!(\OY) \ar[r] \ar[ddd]_{f^\sh} & 0 \\ \\
{\begin{subarray}{l}\ts \uf^*(\EY)=\\
\ts f^{-1}(\EY) \ot_{f^{-1}(\OY)}^{f^\sh}\OX\end{subarray}}
\ar[dr]^(0.7){f''} \\  & \EX \ar[rr]^{\jm_X} && \OXp
\ar[rr]^(0.55){\im_X} && \OX \ar[r] &  {0.\!} }
\end{equation*}
Define morphisms $f^2=\Om_{f'}\ot\id:\uf^*(\FY)\ra\FX$ and
$f^3=\Om_\uf:\uf^*(T^*\uY)\ra T^*\uX$ in $\qcoh(\uX)$. Then the
following commutes in $\qcoh(\uX)$, with exact rows:
\e
\begin{gathered}
\xymatrix@C=20pt@R=12pt{ \uf^*(\EY) \ar[rr]_{\uf^*(\phi_Y)}
\ar[d]^{f''} && \uf^*(\FY) \ar[rr]_{\uf^*(\psi_Y)} \ar[d]^{f^2} &&
\uf^*(T^*\uY) \ar[r] \ar[d]^{f^3} & 0 \\
\EX \ar[rr]^{\phi_X} && \FX \ar[rr]^{\psi_X} && T^*\uX \ar[r] & 0. }
\end{gathered}
\label{sd3eq1}
\e

If $\bX$ is a d-space, the {\it identity $1$-morphism\/}
$\bs\id_\bX:\bX\ra\bX$ is $\bs\id_\bX=\bigl(\uid_\uX,
\id_{\OXp},\id_{\EX}\bigr)$. Let $\bX,\bY,\bZ$ be d-spaces, and $\bs
f:\bX\ra\bY$, $\bs g:\bY\ra\bZ$ be 1-morphisms. Define the {\it
composition of\/ $1$-morphisms\/} $\bs g\ci\bs f:\bX\ra\bZ$ to be
\begin{equation*}
\bs g\ci\bs f=\bigl(\ug\ci\uf,f'\ci f^{-1}(g'),
f''\ci\uf^*(g'')\bigr).
\end{equation*}

Let $\bs f,\bs g:\bX\ra\bY$ be 1-morphisms of d-spaces, where $\bs
f=(\uf,f',f'')$ and $\bs g=(\ug,g',g'')$. Suppose $\uf=\ug$. A
2-{\it morphism\/} $\eta:\bs f\Ra\bs g$ is a morphism
$\eta:\uf^*(\FY)\ra\EX$ in $\qcoh(\uX)$, such that
\begin{gather*}
g'=f'+\jm_X\ci\eta\ci\bigl(\id\ot(f^\sh\ci
f^{-1}(\im_Y))\bigr) \ci\bigl(f^{-1}(\d)\bigr)\\
\text{and}\qquad g''=f''+\eta\ci\uf^*(\phi_Y).
\end{gather*}
Then $g^2=f^2+\phi_X\ci\eta$ and $g^3=f^3$, so \eq{sd3eq1} for $\bs
f,\bs g$ combine to give a diagram
\e
\begin{gathered}
\xymatrix@C=30pt@R=17pt{ \uf^*(\EY) \ar[rr]^{\uf^*(\phi_Y)}
\ar@<-.8ex>[d]_(0.3){f''}
\ar@<.2ex>[d]^(0.3){g''=f''+\eta\ci\uf^*(\phi_Y)} && \uf^*(\FY)
\ar[rr]^{\uf^*(\psi_Y)} \ar@<-.5ex>[d]_(0.45){f^2}
\ar@<.5ex>[d]^(0.45){g^2=f^2+\phi_X\ci\eta}
\ar@<.5ex>@{.>}[dll]^(0.3)\eta
 && \uf^*(T^*\uY) \ar[r] \ar[d]^{f^3=g^3} & 0 \\
\EX \ar@<-.2ex>[rr]^(0.6){\phi_X} && \FX \ar@<-.2ex>[rr]^{\psi_X} &&
T^*\uX \ar@<-.2ex>[r] & {0.\!\!} }
\end{gathered}\!\!\!\!\!
\label{sd3eq2}
\e
That is, $\eta$ is a homotopy between the morphisms of complexes
\eq{sd3eq1} from~$\bs f,\bs g$.

If $\bs f:\bX\ra\bY$ is a 1-morphism, the {\it identity\/
$2$-morphism\/} $\id_{\bs f}:\bs f\Ra\bs f$ is the zero morphism
$0:\uf^*(\FY)\ra\EX$. Suppose $\bX,\bY$ are d-spaces, $\bs f,\bs
g,\bs h:\bX\ra\bY$ are 1-morphisms and $\eta:\bs f\Ra\bs g$,
$\ze:\bs g\Ra\bs h$ are 2-morphisms. The {\it vertical composition
of\/ $2$-morphisms\/} $\ze\od\eta:\bs f\Ra\bs h$ as in \eq{sdAeq1}
is~$\ze\od\eta=\ze+\eta$.

Let $\bX,\bY,\bZ$ be d-spaces, $\bs f,\bs{\ti f}:\bX\ra\bY$ and $\bs
g,\bs{\ti g}:\bY\ra\bZ$ be 1-morphisms, and $\eta:\bs f\Ra\bs{\ti
f}$, $\ze:\bs g\Ra\bs{\ti g}$ be 2-morphisms. The {\it horizontal
composition of\/ $2$-morphisms\/} $\ze*\eta:\bs g\ci\bs f\Ra\bs{\ti
g}\ci\bs{\ti f}$ as in \eq{sdAeq2} is
\begin{equation*}
\ze*\eta=\eta\ci\uf^*(g^2)+f''\ci\uf^*(\ze)+\eta\ci\uf^*(\phi_Y)
\ci\uf^*(\ze).
\end{equation*}

Regard the category $\CSchlfssc$ of separated, second countable,
locally fair $C^\iy$-schemes as a 2-category with only identity
2-morphisms $\id_{\uf}$ for (1-)mor\-phisms $\uf:\uX\ra\uY$. Define
a 2-functor $F_\CSch^\dSpa:\CSchlfssc\ra\dSpa$ to map $\uX$ to
$\bX=(\uX,\OX,0,\id_{\OX},0)$ on objects $\uX$, to map $\uf$ to $\bs
f=(\uf,f^\sh,0)$ on (1-)morphisms $\uf:\uX\ra\uY$, and to map
identity 2-morphisms $\id_\uf:\uf\Ra\uf$ to identity 2-morphisms
$\id_{\bs f}:\bs f\Ra\bs f$. Define a 2-functor
$F_\Man^\dSpa:\Man\ra\dSpa$ by~$F_{\Man}^\dSpa=F_\CSch^\dSpa\ci
F_\Man^\CSch$.

Write $\hCSchlfssc$ for the full 2-subcategory of objects $\bX$ in
$\dSpa$ equivalent to $F_\CSch^\dSpa(\uX)$ for some $\uX$ in
$\CSchlfssc$, and $\hMan$ for the full 2-subcategory of objects
$\bX$ in $\dSpa$ equivalent to $F_\Man^\dSpa(X)$ for some manifold
$X$. When we say that a d-space $\bX$ {\it is a $C^\iy$-scheme}, or
{\it is a manifold}, we mean that $\bX\in\hCSchlfssc$, or
$\bX\in\hMan$, respectively.
\label{sd3def1}
\end{dfn}

In \cite[\S 2.2]{Joyc5} we prove:

\begin{thm}{\bf(a)} Definition\/ {\rm\ref{sd3def1}} defines a
strict\/ $2$-category $\dSpa,$ in which all\/ $2$-morphisms are
$2$-isomorphisms.
\smallskip

\noindent{\bf(b)} For any $1$-morphism $\bs f:\bX\ra\bY$ in\/
$\dSpa$ the\/ $2$-morphisms $\eta:\bs f\Ra\bs f$ form an abelian
group under vertical composition, and in fact a real vector space.
\smallskip

\noindent{\bf(c)} $F_\CSch^\dSpa$ and\/ $F_\Man^\dSpa$ in
Definition\/ {\rm\ref{sd3def1}} are full and faithful strict\/
$2$-functors. Hence $\CSchlfssc,\Man$ and\/ $\hCSchlfssc,\hMan$ are
equivalent\/ $2$-categories.
\label{sd3thm1}
\end{thm}

\begin{rem} One should think of a d-space $\bX=(\uX,\OXp,\EX,\im_X,
\jm_X)$ as being a $C^\iy$-scheme $\uX$, which is the `classical'
part of $\bX$ and lives in a category rather than a 2-category,
together with some extra `derived' information
$\OXp,\EX,\im_X,\jm_X$. 2-morphisms in $\dSpa$ are wholly to do with
this derived part. The sheaf $\EX$ may be thought of as a (dual)
`obstruction sheaf' on~$\uX$.
\label{sd3rem1}
\end{rem}

\subsection{Gluing d-spaces by equivalences}
\label{sd32}

Next we discuss gluing of d-spaces and 1-morphisms on open
d-subspaces.

\begin{dfn} Let $\bX=(\uX,\OXp,\EX,\im_X,\jm_X)$ be a d-space.
Suppose $\uU\subseteq\uX$ is an open $C^\iy$-subscheme. Then $\bU=
\bigl(\uU,\OXp\vert_\uU,\EX\vert_\uU,\im_X\vert_\uU,\jm_X\vert_\uU
\bigr)$ is a d-space. We call $\bU$ an {\it open d-subspace\/} of
$\bX$. An {\it open cover\/} of a d-space $\bX$ is a family
$\{\bU_a:a\in A\}$ of open d-subspaces $\bU_a$ of $\bX$
with~$\uX=\bigcup_{a\in A}\uU_a$.
\label{sd3def2}
\end{dfn}

As in \cite[\S 2.4]{Joyc5}, we can glue 1-morphisms on open
d-subspaces which are 2-isomorphic on the overlap. The proof uses
partitions of unity, as in~\S\ref{sd22}.

\begin{prop} Suppose $\bX,\bY$ are d-spaces, $\bU,\bV\subseteq\bX$
are open d-subspaces with\/ $\bX=\bU\cup\bV,$ $\bs f:\bU\ra\bY$
and\/ $\bs g:\bV\ra\bY$ are $1$-morphisms, and\/ $\eta:\bs
f\vert_{\bU\cap\bV}\Ra\bs g\vert_{\bU\cap\bV}$ is a $2$-morphism.
Then there exists a $1$-morphism $\bs h:\bX\ra\bY$ and\/
$2$-morphisms $\ze:\bs h\vert_\bU\Ra\bs f,$ $\th:\bs
h\vert_\bV\Ra\bs g$ such that\/ $\th\vert_{\bU\cap\bV}=
\eta\od\ze\vert_{\bU\cap\bV}:\bs h\vert_{\bU\cap\bV}\Ra \bs
g\vert_{\bU\cap\bV}$. This $\bs h$ is unique up to $2$-isomorphism,
and independent up to $2$-isomorphism of the choice of\/~$\eta$.
\label{sd3prop}
\end{prop}

{\it Equivalences\/} $\bs f:\bX\ra\bY$ in a 2-category are defined
in Appendix \ref{sdA}, and are the natural notion of when two
objects $\bX,\bY$ are `the same'. In \cite[\S 2.4]{Joyc5} we prove
theorems on gluing d-spaces by equivalences. See Spivak
\cite[Lem.~6.8 \& Prop.~6.9]{Spiv} for results similar to Theorem
\ref{sd3thm2} for his `local $C^\iy$-ringed spaces', an
$\iy$-categorical analogue of our d-spaces.

\begin{thm} Suppose $\bX,\bY$ are d-spaces, $\bU\subseteq\bX,$
$\bV\subseteq\bY$ are open d-subspaces, and\/ $\bs f:\bU\ra\bV$ is
an equivalence in $\dSpa$. At the level of topological spaces, we
have open $U\subseteq X,$ $V\subseteq Y$ with a homeomorphism
$f:U\ra V,$ so we can form the quotient topological space
$Z:=X\amalg_fY=(X\amalg Y)/\sim,$ where the equivalence relation
$\sim$ on $X\amalg Y$ identifies $u\in U\subseteq X$ with\/~$f(u)\in
V\subseteq Y$.

Suppose $Z$ is Hausdorff. Then there exist a d-space $\bZ$ with
topological space $Z,$ open d-subspaces $\bs{\hat X},\bs{\hat Y}$ in
$\bZ$ with\/ $\bZ=\bs{\hat X}\cup\bs{\hat Y},$ equivalences $\bs
g:\bX\ra\bs{\hat X}$ and\/ $\bs h:\bY\ra\bs{\hat Y}$ in $\dSpa$ such
that\/ $\bs g\vert_\bU$ and\/ $\bs h\vert_\bV$ are both equivalences
with\/ $\bs{\hat X}\cap\bs{\hat Y},$ and a $2$-morphism $\eta:\bs
g\vert_\bU\Ra\bs h\ci\bs f:\bU\ra\bs{\hat X}\cap\bs{\hat Y}$.
Furthermore, $\bZ$ is independent of choices up to equivalence.
\label{sd3thm2}
\end{thm}

\begin{thm} Suppose\/ $I$ is an indexing set, and\/ $<$ is a total order
on $I,$ and\/ $\bX_i$ for $i\in I$ are d-spaces, and for all\/ $i<j$
in $I$ we are given open d-subspaces $\bU_{ij}\subseteq\bX_i,$
$\bU_{ji}\subseteq\bX_j$ and an equivalence $\bs
e_{ij}:\bU_{ij}\ra\bU_{ji},$ such that for all\/ $i<j<k$ in $I$ we
have a $2$-commutative diagram
\begin{equation*}
\xymatrix@C=70pt@R=10pt{ & \bU_{ji}\cap\bU_{jk} \ar@<.5ex>[dr]^{\bs
e_{jk}\vert_{\bU_{ji}\cap\bU_{jk}}} \ar@{=>}[d]^{\eta_{ijk}} \\
\bU_{ij}\cap\bU_{ik} \ar@<.5ex>[ur]^{\bs
e_{ij}\vert_{\bU_{ij}\cap\bU_{ik}}} \ar@<-.25ex>[rr]^(0.37){\bs
e_{ik}\vert_{\bU_{ij}\cap\bU_{ik}}} && \bU_{ki}\cap\bU_{kj}}
\end{equation*}
for some $\eta_{ijk},$ where all three $1$-morphisms are
equivalences.

On the level of topological spaces, define the quotient topological
space $Y=(\coprod_{i\in I}X_i)/\sim,$ where $\sim$ is the
equivalence relation generated by $x_i\sim x_j$ if\/ $i<j,$ $x_i\in
U_{ij}\subseteq X_i$ and\/ $x_j\in U_{ji}\subseteq X_j$ with\/
$e_{ij}(x_i)=x_j$. Suppose $Y$ is Hausdorff and second countable.
Then there exist a d-space $\bY$ and a $1$-morphism $\bs
f_i:\bX_i\ra\bY$ which is an equivalence with an open d-subspace
$\bs{\hat X}_i\subseteq\bY$ for all\/ $i\in I,$ where
$\bY=\bigcup_{i\in I}\bs{\hat X}_i,$ such that\/ $\bs
f_i\vert_{\bU_{ij}}$ is an equivalence $\bU_{ij}\ra\bs{\hat
X}_i\cap\bs{\hat X}_j$ for all\/ $i<j$ in $I,$ and there exists a
$2$-morphism\/ $\eta_{ij}:\bs f_j\ci\bs e_{ij}\Ra\bs
f_i\vert_{\bU_{ij}}$. The d-space $\bY$ is unique up to equivalence,
and is independent of choice of\/ $2$-morphisms\/~$\eta_{ijk}$.

Suppose also that\/ $\bZ$ is a d-space, and\/ $\bs g_i:\bX_i\ra\bZ$
are $1$-morphisms for all\/ $i\in I,$ and there exist\/
$2$-morphisms $\ze_{ij}:\bs g_j\ci\bs e_{ij}\Ra\bs
g_i\vert_{\bU_{ij}}$ for all\/ $i<j$ in $I$. Then there exist a
$1$-morphism $\bs h:\bY\ra\bZ$ and\/ $2$-morphisms $\ze_i:\bs
h\ci\bs f_i\Ra\bs g_i$ for all\/ $i\in I$. The $1$-morphism $\bs h$
is unique up to $2$-isomorphism, and is independent of the choice
of\/ $2$-morphisms~$\ze_{ij}$.
\label{sd3thm3}
\end{thm}

\begin{rem} In Proposition \ref{sd3prop}, it is surprising that
$\bs h$ is independent of $\eta$ up to $2$-isomorphism. It holds
because of the existence of {\it partitions of unity\/} on nice
$C^\iy$-schemes, as in Proposition \ref{sd2prop}. Here is a sketch
proof: suppose $\eta,\bs h,\ze,\th$ and $\eta',\bs h',\ze',\th'$ are
alternative choices in Proposition \ref{sd3prop}. Then we have
2-morphisms $(\ze')^{-1}\od\ze:\bs h\vert_\bU\Ra \bs h'\vert_\bU$
and $(\th')^{-1}\od\th:\bs h\vert_\bV\Ra \bs h'\vert_\bV$. Choose a
partition of unity $\{\al,1-\al\}$ on $\uX$ subordinate to
$\{\uU,\uV\}$, so that $\al:\uX\ra\R$ is smooth with $\al$ supported
on $\uU\subseteq\uX$ and $1-\al$ supported on $\uV\subseteq\uX$.
Then $\al\cdot\bigl((\ze')^{-1}\od\ze\bigr)+(1-\al)
\cdot\bigl((\th')^{-1} \od\th\bigr)$ is a 2-morphism $\bs h\Ra\bs
h'$, where $\al\cdot\bigl((\ze')^{-1}\od\ze\bigr)$ makes sense on
all of $\uX$ (rather than just on $\uU$ where $(\ze')^{-1}\od\ze$ is
defined) as $\al$ is supported on $\uU$, so we extend by zero
on~$\uX\sm\uU$.

Similarly, in Theorem \ref{sd3thm3}, the compatibility conditions on
the gluing data $\bX_i,\bU_{ij},\bs e_{ij}$ are significantly weaker
than you might expect, because of the existence of partitions of
unity. The 2-morphisms $\eta_{ijk}$ on overlaps
$\bX_i\cap\bX_j\cap\bX_k$ are only required to exist, not to satisfy
any further conditions. In particular, one might think that on
overlaps $\bX_i\cap\bX_j\cap\bX_k \cap\bX_l$ we should require
\e
\eta_{ikl}\od(\id_{{\bs f}_{kl}}*\eta_{ijk})\vert_{\bU_{ij}\cap
\bU_{ik}\cap\bU_{il}}=\eta_{ijl}\od (\eta_{jkl}*\id_{{\bs f}_{ij}})
\vert_{\bU_{ij}\cap \bU_{ik}\cap\bU_{il}},
\label{sd3eq3}
\e
but we do not. Also, one might expect the $\ze_{ij}$ should satisfy
conditions on triple overlaps $\bX_i\cap\bX_j\cap\bX_k$, but they
need not.

The moral is that constructing d-spaces by gluing together patches
$\bX_i$ is straightforward, as one only has to verify mild
conditions on triple overlaps $\bX_i\cap\bX_j\cap\bX_k$. Again, this
works because of the existence of partitions of unity on nice
$C^\iy$-schemes, which are used to construct the glued d-spaces
$\bZ$ and 1- and 2-morphisms in Theorems \ref{sd3thm2}
and~\ref{sd3thm3}.

In contrast, for gluing d-stacks in \cite[\S 9.4]{Joyc5}, we do need
compatibility conditions of the form \eq{sd3eq3}. The problem of
gluing geometric spaces in an $\iy$-category $\bs{\mathcal{C}}$ by
equivalences, such as Spivak's derived manifolds \cite{Spiv}, is
discussed by To\"en and Vezzosi \cite[\S 1.3.4]{ToVe} and Lurie
\cite[\S 6.1.2]{Luri1}. It requires nontrivial conditions on
overlaps $\bX_{i_1}\cap\cdots\cap\bX_{i_n}$ for all~$n=2,3,\ldots.$
\label{sd3rem2}
\end{rem}

\subsection{Fibre products in $\dSpa$}
\label{sd33}

{\it Fibre products\/} in 2-categories are explained in Appendix
\ref{sdA}. In \cite[\S 2.5--\S 2.6]{Joyc5} we discuss fibre products
in $\dSpa$, and their relation to transverse fibre products
in~$\Man$.

\begin{thm}{\bf(a)} All fibre products exist in the
$2$-category~$\dSpa$.
\smallskip

\noindent{\bf(b)} Let\/ $g:X\ra Z$ and\/ $h:Y\ra Z$ be smooth maps
of manifolds without boundary, and write\/ $\bX=F_\Man^\dSpa(X),$
and similarly for $\bY,\bZ,\bs g,\bs h$. If\/ $g,h$ are transverse,
so that a fibre product $X\t_{g,Z,h}Y$ exists in $\Man,$ then the
fibre product $\bX\t_{\bs g,\bZ,\bs h}\bY$ in $\dSpa$ is equivalent
in $\dSpa$ to\/ $F_\Man^\dSpa(X\t_{g,Z,h}Y)$. If\/ $g,h$ are not
transverse then $\bX\t_{\bs g,\bZ,\bs h}\bY$ exists in $\dSpa,$ but
is not a manifold.
\label{sd3thm4}
\end{thm}

To prove (a), given 1-morphisms $\bs g:\bX\ra\bZ$ and $\bs
h:\bY\ra\bZ$, we write down an explicit d-space
$\bW=(\uW,\OWp,\EW,\im_W,\jm_W)$, 1-morphisms $\bs
e=(\ue,e',e''):\bW\ra\bX$ and $\bs f=(\uf,f',f''):\bW\ra\bY$ and a
2-morphism $\eta:\bs g\ci\bs e\Ra\bs h\ci\bs f$, and verify the
universal property for
\begin{equation*}
\xymatrix@C=60pt@R=10pt{ \bW \ar[r]_(0.25){\bs f} \ar[d]^{\bs e}
\drtwocell_{}\omit^{}\omit{^{\eta}}
 & \bY \ar[d]_{\bs h} \\ \bX \ar[r]^(0.7){\bs g} & \bZ}
\end{equation*}
to be a 2-Cartesian square in $\dSpa$. The underlying $C^\iy$-scheme
$\uW$ is the fibre product $\uW=\uX\t_{\ug,\uZ,\uh}\uY$ in $\CSch$,
and $\ue:\uW\ra\uX$, $\uf:\uW\ra\uY$ are the projections from the
fibre product. The definitions of $\OWp,\im_W,\jm_W,e',f'$ are
complex, and we will not give them here. The remaining data
$\EW,e'',f'',\eta$, as well as the virtual cotangent sheaf
$\phi_W:\EW\ra\FW$, is characterized by the following commutative
diagram in $\qcoh(\uW)$, with exact top row:
\begin{equation*}
\xymatrix@C=15pt@R=30pt{ (\ug\ci\ue)^*(\EZ)
\ar[rrrr]^(0.47){\begin{pmatrix}\st \ue^*(g'') \\ \st -\uf^*(h'') \\
\st (\ug\ci\ue)^*(\phi_Z)\end{pmatrix}}
&&&&
{\vphantom{\EW}\smash{\raisebox{11pt}{$\begin{subarray}{l} \ts \ue^*(\EX)\op
\\ \ts \uf^*(\EY)\op \\ \ts (\ug\ci\ue)^*(\FZ)\end{subarray}$}}}
\ar[d]_(0.45){\begin{pmatrix}
\st -\ue^*(\phi_X)\!{} & \st 0\!{} & \st \ue^*(g^2) \\
\st 0\!{} & \st -\uf^*(\phi_Y)\!{} & \st -\uf^*(h^2) \end{pmatrix}}
\ar[rrr]^(0.55){\begin{pmatrix} \st e''\!\!\!\!{} & \st f''\!\!\!\!{}
& \st \eta \end{pmatrix}} &&& \cE_W \ar[r] \ar[d]^{\phi_W} & 0
\\
&&&& {\begin{subarray}{l} \ts\ue^*(\FX)\op \\
\ts \uf^*(\FY)\end{subarray}}
\ar[rrr]^(0.55){\begin{pmatrix} \st e^2\!\!\!\!{} & \st f^2
\end{pmatrix}}_(0.55)\cong
&&& \FW. }
\end{equation*}

\section{The 2-category of d-manifolds}
\label{sd4}

Sections \ref{sd41}--\ref{sd48} survey the results of \cite[\S 3--\S
4]{Joyc5} on d-manifolds. Section \ref{sd49} briefly describes
extensions to d-manifolds with boundary, d-manifolds with corners,
and d-orbifolds from \cite[\S 6--\S 12]{Joyc5}, and \S\ref{sd410}
discusses d-manifold bordism and virtual classes for d-manifolds and
d-orbifolds following \cite[\S 13]{Joyc5}. Section \ref{sd411}
explains the relationship between d-manifolds and d-orbifolds and
other classes of geometric spaces, summarizing~\cite[\S 14]{Joyc5}.

\subsection{The definition of d-manifolds}
\label{sd41}

\begin{dfn} A d-space $\bU$ is called a {\it principal d-manifold\/}
if is equivalent in $\dSpa$ to a fibre product $\bX\t_{\bs g,\bZ,\bs
h}\bY$ with $\bX,\bY,\bZ\in\hMan$. That is,
\begin{equation*}
\bU\simeq F_\Man^\dSpa(X)\t_{F_\Man^\dSpa(g),F_\Man^\dSpa(Z),
F_\Man^\dSpa(h)}F_\Man^\dSpa(Y)
\end{equation*}
for manifolds $X,Y,Z$ and smooth maps $g:X\ra Z$ and $h:Y\ra Z$. The
{\it virtual dimension\/} $\vdim\bU$ of $\bU$ is defined to be
$\vdim\bU=\dim X+\dim Y-\dim Z$. Proposition \ref{sd4prop2}(b) below
shows that if $\bU\ne\bs\es$ then $\vdim\bU$ depends only on the
d-space $\bU$, and not on the choice of $X,Y,Z,g,h$, and so is well
defined.

A d-space $\bW$ is called a {\it d-manifold of virtual dimension\/}
$n\in\Z$, written $\vdim\bW=n$, if $\bW$ can be covered by nonempty
open d-subspaces $\bU$ which are principal d-manifolds
with~$\vdim\bU=n$.

Write $\dMan$ for the full 2-subcategory of d-manifolds in $\dSpa$.
If $\bX\in\hMan$ then $\bX\simeq\bX\t_{\bs *}\bs *$, so $\bX$ is a
principal d-manifold, and thus a d-manifold. Therefore $\hMan$ is a
2-subcategory of $\dMan$. We say that a d-manifold $\bX$ {\it is a
manifold\/} if it lies in $\hMan$. The 2-functor
$F_\Man^\dSpa:\Man\ra\dSpa$ maps into $\dMan$, and we will write
$F_\Man^\dMan=F_\Man^\dSpa:\Man\ra\dMan$.
\label{sd4def1}
\end{dfn}

Here \cite[\S 3.2]{Joyc5} are alternative descriptions of principal
d-manifolds:

\begin{prop} The following are equivalent characterizations of when
a d-space $\bW$ is a principal d-manifold:
\begin{itemize}
\setlength{\itemsep}{0pt}
\setlength{\parsep}{0pt}
\item[{\bf(a)}] $\bW\simeq \bX\t_{\bs g,\bZ,\bs
h}\bY$ for $\bX,\bY,\bZ\in\hMan$.
\item[{\bf(b)}] $\bW\simeq \bX\t_{\bs i,\bZ,\bs
j}\bY,$ where $X,Y,Z$ are manifolds, $i:X\ra Z,$ $j:Y\ra Z$ are
embeddings, $\bX=F_\Man^\dSpa(X),$ and similarly for $Y,Z,i,j$.
That is, $\bW$ is an intersection of two submanifolds $X,Y$ in
$Z,$ in the sense of d-spaces.
\item[{\bf(c)}] $\bW\simeq \bV\t_{\bs s,\bs E,\bs 0}\bV,$ where $V$
is a manifold, $E\ra V$ is a vector bundle, $s:V\ra E$ is a
smooth section, $0:V\ra E$ is the zero section,
$\bV=F_\Man^\dSpa(V),$ and similarly for $E,s,0$. That is, $\bW$
is the zeroes $s^{-1}(0)$ of a smooth section $s$ of a vector
bundle $E,$ in the sense of d-spaces.
\end{itemize}
\label{sd4prop1}
\end{prop}

\subsection{`Standard model' d-manifolds, 1- and 2-morphisms}
\label{sd42}

The next three examples, taken from \cite[\S 3.2 \& \S 3.4]{Joyc5},
give explicit models for principal d-manifolds in the form
$\bV\t_{\bs s,\bs E,\bs 0}\bV$ from Proposition \ref{sd4prop1}(c)
and their 1- and 2-morphisms, which we call {\it standard models}.

\begin{ex} Let $V$ be a manifold, $E\ra V$ a vector bundle (which we
sometimes call the {\it obstruction bundle\/}), and $s\in C^\iy(E)$.
We will write down an explicit principal d-manifold
$\bS=(\uS,\OSp,\ES,\im_S,\jm_S)$ which is equivalent to $\bV\t_{\bs
s,\bs E,\bs 0}\bV$ in Proposition \ref{sd4prop1}(c). We call $\bS$
the {\it standard model\/} of $(V,E,s)$, and also write it
$\bS_{V,E,s}$. Proposition \ref{sd4prop1} shows that every principal
d-manifold $\bW$ is equivalent to $\bS_{V,E,s}$ for some~$V,E,s$.

Write $C^\iy(V)$ for the $C^\iy$-ring of smooth functions
$c:V\ra\R$, and $C^\iy(E),\ab C^\iy(E^*)$ for the vector spaces of
smooth sections of $E,E^*$ over $V$. Then $s$ lies in $C^\iy(E)$,
and $C^\iy(E),C^\iy(E^*)$ are modules over $C^\iy(V)$, and there is
a natural bilinear product $\cdot\,:C^\iy(E^*)\t C^\iy(E)\ra
C^\iy(V)$. Define $I_s\subseteq C^\iy(V)$ to be the ideal generated
by $s$. That is,
\e
I_s=\bigl\{\al\cdot s:\al\in C^\iy(E^*)\bigr\}\subseteq C^\iy(V).
\label{sd4eq1}
\e

Let $I_s^2=\an{fg:f,g\in I_s}_\R$ be the square of $I_s$. Then
$I_s^2$ is an ideal in $C^\iy(V)$, the ideal generated by $s\ot s\in
C^\iy(E\ot E)$. That is,
\begin{equation*}
I_s^2=\bigl\{\be\cdot (s\ot s):\be\in C^\iy(E^*\ot E^*)\bigr\}
\subseteq C^\iy(V).
\end{equation*}
Define $C^\iy$-rings $\fC=C^\iy(V)/I_s$, $\fC'=C^\iy(V)/I_s^2$, and
let $\pi:\fC'\ra\fC$ be the natural projection from the inclusion
$I_s^2\subseteq I_s$. Define a topological space $S=\{v\in
V:s(v)=0\}$, as a subspace of $V$. Now $s(v)=0$ if and only if
$(s\ot s)(v)=0$. Thus $S$ is the underlying topological space for
both $\Spec\fC$ and $\Spec\fC'$. So $\Spec\fC=\uS=(S,\OS)$,
$\Spec\fC'=\uS'=(S,\OSp)$, and $\Spec\pi=\uim_S=(\id_S,\im_S):
\uS'\ra\uS$, where $\uS,\uS'$ are fair affine $C^\iy$-schemes, and
$\OS,\OSp$ are sheaves of $C^\iy$-rings on $S$, and
$\im_S:\OSp\ra\OS$ is a morphism of sheaves of $C^\iy$-rings. Since
$\pi$ is surjective with kernel the square zero ideal $I_s/I_s^2$,
$\im_S$ is surjective, with kernel $\IS$ a sheaf of square zero
ideals in~$\OSp$.

From \eq{sd4eq1} we have a surjective $C^\iy(V)$-module morphism
$C^\iy(E^*)\ra I_s$ mapping $\al\mapsto \al\cdot s$. Applying
$\ot_{C^\iy(V)}\fC$ gives a surjective $\fC$-module morphism
\begin{equation*}
\si:C^\iy(E^*)/(I_s\cdot C^\iy(E^*))\longra I_s/I_s^2,\quad
\si: \al+(I_s\cdot C^\iy(E^*))\longmapsto \al\cdot s+I_s^2.
\end{equation*}
Define $\ES=\MSpec\bigl(C^\iy(E^*)/(I_s\cdot C^\iy(E^*))\bigr)$.
Also $\MSpec(I_s/I_s^2)=\IS$, so $\jm_S=\MSpec\si$ is a surjective
morphism $\jm_S:\ES\ra\IS$ in $\qcoh(\uS)$. Therefore
$\bS_{V,E,s}=\bS=(\uS,\OSp,\ES,\im_S,\jm_S)$ is a d-space.

In fact $\cE_S$ is a vector bundle on $\uS$ naturally isomorphic to
$\uE^*\vert_\uS$, where $\uE$ is the vector bundle on
$\uV=F_\Man^\CSch(V)$ corresponding to $E\ra V$. Also $\FS\cong
T^*\uV\vert_\uS$. The morphism $\phi_S:\ES\ra\FS$ can be interpreted
as follows: choose a connection $\nabla$ on $E\ra V$. Then $\nabla
s\in C^\iy(E\ot T^*V)$, so we can regard $\nabla s$ as a morphism of
vector bundles $E^*\ra T^*V$ on $V$. This lifts to a morphism of
vector bundles $\hat\nabla s:\uE^*\ra T^*\uV$ on the $C^\iy$-scheme
$\uV$, and $\phi_S$ is identified with $\hat\nabla
s\vert_\uS:\uE^*\vert_\uS\ra T^*\uV\vert_\uS$ under the isomorphisms
$\ES\cong\uE^*\vert_\uS$, $\FS\cong T^*\uV\vert_\uS$.
\label{sd4ex1}
\end{ex}

Proposition \ref{sd4prop1} implies that every principal d-manifold
$\bW$ is equivalent to $\bS_{V,E,s}$ for some $V,E,s$. The notation
$O(s)$ and $O(s^2)$ used below should be interpreted as follows. Let
$V$ be a manifold, $E\ra V$ a vector bundle, and $s\in C^\iy(E)$. If
$F\ra V$ is another vector bundle and $t\in C^\iy(F)$, then we write
$t=O(s)$ if $t=\al\cdot s$ for some $\al\in C^\iy(F\ot E^*)$, and
$t=O(s^2)$ if $t=\be\cdot (s\ot s)$ for some $\be\in C^\iy(F\ot
E^*\ot E^*)$. Similarly, if $W$ is a manifold and $f,g:V\ra W$ are
smooth then we write $f=g+O(s)$ if $c\ci f-c\ci g=O(s)$ for all
smooth $c:W\ra\R$, and $f=g+O(s^2)$ if $c\ci f-c\ci g=O(s^2)$ for
all~$c$.

\begin{ex} Let $V,W$ be manifolds, $E\ra V$, $F\ra W$ be vector
bundles, and $s\in C^\iy(E)$, $t\in C^\iy(F)$. Write
$\bX=\bS_{V,E,s}$, $\bY=\bS_{W,F,t}$ for the `standard model'
principal d-manifolds from Example \ref{sd4ex1}. Suppose $f:V\ra W$
is a smooth map, and $\hat f:E\ra f^*(F)$ is a morphism of vector
bundles on $V$ satisfying
\e
\hat f\ci s= f^*(t)+O(s^2)\quad\text{in $C^\iy\bigl(f^*(F)\bigr)$.}
\label{sd4eq2}
\e

We will define a 1-morphism $\bs g=(\ug,g',g''):\bX\ra\bY$ in
$\dMan$ using $f,\hat f$. We will also write $\bs g:\bX\ra\bY$ as
$\bS_{\smash{f,\hat f}}:\bS_{V,E,s}\ra \bS_{W,F,t}$, and call it a
{\it standard model\/ $1$-morphism}. If $x\in X$ then $x\in V$ with
$s(x)=0$, so \eq{sd4eq2} implies that
\begin{equation*}
t\bigl(f(x)\bigr)=\bigl(f^*(t)\bigr)(x)=
\hat f\bigl(s(x)\bigr)+O\bigl(s(x)^2\bigr)=0,
\end{equation*}
so $f(x)\in Y\subseteq W$. Thus $g:=f\vert_X$ maps~$X\ra Y$.

Define morphisms of $C^\iy$-rings
\begin{gather*}
\phi:C^\iy(W)/I_t\longra C^\iy(V)/I_s,\quad
\phi':C^\iy(W)/I_t^2\longra C^\iy(V)/I_s^2,\\
\text{by}\quad \phi:c+I_t\longmapsto c\ci f+I_s,\quad
\phi':c+I_t^2\longmapsto c\ci f+I_s^2.
\end{gather*}
Here $\phi$ is well-defined since if $c\in I_t$ then $c=\ga\cdot t$
for some $\ga\in C^\iy(F^*)$, so
\begin{equation*}
c\ci f\!=\!(\ga\cdot t)\ci f\!=\!f^*(\ga)\cdot f^*(t)\!=\!f^*(\ga)
\cdot\bigl(\hat f\ci s+O(s^2)\bigr)\!=\!\bigl(\hat f\ci f^*(\ga)
\bigr)\cdot s+O(s^2)\in I_s.
\end{equation*}
Similarly if $c\in I_t^2$ then $c\ci f\in I_s^2$, so $\phi'$ is
well-defined. Thus we have $C^\iy$-scheme morphisms
$\ug=(g,g^\sh)=\Spec\phi:\uX\ra\uY$, and
$(g,g')=\Spec\phi':(X,\OXp)\ra(Y,\OYp)$, both with underlying map
$g$. Hence $g^\sh:g^{-1}(\OY)\ra\OX$ and $g':g^{-1}(\OYp)\ra\OXp$
are morphisms of sheaves of $C^\iy$-rings on~$X$.

Since $\ug^*(\EY)=\MSpec\bigl(C^\iy(f^*(F^*))/(I_s\cdot
C^\iy(f^*(F^*))\bigr)$, we may define $g'':\ug^*(\EY)\ra\EX$ by
$g''=\MSpec(G'')$, where
\begin{align*}
G'':C^\iy(f^*(F^*))/(I_s\cdot C^\iy(f^*(F^*))\longra
C^\iy(E^*)/(I_s\cdot C^\iy(E^*))\\
\text{is defined by}\quad G'':\ga+I_s\cdot C^\iy(f^*(F^*))\longmapsto
\ga\ci\hat f+I_s\cdot C^\iy(E^*).
\end{align*}
This defines $\bs g=(\ug,g',g'')$. One can show it is a 1-morphism
$\bs g:\bX\ra\bY$ in $\dSpa$, which we also write
as~$\bS_{\smash{f,\hat f}}:\bS_{V,E,s}\ra \bS_{W,F,t}$.

Now suppose $\ti V$ is an open neighbourhood of $s^{-1}(0)$ in $V$,
and let $\ti E=E\vert_{\smash{\ti V}}$ and $\ti s=s\vert_{\smash{\ti
V}}$. Write $i_{\smash{\ti V}}:\ti V\ra V$ for the inclusion. Then
$i_{\ti V}^*(E)=\ti E$, and $\id_{\ti E}\ci\ti s=\ti s= i_{\ti
V}^*(s)$. Thus we have a 1-morphism $\bs i_{\smash{\ti
V,V}}=\bS_{\smash{i_{\smash{\ti V}},\id_{\smash{\ti
E}}}}:\bS_{\smash{\ti V,\ti E,\ti s}}\ra\bS_{V,E,s}$. It is easy to
show that $\bs i_{\smash{\ti V,V}}$ is a 1-{\it isomorphism}, with
an inverse $\bs i_{\smash{\ti V,V}}^{-1}$. That is, making $V$
smaller without making $s^{-1}(0)$ smaller does not really change
$\bS_{V,E,s}$; the d-manifold $\bS_{V,E,s}$ depends only on $E,s$ on
an arbitrarily small open neighbourhood of $s^{-1}(0)$ in~$V$.
\label{sd4ex2}
\end{ex}

\begin{ex} Let $V,W$ be manifolds, $E\ra V$, $F\ra W$ be vector
bundles, and $s\in C^\iy(E)$, $t\in C^\iy(F)$. Suppose $f,g:V\ra W$
are smooth and $\hat f:E\ra f^*(F),$ $\hat g:E\ra g^*(F)$ are vector
bundle morphisms with $\hat f\ci s= f^*(t)+O(s^2)$ and $\hat g\ci s=
g^*(t)+O(s^2),$ so we have 1-morphisms $\bS_{\smash{f,\hat
f}},\bS_{\smash{g,\hat g}}:\bS_{V,E,s}\ra \bS_{W,F,t}$. It is easy
to show that $\bS_{\smash{f,\hat f}}=\bS_{\smash{g,\hat g}}$ if and
only if $g=f+O(s^2)$ and~$\hat g=\hat f+O(s)$.

Now suppose $\La:E\ra f^*(TW)$ is a morphism of vector bundles on
$V$. Taking the dual of $\La$ and lifting to $\uV$ gives
$\La^*:\uf^*(T^*\uW)\ra\cE^*$. Restricting to the $C^\iy$-subscheme
$\uX=s^{-1}(0)$ in $\uV$ gives $\la=\La^*\vert_\uX:
\uf^*(\FY)\cong\uf^*(T^*\uW) \vert_\uX\ra\cE^*\vert_\uX=\EX$. One
can show that $\la$ is a 2-morphism $\bS_{\smash{f,\hat
f}}\Ra\bS_{\smash{g,\hat g}}$ if and only if
\begin{equation*}
g= f+\La\ci s+O(s^2)\quad\text{and}\quad \hat g=\hat f+f^*(\d
t)\ci\La+O(s).
\end{equation*}
We write $\la$ as $\bS_\La:\bS_{\smash{f,\hat f}}\Ra\bS_{g,\hat g}$,
and call it a {\it standard model\/ $2$-morphism}. Every 2-morphism
$\eta:\bS_{\smash{f,\hat f}}\Ra\bS_{\smash{g,\hat g}}$ is $S_\La$
for some $\La$. Two vector bundle morphisms $\La,\La':E\ra f^*(TW)$
have $S_\La=S_{\smash{\La'}}$ if and only if~$\La=\La'+O(s)$.
\label{sd4ex3}
\end{ex}

If $\bX$ is a d-manifold and $x\in\bX$ then $x$ has an open
neighbourhood $\bU$ in $\bX$ equivalent in $\dSpa$ to $\bS_{V,E,s}$
for some manifold $V$, vector bundle $E\ra V$ and $s\in C^\iy(E)$.
In \cite[\S 3.3]{Joyc5} we investigate the extent to which $\bX$
determines $V,E,s$ near a point in $\bX$ and $V$, and prove:

\begin{thm} Let\/ $\bX$ be a d-manifold, and\/ $x\in\bX$. Then there
exists an open neighbourhood\/ $\bU$ of\/ $x$ in $\bX$ and an
equivalence $\bU\simeq\bS_{V,E,s}$ in $\dMan$ for some manifold\/
$V,$ vector bundle $E\ra V$ and\/ $s\in C^\iy(E)$ which identifies
$x\in\bU$ with a point\/ $v\in V$ such that\/ $s(v)=\d s(v)=0,$
where $\bS_{V,E,s}$ is as in Example\/ {\rm\ref{sd4ex1}}. These
$V,E,s$ are determined up to non-canonical isomorphism near $v$ by
$\bX$ near $x,$ and in fact they depend only on the underlying\/
$C^\iy$-scheme $\uX$ and the integer\/~$\vdim\bX$.
\label{sd4thm1}
\end{thm}

Thus, if we impose the extra condition $\d s(v)=0$, which is in fact
equivalent to choosing $V,E,s$ with $\dim V$ as small as possible,
then $V,E,s$ are determined uniquely near $v$ by $\bX$ near $x$
(that is, $V,E,s$ are determined locally up to isomorphism, but not
up to canonical isomorphism). If we drop the condition $\d s(v)=0$
then $V,E,s$ are determined uniquely near $v$ by $\bX$ near $x$
and~$\dim V$.

Theorem \ref{sd4thm1} shows that any d-manifold
$\bX=(\uX,\OXp,\EX,\im_X,\jm_X)$ is determined up to equivalence in
$\dSpa$ near any point $x\in\bX$ by the `classical' underlying
$C^\iy$-scheme $\uX$ and the integer $\vdim\bX$. So we can ask: what
extra information about $\bX$ is contained in the `derived' data
$\OXp,\EX,\im_X,\jm_X$? One can think of this extra information as
like a vector bundle $\cE$ over $\uX$. The only local information in
a vector bundle $\cE$ is $\rank\cE\in\Z$, but globally it also
contains nontrivial algebraic-topological information.

Suppose now that $\bs f:\bX\ra\bY$ is a 1-morphism in $\dMan$, and
$x\in\bX$ with $\bs f(x)=y\in\bY$. Then by Theorem \ref{sd4thm1} we
have $\bX\simeq\bS_{V,E,s}$ near $x$ and $\bY\simeq\bS_{W,F,t}$ near
$y$. So up to composition with equivalences, we can identify $\bs f$
near $x$ with a 1-morphism $\bs g:\bS_{V,E,s}\ra\bS_{W,F,t}$. Thus,
to understand arbitrary 1-morphisms $\bs f$ in $\dMan$ near a point,
it is enough to study 1-morphisms $\bs g:\bS_{V,E,s}\ra\bS_{W,F,t}$.
Our next theorem, proved in \cite[\S 3.4]{Joyc5}, shows that after
making $V$ smaller, every 1-morphism $\bs g:\bS_{V,E,s}
\ra\bS_{W,F,t}$ is of the form~$\bS_{\smash{f,\hat f}}$.

\begin{thm} Let\/ $V,W$ be manifolds, $E\ra V,$ $F\ra W$ be vector
bundles, and\/ $s\in C^\iy(E),$ $t\in C^\iy(F)$. Define principal
d-manifolds\/ $\bX=\bS_{V,E,s},$ $\bY=\bS_{W,F,t},$ with topological
spaces $X=\{v\in V:s(v)=0\}$ and\/ $Y=\{w\in W:t(w)=0\}$. Suppose
$\bs g:\bX\ra\bY$ is a $1$-morphism. Then there exist an open
neighbourhood\/ $\ti V$ of\/ $X$ in $V,$ a smooth map $f:\ti V\ra
W,$ and a morphism of vector bundles $\hat f:\ti E\ra f^*(F)$ with\/
$\hat f\ci\ti s= f^*(t),$ where $\ti E=E\vert_{\ti V},$ $\ti
s=s\vert_{\ti V},$ such that\/ $\bs g=\bS_{\smash{f,\hat f}}\ci\bs
i_{\ti V,V}^{-1},$ where $\bs i_{\ti V,V}=\bS_{\smash{\id_{\ti
V},\id_{\ti E}}}:\bS_{\smash{\ti V,\ti E,\ti s}}\ra \bS_{V,E,s}$ is
a $1$-isomorphism, and\/~$\bS_{\smash{f,\hat f}}:\bS_{\smash{\ti
V,\ti E,\ti s}}\ra \bS_{W,F,t}$.
\label{sd4thm2}
\end{thm}

These results give a good differential-geometric picture of
d-manifolds and their 1- and 2-morphisms near a point. The $O(s)$
and $O(s^2)$ notation helps keep track of what information from
$V,E,s$ and $f,\hat f$ and $\La$ is remembered and what forgotten by
the d-manifolds $\bS_{V,E,s}$, 1-morphisms $\bS_{\smash{f,\hat f}}$
and 2-morphisms~$S_\La$.

\subsection{The 2-category of virtual vector bundles}
\label{sd43}

In our theory of derived differential geometry, it is a general
principle that categories in classical differential geometry should
often be replaced by 2-categories, and classical concepts be
replaced by 2-categorical analogues.

In classical differential geometry, if $X$ is a manifold, the vector
bundles $E\ra X$ and their morphisms form a category $\vect(X)$. The
cotangent bundle $T^*X$ is an important example of a vector bundle.
If $f:X\ra Y$ is smooth then pullback $f^*:\vect(Y)\ra\vect(X)$ is a
functor. There is a natural morphism $\d f^*:f^*(T^*Y)\ra T^*X$. We
now explain 2-categorical analogues of all this for d-manifolds,
following~\cite[\S 3.1--\S 3.2]{Joyc5}.

\begin{dfn} Let $\uX$ be a $C^\iy$-scheme, which will
usually be the $C^\iy$-scheme underlying a d-manifold $\bX$. We will
define a 2-category $\vqcoh(\uX)$ of {\it virtual quasicoherent
sheaves\/} on $\uX$. {\it Objects\/} of $\vqcoh(\uX)$ are morphisms
$\phi:\cE^1\ra\cE^2$ in $\qcoh(\uX)$, which we also may write as
$(\cE^1,\cE^2,\phi)$ or $(\cE^\bu,\phi)$. Given objects
$\phi:\cE^1\ra\cE^2$ and $\psi:\cF^1\ra\cF^2$, a 1-{\it morphism\/}
$(f^1,f^2):(\cE^\bu,\phi)\ra(\cF^\bu,\psi)$ is a pair of morphisms
$f^1:\cE^1\ra\cF^1$, $f^2:\cE^2\ra\cF^2$ in $\qcoh(\uX)$ such that
$\psi\ci f^1=f^2\ci\phi$. We write $f^\bu$ for~$(f^1,f^2)$.

The {\it identity\/ $1$-morphism\/} of $(\cE^\bu,\phi)$ is
$(\id_{\cE^1},\id_{\cE^2})$. The {\it composition\/} of 1-morphisms
$f^\bu:(\cE^\bu,\phi)\ra(\cF^\bu,\psi)$ and
$g^\bu:(\cF^\bu,\psi)\ra(\cG^\bu,\xi)$ is $g^\bu\ci f^\bu=(g^1\ci
f^1,g^2\ci f^2):(\cE^\bu,\phi)\ra(\cG^\bu,\xi)$.

Given $f^\bu,g^\bu:(\cE^\bu,\phi)\ra(\cF^\bu,\psi)$, a 2-{\it
morphism\/} $\eta:f^\bu\Ra g^\bu$ is a morphism $\eta:\cE^2\ra\cF^1$
in $\qcoh(\uX)$ such that $g^1=f^1+\eta\ci\phi$ and
$g^2=f^2+\psi\ci\eta$. The {\it identity\/ $2$-morphism\/} for
$f^\bu$ is $\id_{f^\bu}=0$. If $f^\bu,g^\bu,h^\bu:(\cE^\bu,\phi)\ra
(\cF^\bu,\psi)$ are 1-morphisms and $\eta:f^\bu\Ra g^\bu$,
$\ze:g^\bu\Ra h^\bu$ are 2-morphisms, the {\it vertical composition
of\/ $2$-morphisms\/} $\ze\od\eta:f^\bu\Ra h^\bu$ is
$\ze\od\eta=\ze+\eta$. If $f^\bu,\ti
f{}^\bu:(\cE^\bu,\phi)\ra(\cF^\bu, \psi)$ and $g^\bu,\ti
g{}^\bu:(\cF^\bu,\psi)\ra(\cG^\bu,\xi)$ are 1-morphisms and
$\eta:f^\bu\Ra\ti f{}^\bu$, $\ze:g^\bu\Ra\ti g{}^\bu$ are
2-morphisms, the {\it horizontal composition of\/ $2$-morphisms\/}
$\ze*\eta:g^\bu\ci f^\bu\Ra\ti g{}^\bu\ci\ti f{}^\bu$ is
$\ze*\eta=g^1\ci\eta+\ze\ci f^2+\ze\ci\psi\ci\eta$. This defines a
strict 2-category $\vqcoh(\uX)$, the obvious 2-category of 2-term
complexes in~$\qcoh(\uX)$.

If $\uU\subseteq\uX$ is an open $C^\iy$-subscheme then restriction
from $\uX$ to $\uU$ defines a strict 2-functor
$\vert_\uU:\vqcoh(\uX)\ra\vqcoh(\uU)$. An object $(\cE^\bu,\phi)$ in
$\vqcoh(\uX)$ is called a {\it virtual vector bundle of rank\/}
$d\in\Z$ if $\uX$ may be covered by open $\uU\subseteq\uX$ such that
$(\cE^\bu,\phi)\vert_\uU$ is equivalent in $\vqcoh(\uU)$ to some
$(\cF^\bu,\psi)$ for $\cF^1,\cF^2$ vector bundles on $\uU$ with
$\rank\cF^2-\rank\cF^1=d$. We write $\rank(\cE^\bu,\phi)=d$. If
$\uX\ne\es$ then $\rank(\cE^\bu,\phi)$ depends only on
$\cE^1,\cE^2,\phi$, so it is well-defined. Write $\vvect(\uX)$ for
the full 2-subcategory of virtual vector bundles in~$\vqcoh(\uX)$.

If $\uf:\uX\ra\uY$ is a $C^\iy$-scheme morphism then pullback gives
a strict 2-functor $\uf^*:\vqcoh(\uY)\ra\vqcoh(\uX)$, which maps
$\vvect(\uY)\ra\vvect(\uX)$.
\label{sd4def2}
\end{dfn}

We apply these ideas to d-spaces.

\begin{dfn} Let $\bX=(\uX,\OXp,\EX,\im_X,\jm_X)$ be a d-space. Define
the {\it virtual cotangent sheaf\/} $T^*\bX$ of $\bX$ to be the
morphism $\phi_X:\EX\ra\FX$ in $\qcoh(\uX)$ from Definition
\ref{sd3def1}, regarded as a virtual quasicoherent sheaf on~$\uX$.

Let $\bs f\!=\!(\uf,f',f''):\bX\ra\bY$ be a 1-morphism in $\dSpa$.
Then $T^*\bX\!=\!(\EX,\FX,\phi_X)$ and
$\uf^*(T^*\bY)\!=\!\bigl(\uf^*(\EY), \uf^*(\FY),\uf^*(\phi_Y)\bigr)$
are virtual quasicoherent sheaves on $\uX$, and $\Om_{\bs
f}:=\!(f'',f^2)$ is a 1-morphism $\uf^*(T^*\bY)\!\ra\! T^*\bX$ in
$\vqcoh(\uX)$, as \eq{sd3eq1} commutes.

Let $\bs f,\bs g:\bX\ra\bY$ be 1-morphisms in $\dSpa$, and $\eta:\bs
f\Ra\bs g$ a 2-morphism. Then $\eta:\uf^*(\FY)\ra\EX$ with
$g''=f''+\eta\ci\uf^*(\phi_Y)$ and $g^2=f^2+\phi_X\ci\eta$, as in
\eq{sd3eq2}. It follows that $\eta$ is a 2-morphism $\Om_{\bs
f}\Ra\Om_{\bs g}$ in $\vqcoh(\uX)$. Thus, objects, 1-morphisms and
2-morphisms in $\dSpa$ lift to objects, 1-morphisms and 2-morphisms
in $\vqcoh(\uX)$.
\label{sd4def3}
\end{dfn}

The next proposition justifies the definition of virtual vector
bundle. Because of part (b), if $\bW$ is a d-manifold we call
$T^*\bW$ the {\it virtual cotangent bundle\/} of $\bW$, rather than
the virtual cotangent sheaf.

\begin{prop}{\bf(a)} Let\/ $V$ be a manifold, $E\ra V$ a vector
bundle, and\/ $s\in C^\iy(E)$. Then Example\/ {\rm\ref{sd4ex1}}
defines a principal d-manifold\/ $\bS_{V,E,s}$. Its cotangent bundle
$T^*\bS_{V,E,s}$ is a virtual vector bundle on $\uS_{V,E,s}$ of
rank\/~$\dim V-\rank E$.
\smallskip

\noindent{\bf(b)} Let\/ $\bW$ be a d-manifold. Then $T^*\bW$ is a
virtual vector bundle on $\uW$ of rank\/ $\vdim\bW$. Hence if\/
$\bW\ne\bs\es$ then $\vdim\bW$ is well-defined.
\label{sd4prop2}
\end{prop}

The virtual cotangent bundle $T^*\bX$ of a d-manifold $\bX$ contains
only a fraction of the information in $\bX=(\uX,\OXp,\EX,\im_X,
\jm_X)$, but many interesting properties of d-manifolds $\bX$ and
1-morphisms $\bs f:\bX\ra\bY$ can be expressed solely in terms of
virtual cotangent bundles $T^*\bX,T^*\bY$ and 1-morphisms $\Om_{\bs
f}:\uf^*(T^*\bY)\ra T^*\bX$. Here is an example of this.

\begin{dfn} Let $\uX$ be a $C^\iy$-scheme. We say that a
virtual vector bundle $(\cE^1,\cE^2,\phi)$ on $\uX$ {\it is a vector
bundle\/} if it is equivalent in $\vvect(\uX)$ to $(0,\cE,0)$ for
some vector bundle $\cE$ on $\uX$. One can show $(\cE^1,\cE^2,\phi)$
is a vector bundle if and only if $\phi$ has a left inverse
in~$\qcoh(\uX)$.
\label{sd4def4}
\end{dfn}

\begin{prop} Let\/ $\bX$ be a d-manifold. Then $\bX$ is a manifold
(that is, $\bX\in\hMan$) if and only if\/ $T^*\bX$ is a vector
bundle, or equivalently, if\/ $\phi_X:\EX\ra\FX$ has a left inverse
in\/~$\qcoh(\uX)$.
\label{sd4prop3}
\end{prop}

\subsection{Equivalences in $\dMan$, and gluing by equivalences}
\label{sd44}

Equivalences in a 2-category are defined in Appendix \ref{sdA}.
Equivalences in $\dMan$ are the best derived analogue of
isomorphisms in $\Man$, that is, of diffeomorphisms of manifolds. A
smooth map of manifolds $f:X\ra Y$ is called {\it \'etale\/} if it
is a local diffeomorphism. Here is the derived analogue.

\begin{dfn} Let $\bs f:\bX\ra\bY$ be a 1-morphism in $\dMan$. We
call $\bs f$ {\it \'etale\/} if it is a {\it local equivalence},
that is, if for each $x\in\bX$ there exist open
$x\in\bU\subseteq\bX$ and $\bs f(x)\in\bV\subseteq\bY$ such that
$\bs f(\bU)=\bV$ and $\bs f\vert_\bU:\bU\ra\bV$ is an equivalence.
\label{sd4def5}
\end{dfn}

If $f:X\ra Y$ is a smooth map of manifolds, then $f$ is \'etale if
and only if $\d f^*:f^*(T^*Y)\ra T^*X$ is an isomorphism of vector
bundles. (The analogue is false for schemes.) In \cite[\S
3.5]{Joyc5} we prove a version of this for d-manifolds:

\begin{thm} Suppose $\bs f:\bX\ra\bY$ is a $1$-morphism of
d-manifolds. Then the following are equivalent:
\begin{itemize}
\setlength{\itemsep}{0pt}
\setlength{\parsep}{0pt}
\item[{\rm(i)}] $\bs f$ is \'etale;
\item[{\rm(ii)}] $\Om_{\bs f}:\uf^*(T^*\bY)\ra T^*\bX$ is an
equivalence in $\vqcoh(\uX);$ and
\item[{\rm(iii)}] the following is a split short
exact sequence in\/~$\qcoh(\uX)\!:$
\begin{equation*}
\smash{\xymatrix@C=25pt{ 0 \ar[r] & \uf^*(\EY) \ar[rr]^(0.45){f''
\op -\uf^*(\phi_Y)} && \EX\op \uf^*(\FY) \ar[rr]^(0.6){\phi_X\op
f^2} && \FX \ar[r] & 0.}}
\end{equation*}
\end{itemize}
If in addition $f:X\ra Y$ is a bijection, then $\bs f$ is an
equivalence in\/~$\dMan$.
\label{sd4thm3}
\end{thm}

The analogue of Theorem \ref{sd4thm3} for d-spaces is false. When
$\bs f:\bX\ra\bY$ is a `standard model' 1-morphism
$\bS_{\smash{f,\hat f}}:\bS_{V,E,s}\ra\bS_{W,F,t}$, as in
\S\ref{sd42}, we can express the conditions for $\bS_{\smash{f,\hat
f}}$ to be \'etale or an equivalence in terms of~$f,\hat f$.

\begin{thm} Let\/ $V,W$ be manifolds, $E\ra V,$ $F\ra W$ be vector
bundles, $s\in C^\iy(E),$ $t\in C^\iy(F),$ $f:V\ra W$ be smooth,
and\/ $\hat f:E\ra f^*(F)$ be a morphism of vector bundles on $V$
with\/ $\hat f\ci s= f^*(t)+O(s^2)$. Then Example\/
{\rm\ref{sd4ex2}} defines a $1$-morphism\/ $\bS_{\smash{f,\hat
f}}:\bS_{V,E,s}\ra\bS_{W,F,t}$ in $\dMan$. This $\bS_{\smash{f,\hat
f}}$ is \'etale if and only if for each\/ $v\in V$ with\/ $s(v)=0$
and\/ $w=f(v)\in W,$ the following sequence of vector spaces is
exact:
\begin{equation*}
\smash{\xymatrix@C=18pt{ 0 \ar[r] & T_vV \ar[rrr]^(0.42){\d s(v)\op
\,\d f(v)} &&& E_v\op T_wW \ar[rrr]^(0.57){\hat f(v)\op\, -\d t(w)}
&&& F_w \ar[r] & 0.}}
\end{equation*}
Also $\bS_{\smash{f,\hat f}}$ is an equivalence if and only if in
addition\/ $f\vert_{s^{-1}(0)}:s^{-1}(0)\!\ra\! t^{-1}(0)$ is a
bijection, where $s^{-1}(0)\!=\!\{v\in V:s(v)\!=\!0\},$
$t^{-1}(0)\!=\!\{w\in W:t(w)\!=\!0\}$.
\label{sd4thm4}
\end{thm}

Section \ref{sd32} discussed gluing d-spaces by equivalences on open
d-subspaces. It generalizes immediately to d-manifolds: if in
Theorem \ref{sd3thm3} we fix $n\in\Z$ and take the initial d-spaces
$\bX_i$ to be d-manifolds with $\vdim\bX_i=n$, then the glued
d-space $\bY$ is also a d-manifold with $\vdim\bY=n$.

Here is an analogue of Theorem \ref{sd3thm3}, taken from \cite[\S
3.6]{Joyc5}, in which we take the d-spaces $\bX_i$ to be `standard
model' d-manifolds $\bS_{V_i,E_i,s_i}$, and the 1-morphisms $\bs
e_{ij}$ to be `standard model' 1-morphisms $\bS_{e_{ij},\hat
e_{ij}}$. We also use Theorem \ref{sd4thm4} in (iii) to characterize
when $\bs e_{ij}$ is an equivalence.

\begin{thm} Suppose we are given the following data:
\begin{itemize}
\setlength{\itemsep}{0pt}
\setlength{\parsep}{0pt}
\item[{\rm(a)}] an integer $n;$
\item[{\rm(b)}] a Hausdorff, second countable topological space $X;$
\item[{\rm(c)}] an indexing set\/ $I,$ and a total order $<$ on $I;$
\item[{\rm(d)}] for each\/ $i$ in $I,$ a manifold\/ $V_i,$ a vector
bundle $E_i\ra V_i$ with\/ $\dim V_i-\rank E_i=n,$ a smooth
section $s_i:V_i\ra E_i,$ and a homeomorphism $\psi_i:X_i\ra\hat
X_i,$ where $X_i=\{v_i\in V_i:s_i(v_i)=0\}$ and\/ $\hat
X_i\subseteq X$ is open; and
\item[{\rm(e)}] for all\/ $i<j$ in $I,$ an open submanifold
$V_{ij}\subseteq V_i,$ a smooth map $e_{ij}:V_{ij}\ra V_j,$ and
a morphism of vector bundles $\hat e_{ij}:E_i\vert_{V_{ij}}\ra
e_{ij}^*(E_j)$.
\end{itemize}
Using notation $O(s_i),O(s_i^2)$ as in\/ {\rm\S\ref{sd42},} let this
data satisfy the conditions:
\begin{itemize}
\setlength{\itemsep}{0pt}
\setlength{\parsep}{0pt}
\item[{\rm(i)}] $X=\bigcup_{i\in I}\hat X_i;$
\item[{\rm(ii)}] if\/ $i<j$ in $I$ then $\hat e_{ij}\ci
s_i\vert_{V_{ij}}= e_{ij}^*(s_j)+O(s_i^2),$ $\psi_i(X_i\cap
V_{ij})=\hat X_i\cap\hat X_j,$ and\/ $\psi_i\vert_{X_i\cap
V_{ij}}=\psi_j\ci e_{ij}\vert_{X_i\cap V_{ij}},$ and if\/
$v_i\in V_{ij}$ with\/ $s_i(v_i)=0$ and\/ $v_j=e_{ij}(v_i)$ then
the following is exact:
\begin{equation*}
\smash{\xymatrix@C=18pt{ 0 \ar[r] & T_{v_i}V_i \ar[rrr]^(0.42){\d
s_i(v_i)\op \,\d e_{ij}(v_i)} &&& E_i\vert_{v_i}\!\op\! T_{v_j}V_j
\ar[rrr]^(0.57){\hat e_{ij}(v_i)\op\, -\d s_j(v_j)} &&&
E_j\vert_{v_j} \ar[r] & 0;}}
\end{equation*}
\item[{\rm(iii)}] if\/ $i<j<k$ in $I$ then
\begin{align*}
e_{ik}\vert_{V_{ij}\cap V_{ik}}&= e_{jk}\ci
e_{ij}\vert_{V_{ij}\cap V_{ik}}+O(s_i^2)\qquad\text{and}\\
\hat e_{ik}\vert_{V_{ij}\cap V_{ik}}&= e_{ij}\vert_{V_{ij}\cap
V_{ik}}^*(\hat e_{jk})\ci \hat e_{ij}\vert_{V_{ij}\cap
V_{ik}}+O(s_i).
\end{align*}
\end{itemize}

Then there exist a d-manifold\/ $\bX$ with\/ $\vdim\bX=n$ and
underlying topological space $X,$ and a $1$-morphism
$\bs\psi_i:\bS_{V_i,E_i,s_i}\ra\bX$ with underlying continuous map
$\psi_i$ which is an equivalence with the open d-submanifold\/
$\bs{\hat X}_i\subseteq\bX$ corresponding to $\hat X_i\subseteq X$
for all\/ $i\in I,$ such that for all\/ $i<j$ in $I$ there exists a
$2$-morphism\/ $\eta_{ij}:\bs\psi_j\ci\bS_{e_{ij},\hat
e_{ij}}\Ra\bs\psi_i\ci\bs i_{V_{ij},V_i},$ where $\bS_{e_{ij},\hat
e_{ij}}:\bS_{V_{ij},E_i\vert_{V_{ij}},s_i\vert_{V_{ij}}}\ra
\bS_{V_j,E_j,s_j}$ and\/ $\bs i_{V_{ij},V_i}:\bS_{V_{ij},
E_i\vert_{V_{ij}},s_i\vert_{V_{ij}}}\ra \bS_{V_i,E_i,s_i}$. This
d-manifold\/ $\bX$ is unique up to equivalence in~$\dMan$.

Suppose also that\/ $Y$ is a manifold, and\/ $g_i:V_i\ra Y$ are
smooth maps for all\/ $i\in I,$ and\/ $g_j\ci
e_{ij}=g_i\vert_{V_{ij}}+O(s_i)$ for all\/ $i<j$ in $I$. Then there
exist a $1$-morphism $\bs h:\bX\ra\bY$ unique up to $2$-isomorphism,
where $\bY=F_\Man^\dMan(Y)=\bS_{Y,0,0},$ and\/ $2$-morphisms
$\ze_i:\bs h\ci\bs\psi_i\Ra\bS_{g_i,0}$ for all\/ $i\in I$. Here
$\bS_{Y,0,0}$ is from Example\/ {\rm\ref{sd4ex1}} with vector bundle
$E$ and section $s$ both zero, and\/
$\bS_{g_i,0}:\bS_{V_i,E_i,s_i}\ra \bS_{Y,0,0}=\bY$ is from Example
{\rm\ref{sd4ex2}} with\/~$\hat g_i=0$.
\label{sd4thm5}
\end{thm}

The hypotheses of Theorem \ref{sd4thm5} are similar to the notion of
{\it good coordinate system\/} in the theory of Kuranishi spaces of
Fukaya and Ono \cite[Def.~6.1]{FuOn}. The importance of Theorem
\ref{sd4thm5} is that all the ingredients are described wholly in
differential-geometric or topological terms. So we can use the
theorem as a tool to prove the existence of d-manifold structures on
spaces coming from other areas of geometry, for instance, on moduli
spaces.

\subsection{Submersions, immersions and embeddings}
\label{sd45}

Let $f:X\ra Y$ be a smooth map of manifolds. Then $\d
f^*:f^*(T^*Y)\ra T^*X$ is a morphism of vector bundles on $X$, and
$f$ is a {\it submersion\/} if $\d f^*$ is injective, and $f$ is an
{\it immersion\/} if $\d f^*$ is surjective. Here the appropriate
notions of injective and surjective for morphisms of vector bundles
are stronger than the corresponding notions for sheaves: $\d f^*$ is
{\it injective\/} if it has a left inverse, and {\it surjective\/}
if it has a right inverse.

In a similar way, if $\bs f:\bX\ra\bY$ is a 1-morphism of
d-manifolds, we would like to define $\bs f$ to be a submersion or
immersion if the 1-morphism $\Om_{\bs f}:\uf^*(T^*\bY)\ra T^*\bX$ in
$\vvect(\uX)$ is injective or surjective in some suitable sense. It
turns out that there are two different notions of injective and
surjective 1-morphisms in the 2-category $\vvect(\uX)$, a weak and a
strong:

\begin{dfn} Let $\uX$ be a $C^\iy$-scheme, $(\cE^1,\cE^2,\phi)$ and
$(\cF^1,\cF^2,\psi)$ be virtual vector bundles on $\uX$, and
$(f^1,f^2):(\cE^\bu,\phi)\ra(\cF^\bu,\psi)$ be a 1-morphism in
$\vvect(\uX)$. Then we have a complex in~$\qcoh(\uX)$:
\e
\xymatrix@C=25pt{ 0 \ar[r] & \cE^1 \ar@<.5ex>[rr]^(0.45){f^1 \op
-\phi} && \cF^1\op \cE^2 \ar@<.5ex>[rr]^(0.6){\psi\op f^2}
\ar@<.5ex>@{.>}[ll]^(0.55)\ga && \cF^2 \ar@<.5ex>@{.>}[ll]^(0.4)\de
\ar[r] & 0.}
\label{sd4eq3}
\e
One can show that $f^\bu$ is an equivalence in $\vvect(\uX)$ if and
only if \eq{sd4eq3} is a {\it split short exact sequence\/} in
$\qcoh(\uX)$. That is, $f^\bu$ is an equivalence if and only if
there exist morphisms $\ga,\de$ as shown in \eq{sd4eq3} satisfying
the conditions:
\e
\begin{aligned}
\ga\ci\de&=0,& \ga\ci(f^1 \op -\phi)&=\id_{\cE^1},\\
(f^1 \op -\phi)\ci\ga+\de\ci(\psi\op f^2)&=\id_{\cF^1\op \cE^2}, &
(\psi\op f^2)\ci\de&=\id_{\cF^2}.
\end{aligned}
\label{sd4eq4}
\e

Our notions of $f^\bu$ injective or surjective impose some but not
all of \eq{sd4eq4}:
\begin{itemize}
\setlength{\itemsep}{0pt}
\setlength{\parsep}{0pt}
\item[(a)] We call $f^\bu$ {\it weakly injective\/} if there
exists $\ga:\cF^1\op\cE^2\ra\cE^1$ in $\qcoh(\uX)$
with~$\ga\ci(f^1 \op -\phi)=\id_{\cE^1}$.
\item[(b)] We call $f^\bu$ {\it injective\/} if there exist
$\ga:\cF^1\op\cE^2\ra\cE^1$ and $\de:\cF^2\ra\cF^1\op\cE^2$ with
$\ga\ci\de=0$, $\ga\ci(f^1 \op -\phi)= \id_{\cE^1}$ and~$(f^1
\op -\phi)\ci\ga+\de\ci(\psi\op f^2)=\id_{\cF^1\op\cE^2}$.
\item[(c)] We call $f^\bu$ {\it weakly surjective\/} if there exists
$\de:\cF^2\ra\cF^1\op\cE^2$ in $\qcoh(\uX)$ with~$(\psi\op
f^2)\ci\de=\id_{\cF^2}$.
\item[(d)] We call $f^\bu$ {\it surjective\/} if there exist
$\ga:\cF^1\op\cE^2\ra\cE^1$ and $\de:\cF^2\ra\cF^1\op\cE^2$ with
$\ga\ci\de=0$, $\ga\ci(f^1 \op -\phi)= \id_{\cE^1}$
and~$(\psi\op f^2)\ci\de=\id_{\cF^2}$.
\end{itemize}
\label{sd4def6}
\end{dfn}

Using these we define weak and strong forms of submersions,
immersions, and embeddings for d-manifolds.

\begin{dfn} Let $\bs f:\bX\ra\bY$ be a 1-morphism of d-manifolds.
Definition \ref{sd4def3} defines a 1-morphism $\Om_{\bs
f}:\uf^*(T^*\bY)\ra T^*\bX$ in $\vvect(\uX)$. Then:
\begin{itemize}
\setlength{\itemsep}{0pt}
\setlength{\parsep}{0pt}
\item[(a)] We call $\bs f$ a {\it w-submersion\/} if $\Om_{\bs
f}$ is weakly injective.
\item[(b)] We call $\bs f$ a {\it submersion\/} if $\Om_{\bs
f}$ is injective.
\item[(c)] We call $\bs f$ a {\it w-immersion\/} if $\Om_{\bs
f}$ is weakly surjective.
\item[(d)] We call $\bs f$ an {\it immersion\/} if $\Om_{\bs
f}$ is surjective.
\item[(e)] We call $\bs f$ a {\it w-embedding\/} if it is a
w-immersion and $f:X\ra f(X)$ is a homeomorphism, so in
particular $f$ is injective.
\item[(f)] We call $\bs f$ an {\it embedding\/} if it is an
immersion and $f$ is a homeomorphism with its image.
\end{itemize}
Here w-submersion is short for {\it weak submersion}, etc.
Conditions (a)--(d) all concern the existence of morphisms $\ga,\de$
in the next equation satisfying identities.
\begin{equation*}
\xymatrix@C=25pt{ 0 \ar[r] & \uf^*(\EY) \ar@<.5ex>[rr]^(0.45){f''
\op -\uf^*(\phi_Y)} && \EX\op \uf^*(\FY)
\ar@<.5ex>@{.>}[ll]^(0.55)\ga \ar@<.5ex>[rr]^(0.6){\phi_X\op f^2} &&
\FX \ar@<.5ex>@{.>}[ll]^(0.4)\de \ar[r] & 0.}
\end{equation*}

Parts (c)--(f) enable us to define {\it d-submanifolds\/} of
d-manifolds. {\it Open d-submanifolds\/} are open d-subspaces of a
d-manifold. More generally, we call $\bs i:\bX\ra\bY$ a {\it
w-immersed}, or {\it immersed}, or {\it w-embedded}, or {\it
embedded d-submanifold}, of $\bY$, if $\bX,\bY$ are d-manifolds and
$\bs i$ is a w-immersion, immersion, w-embedding, or embedding,
respectively.
\label{sd4def7}
\end{dfn}

Here are some properties of these, taken from~\cite[\S 4.1--\S
4.2]{Joyc5}:

\begin{thm}{\bf(i)} Any equivalence of d-manifolds is a
w-submersion, submersion, w-immersion, immersion, w-embedding and
embedding.
\smallskip

\noindent{\bf(ii)} If\/ $\bs f,\bs g:\bX\ra\bY$ are $2$-isomorphic
$1$-morphisms of d-manifolds then $\bs f$ is a w-submersion,
submersion, \ldots, embedding, if and only if\/ $\bs g$ is.
\smallskip

\noindent{\bf(iii)} Compositions of w-submersions, submersions,
w-immersions, immersions, w-embeddings, and embeddings are
$1$-morphisms of the same kind.
\smallskip

\noindent{\bf(iv)} The conditions that a $1$-morphism of
d-manifolds\/ $\bs f:\bX\ra\bY$ is a w-submersion, submersion,
w-immersion or immersion are local in $\bX$ and\/ $\bY$. That is,
for each\/ $x\in\bX$ with\/ $\bs f(x)=y\in\bY,$ it suffices to check
the conditions for $\bs f\vert_\bU:\bU\ra\bV$ with $\bV$ an open
neighbourhood of\/ $y$ in $\bY,$ and\/ $\bU$ an open neighbourhood
of\/ $x$ in\/~${\bs f}^{-1}(\bV)\!\subseteq\!\bX$.
\smallskip

\noindent{\bf(v)} Let\/ $\bs f:\bX\ra\bY$ be a submersion of
d-manifolds. Then $\vdim\bX\ge\vdim\bY,$ and if\/
$\vdim\bX=\vdim\bY$ then\/ $\bs f$ is \'etale.
\smallskip

\noindent{\bf(vi)} Let\/ $\bs f:\bX\ra\bY$ be an immersion of
d-manifolds. Then $\vdim\bX\le\vdim\bY,$ and if\/
$\vdim\bX=\vdim\bY$ then\/ $\bs f$ is \'etale.
\smallskip

\noindent{\bf(vii)} Let\/ $f:X\ra Y$ be a smooth map of manifolds,
and\/ $\bs f=F_\Man^\dMan(f)$. Then $\bs f$ is a submersion,
immersion, or embedding in $\dMan$ if and only if\/ $f$ is a
submersion, immersion, or embedding in $\Man,$ respectively. Also
$\bs f$ is a w-immersion or w-embedding if and only if\/ $f$ is an
immersion or embedding.
\smallskip

\noindent{\bf(viii)} Let\/ $\bs f:\bX\ra\bY$ be a $1$-morphism of
d-manifolds, with\/ $\bY$ a manifold. Then $\bs f$ is a
w-submersion.
\smallskip

\noindent{\bf(ix)} Let\/ $\bX,\bY$ be d-manifolds, with\/ $\bY$ a
manifold. Then $\bs\pi_\bX:\bX\t\bY\ra\bX$ is a submersion.
\smallskip

\noindent{\bf(x)} Let\/ $\bs f:\bX\ra\bY$ be a submersion of
d-manifolds, and\/ $x\in\bX$ with\/ $\bs f(x)=y\in\bY$. Then there
exist open $x\in\bU\subseteq\bX$ and\/ $y\in\bV\subseteq\bY$ with\/
$\bs f(\bU)=\bV,$ a manifold\/ $\bZ,$ and an equivalence $\bs
i:\bU\ra\bV\t\bZ,$ such that\/ $\bs f\vert_\bU:\bU\ra\bV$ is
$2$-isomorphic to $\bs\pi_\bV\ci\bs i,$ where
$\bs\pi_\bV:\bV\t\bZ\ra\bV$ is the projection.
\smallskip

\noindent{\bf(xi)} Let\/ $\bs f:\bX\ra\bY$ be a submersion of
d-manifolds with\/ $\bY$ a manifold. Then\/ $\bX$ is a manifold.
\label{sd4thm6}
\end{thm}

\subsection{D-transversality and fibre products}
\label{sd46}

From \S\ref{sd33}, if $\bs g:\bX\ra\bZ$ and $\bs h:\bY\ra\bZ$ are
1-morphisms of d-manifolds then a fibre product $\bW=\bX_{\bs
g,\bZ,\bs h}\bY$ exists in $\dSpa$, and is unique up to equivalence.
We want to know whether $\bW$ is a d-manifold. We will define when
$\bs g,\bs h$ are {\it d-transverse}, which is a sufficient
condition for $\bW$ to be a d-manifold.

Recall that if $g:X\ra Z$, $h:Y\ra Z$ are smooth maps of manifolds,
then a fibre product $W=X\t_{g,Z,h}Y$ in $\Man$ exists if $g,h$ are
{\it transverse}, that is, if $T_zZ=\d g\vert_x(T_xX)+\d
h\vert_y(T_yY)$ for all $x\in X$ and $y\in Y$ with $g(x)=h(y)=z\in
Z$. Equivalently, $\d g\vert_x^*\op\d h\vert_y^*:T_zZ^*\ra T_x^*X\op
T_y^*Y$ should be injective. Writing $W=X\t_ZY$ for the topological
fibre product and $e:W\ra X$, $f:W\ra Y$ for the projections, with
$g\ci e=h\ci f$, we see that $g,h$ are transverse if and only if
\e
e^*(\d g^*)\op f^*(\d h^*):(g\ci e)^*(T^*Z)\ra e^*(T^*X)\op
f^*(T^*Y)
\label{sd4eq5}
\e
is an injective morphism of vector bundles on the topological space
$W$, that is, it has a left inverse. The condition that \eq{sd4eq6}
has a left inverse is an analogue of this, but on (dual) obstruction
rather than cotangent bundles.

\begin{dfn} Let $\bX,\bY,\bZ$ be d-manifolds and $\bs g:\bX\ra\bZ,$
$\bs h:\bY\ra\bZ$ be 1-morphisms. Let $\uW=\uX\t_{\ug,\uZ,\uh}\uY$
be the $C^\iy$-scheme fibre product, and write $\ue:\uW\ra\uX$,
$\uf:\uW\ra\uY$ for the projections. Consider the morphism
\e
\begin{split}
\al=\ue^*(g'') \op -\uf^*(h'') &\op (\ug\ci\ue)^*(\phi_Z):
(\ug\ci\ue)^*(\EZ)\longra\\
&\qquad\qquad\ue^*(\EX)\op\uf^*(\EY)\op(\ug\ci\ue)^*(\FZ)
\end{split}
\label{sd4eq6}
\e
in $\qcoh(\uW)$. We call $\bs g,\bs h$ {\it d-transverse\/} if $\al$
has a left inverse.

In the notation of \S\ref{sd43} and \S\ref{sd45}, we have
1-morphisms $\Om_{\bs g}:\ug^*(T^*\bZ)\ra T^*\bX$ in $\vvect(\uX)$
and $\Om_{\bs h}:\uh^*(T^*\bZ)\ra T^*\bY$ in $\vvect(\uY)$. Pulling
these back to $\vvect(\uW)$ using $\ue^*,\uf^*$ we form the
1-morphism in~$\vvect(\uW)$:
\e
\ue^*(\Om_{\bs g}) \op \uf^*(\Om_{\bs h}): (\ug\ci\ue)^*(T^*\bZ)
\longra\ue^*(T^*\bX)\op\uf^*(T^*\bY).
\label{sd4eq7}
\e
For \eq{sd4eq6} to have a left inverse is equivalent to \eq{sd4eq7}
being weakly injective, as in Definition \ref{sd4def6}. This is the
d-manifold analogue of \eq{sd4eq5} being injective.
\label{sd4def8}
\end{dfn}

Here are the main results of~\cite[\S 4.3]{Joyc5}:

\begin{thm} Suppose\/ $\bX,\bY,\bZ$ are d-manifolds and\/ $\bs
g:\bX\ra\bZ,$ $\bs h:\bY\ra\bZ$ are d-transverse $1$-morphisms, and
let\/ $\bW=\bX\t_{\bs g,\bZ,\bs h}\bY$ be the d-space fibre product.
Then $\bW$ is a d-manifold, with
\e
\vdim\bW=\vdim\bX+\vdim\bY-\vdim\bZ.
\label{sd4eq8}
\e
\label{sd4thm7}
\end{thm}

\begin{thm} Suppose\/ $\bs g:\bX\ra\bZ,$ $\bs h:\bY\ra\bZ$ are
$1$-morphisms of d-manifolds. The following are sufficient
conditions for $\bs g,\bs h$ to be d-transverse, so that\/
$\bW=\bX\t_{\bs g,\bZ,\bs h}\bY$ is a d-manifold of virtual
dimension\/ {\rm\eq{sd4eq8}:}
\begin{itemize}
\setlength{\itemsep}{0pt}
\setlength{\parsep}{0pt}
\item[{\bf(a)}] $\bZ$ is a manifold, that is, $\bZ\in\hMan;$ or
\item[{\bf(b)}] $\bs g$ or $\bs h$ is a w-submersion.
\end{itemize}
\label{sd4thm8}
\end{thm}

The point here is that roughly speaking, $\bs g,\bs h$ are
d-transverse if they map the direct sum of the obstruction spaces of
$\bX,\bY$ surjectively onto the obstruction spaces of $\bZ$. If
$\bZ$ is a manifold its obstruction spaces are zero. If $\bs g$ is a
w-submersion it maps the obstruction spaces of $\bX$ surjectively
onto the obstruction spaces of $\bZ$. In both cases,
d-transversality follows. See \cite[Th.~8.15]{Spiv} for the analogue
of Theorem \ref{sd4thm8}(a) for Spivak's derived manifolds.

\begin{thm} Let\/ $\bX,\bZ$ be d-manifolds, $\bY$ a manifold, and\/
$\bs g:\bX\ra\bZ,$ $\bs h:\bY\ra\bZ$ be $1$-morphisms with\/ $\bs g$
a submersion. Then\/ $\bW=\bX\t_{\bs g,\bZ,\bs h}\bY$ is a manifold,
with\/~$\dim\bW=\vdim\bX+\dim\bY-\vdim\bZ$.
\label{sd4thm9}
\end{thm}

Theorem \ref{sd4thm9} shows that we may think of submersions as {\it
representable\/ $1$-morphisms\/} in $\dMan$. We can locally
characterize embeddings and immersions in $\dMan$ in terms of fibre
products with $\bR^{\bs n}$ in~$\dMan$.

\begin{thm}{\bf(a)} Let\/ $\bX$ be a d-manifold and\/ $\bs
g:\bX\ra\bR^{\bs n}$ a $1$-morphism in $\dMan$. Then the fibre
product\/ $\bW=\bX\t_{\bs g,\bR^{\bs n},\bs 0}\bs *$ exists in
$\dMan$ by Theorem\/ {\rm\ref{sd4thm8}(a),} and the projection $\bs
e:\bW\ra\bX$ is an embedding.
\smallskip

\noindent{\bf(b)} Suppose\/ $\bs f:\bX\ra\bY$ is an immersion of
d-manifolds, and\/ $x\in\bX$ with\/ $\bs f(x)=y\in\bY$. Then there
exist open d-submanifolds $x\in\bU\subseteq\bX$ and\/
$y\in\bV\subseteq\bY$ with\/ $\bs f(\bU)\subseteq\bV,$ and a
$1$-morphism $\bs g:\bV\ra\bR^{\bs n}$ with\/ $\bs g(y)=0,$ where
$n=\vdim\bY-\vdim\bX\ge 0,$ fitting into a $2$-Cartesian square
in~$\dMan:$
\begin{equation*}
\xymatrix@C=60pt@R=10pt{\bU \ar[d]^{\bs f\vert_\bU} \ar[r]
\drtwocell_{}\omit^{}\omit{^{}} & \bs{*} \ar[d]_{\bs 0} \\
\bV \ar[r]^(0.7){\bs g} & \bR^{\bs n}. }
\end{equation*}
If\/ $\bs f$ is an embedding we may take\/~$\bU=\bs f^{-1}(\bV)$.
\label{sd4thm10}
\end{thm}

\begin{rem} For the applications the author has in
mind, it will be crucial that if $\bs g:\bX\ra\bZ$ and $\bs
h:\bY\ra\bZ$ are 1-morphisms with $\bX,\bY$ d-manifolds and $\bZ$ a
manifold then $\bW=\bX\t_\bZ\bY$ is a d-manifold, with
$\vdim\bW=\vdim\bX+\vdim\bY-\dim\bZ$, as in Theorem
\ref{sd4thm8}(a). We will show by example, following Spivak
\cite[Prop.~1.7]{Spiv}, that if d-manifolds $\dMan$ were an ordinary
category containing manifolds as a full subcategory, then this would
be false.

Consider the fibre product $\bs{*}\t_{\bs 0,\bR,\bs 0}\bs{*}$ in
$\dMan$. If $\dMan$ were a category then as $\bs{*}$ is a terminal
object, the fibre product would be $\bs{*}$. But then
\begin{equation*}
\vdim(\bs{*}\t_{\bs 0,\bR,\bs 0}\bs{*})=\vdim\bs{*}=0\ne -1=\vdim\bs{*}+
\vdim\bs{*}-\vdim\bR,
\end{equation*}
so equation \eq{sd4eq8} and Theorem \ref{sd4thm8}(a) would be false.
Thus, if we want fibre products of d-manifolds over manifolds to be
well behaved, then $\dMan$ must be at least a 2-category. It could
be an $\iy$-category, as for Spivak's derived manifolds \cite{Spiv},
or some other kind of higher category. Making d-manifolds into a
2-category, as we have done, is the simplest of the available
options.
\label{sd4rem1}
\end{rem}

\subsection{Embedding d-manifolds into manifolds}
\label{sd47}

Let $V$ be a manifold, $E\ra V$ a vector bundle, and $s\in
C^\iy(E)$. Then Example \ref{sd4ex1} defines a `standard model'
principal d-manifold $\bS_{V,E,s}$. When $E$ and $s$ are zero, we
have $\bS_{V,0,0}=\bV=F_\Man^\dMan(V)$, so that $\bS_{V,0,0}$ is a
manifold. For general $V,E,s$, taking $f=\id_V:V\ra V$ and $\hat
f=0:0\ra E$ in Example \ref{sd4ex2} gives a `standard model'
1-morphism $\bS_{\id_V,0}:\bS_{V,E,s}\ra \bS_{V,0,0}=\bV$. One can
show $\bS_{\id_V,0}$ is an embedding, in the sense of Definition
\ref{sd4def7}. Any principal d-manifold $\bU$ is equivalent to some
$\bS_{V,E,s}$. Thus we deduce:

\begin{lem} Any principal d-manifold\/ $\bU$ admits an embedding\/
$\bs i:\bU\ra\bV$ into a manifold\/~$\bV$.
\label{sd4lem1}
\end{lem}

Theorem \ref{sd4thm13} below is a converse to this: if a d-manifold
$\bX$ can be embedded into a manifold $\bY$, then $\bX$ is
principal. So it will be useful to study embeddings of d-manifolds
into manifolds. The following facts are due to Whitney~\cite{Whitn}.

\begin{thm}{\bf(a)} Let\/ $X$ be an $m$-manifold and\/ $n\ge 2m$.
Then a generic smooth map $f:X\ra\R^n$ is an immersion.
\smallskip

\noindent{\bf(b)} Let\/ $X$ be an $m$-manifold and\/ $n\ge 2m+1$.
Then there exists an embedding $f:X\ra\R^n,$ and we can choose
such\/ $f$ with\/ $f(X)$ closed in\/ $\R^n$. Generic smooth maps\/
$f:X\ra\R^n$ are embeddings.
\label{sd4thm11}
\end{thm}

In \cite[\S 4.4]{Joyc5} we generalize Theorem \ref{sd4thm11} to
d-manifolds.

\begin{thm} Let\/ $\bX$ be a d-manifold. Then there exist
immersions and/or embeddings $\bs f:\bX\ra\bR^{\bs n}$ for some
$n\gg 0$ if and only if there is an upper bound for\/ $\dim
T^*_x\uX$ for all\/ $x\in\uX$. If there is such an upper bound, then
immersions $\bs f:\bX\ra\bR^{\bs n}$ exist provided\/ $n\ge 2\dim
T_x^*\uX$ for all\/ $x\in\uX,$ and embeddings $\bs f:\bX\ra\bR^{\bs
n}$ exist provided\/ $n\ge 2\dim T_x^*\uX+1$ for all\/ $x\in\uX$.
For embeddings we may also choose $\bs f$ with\/ $f(X)$ closed
in\/~$\R^n$.
\label{sd4thm12}
\end{thm}

Here is an example in which the condition does not hold.

\begin{ex} $\bR^{\bs k}\t_{\bs 0,\bR^{\bs k},\bs 0}\bs *$ is a
principal d-manifold of virtual dimension 0, with $C^\iy$-scheme
$\ul{\R}^k$, and obstruction bundle $\R^k$. Thus $\bX=\coprod_{k\ge
0}\bR^{\bs k}\t_{\bs 0,\bR^{\bs k},\bs 0}\bs *$ is a d-manifold of
virtual dimension 0, with $C^\iy$-scheme $\uX=\coprod_{k\ge
0}\ul{\R}^k$. Since $T^*_x\uX\cong\R^n$ for
$x\in\R^n\subset\coprod_{k\ge 0}\R^k$, $\dim T^*_x\uX$ realizes all
values $n\ge 0$. Hence there cannot exist immersions or embeddings
$\bs f:\bX\ra\bR^{\bs n}$ for any~$n\ge 0$.
\label{sd4ex4}
\end{ex}

As $x\mapsto\dim T_x^*\uX$ is an upper semicontinuous map $X\ra\N$,
if $\bX$ is compact then $\dim T_x^*\uX$ is bounded above, giving:

\begin{cor} Let\/ $\bX$ be a compact d-manifold. Then there exists
an embedding\/ $\bs f:\bX\ra\bR^{\bs n}$ for some\/~$n\gg 0$.
\label{sd4cor1}
\end{cor}

If a d-manifold $\bX$ can be embedded into a manifold $Y$, we show
in \cite[\S 4.4]{Joyc5} that we can write $\bX$ as the zeroes of a
section of a vector bundle over $Y$ near its image. See
\cite[Prop.~9.5]{Spiv} for the analogue for Spivak's derived
manifolds.

\begin{thm} Suppose $\bX$ is a d-manifold, $Y$ a manifold, and\/
$\bs f:\bX\ra\bY$ an embedding, in the sense of Definition\/
{\rm\ref{sd4def7}}. Then there exist an open subset $V$ in $Y$ with
$\bs f(\bX)\subseteq\bV,$ a vector bundle $E\ra V,$ and\/ $s\in
C^\iy(E)$ fitting into a $2$-Cartesian diagram in $\dSpa\!:$
\begin{equation*}
\xymatrix@C=60pt@R=10pt{ \bX \ar[r]_(0.25){\bs f} \ar[d]^{\bs f}
\drtwocell_{}\omit^{}\omit{^{\eta}}
 & \bV \ar[d]_{\bs 0} \\ \bV \ar[r]^(0.7){\bs s} & \bE.}
\end{equation*}
Here \/ $\bY=F_\Man^\dMan(Y),$ and similarly for $\bV,\bE,\bs s,\bs
0,$ with\/ $0:V\ra E$ the zero section. Hence $\bX$ is equivalent to
the `standard model' d-manifold $\bS_{V,E,s}$ of Example\/
{\rm\ref{sd4ex1},} and is a principal d-manifold.
\label{sd4thm13}
\end{thm}

Combining Theorems \ref{sd4thm12} and \ref{sd4thm13}, Lemma
\ref{sd4lem1}, and Corollary \ref{sd4cor1} yields:

\begin{cor} Let\/ $\bX$ be a d-manifold. Then $\bX$ is a principal
d-manifold if and only if\/ $\dim T^*_x\uX$ is bounded above for
all\/ $x\in\uX$. In particular, if\/ $\bX$ is compact, then $\bX$ is
principal.
\label{sd4cor2}
\end{cor}

Corollary \ref{sd4cor2} suggests that most interesting d-manifolds
are principal, in a similar way to most interesting $C^\iy$-schemes
being affine in Remark \ref{sd2rem1}(ii). Example \ref{sd4ex4} gives
a d-manifold which is not principal.

\subsection{Orientations on d-manifolds}
\label{sd48}

Let $X$ be an $n$-manifold. Then $T^*X$ is a rank $n$ vector bundle
on $X$, so its top exterior power $\La^nT^*X$ is a line bundle (rank
1 vector bundle) on $X$. In algebraic geometry, $\La^nT^*X$ would be
called the canonical bundle of $X$. We define an {\it orientation\/}
$\om$ on $X$ to be an {\it orientation on the fibres of\/}
$\La^nT^*X$. That is, $\om$ is an equivalence class $[\tau]$ of
isomorphisms of line bundles $\tau:O_X\ra\La^nT^*X$, where $O_X$ is
the trivial line bundle $\R\t X\ra X$, and $\tau,\tau'$ are
equivalent if $\tau'=\tau\cdot c$ for some smooth~$c:X\ra(0,\iy)$.

To generalize all this to d-manifolds, we will need a notion of the
`top exterior power' $\cL_{\smash{ (\cE^\bu,\phi)}}$ of a virtual
vector bundle $(\cE^\bu,\phi)$ in \S\ref{sd43}. As the definition is
long and complicated, we will not give it, but just state its
important properties.

\begin{thm} Let\/ $\uX$ be a $C^\iy$-scheme, and\/ $(\cE^\bu,\phi)$
a virtual vector bundle on $\uX$. Then in {\rm\cite[\S 4.5]{Joyc5}}
we define a line bundle (rank\/ $1$ vector bundle)
$\cL_{\smash{(\cE^\bu,\phi)}}$ on $\uX,$ which we call the
\begin{bfseries}orientation line bundle\end{bfseries} of\/
$(\cE^\bu,\phi)$. This satisfies:
\begin{itemize}
\setlength{\itemsep}{0pt}
\setlength{\parsep}{0pt}
\item[{\bf(a)}] Suppose $\cE^1,\cE^2$ are vector bundles on $\uX$ with
ranks $k_1,k_2,$ and\/ $\phi:\cE^1\ra\cE^2$ is a morphism. Then
$(\cE^\bu,\phi)$ is a virtual vector bundle of rank\/ $k_2-k_1,$
and there is a canonical isomorphism\/~$\cL_{\smash{
(\cE^\bu,\phi)}}\cong \La^{k_1}(\cE^1)^*\ot\La^{k_2}\cE^2$.
\item[{\bf(b)}] Let\/ $f^\bu: (\cE^\bu,\phi)\ra(\cF^\bu,\psi)$
be an equivalence in $\vvect(\uX)$. Then there is a canonical
isomorphism\/ $\cL_{f^\bu}:\cL_{\smash{(\cE^\bu,\phi)}}\ra
\cL_{\smash{ (\cF^\bu,\psi)}}$ in\/~$\qcoh(\uX)$.
\item[{\bf(c)}] If\/ $(\cE^\bu,\phi)\in\vvect(\uX)$
then\/~$\cL_{\id_\phi}=\id_{\cL_{\smash{(\cE^\bu,\phi)}}}:
\cL_{\smash{(\cE^\bu,\phi)}}\ra\cL_{\smash{(\cE^\bu,\phi)}}$.
\item[{\bf(d)}] If\/ $f^\bu:(\cE^\bu,\phi)\ra(\cF^\bu,\psi)$ and\/
$g^\bu:(\cF^\bu,\psi)\ra(\cG^\bu,\xi)$ are equivalences in
$\vvect(\uX)$ then $\cL_{g^\bu\ci f^\bu}=\cL_{g^\bu}\ci\cL_{
f^\bu}:\cL_{\smash{(\cE^\bu,\phi)}}\ra\cL_{\smash{(\cG^\bu,\xi)}}$.
\item[{\bf(e)}] If\/ $f^\bu,g^\bu:(\cE^\bu,\phi)\ra(\cF^\bu,\psi)$
are $2$-isomorphic equivalences in $\vvect(\uX)$ then\/~$\cL_{
f^\bu}=\cL_{g^\bu} :\cL_{\smash{(\cE^\bu,\phi)}}\ra\cL_{\smash{
(\cF^\bu,\psi)}}$.
\item[{\bf(f)}] Let\/ $\uf:\uX\ra\uY$ be a morphism of $C^\iy$-schemes,
and\/ $(\cE^\bu,\phi)\in\vvect(\uY)$. Then there is a canonical
isomorphism\/~$I_{\uf,(\cE^\bu,\phi)}:\uf^*(\cL_{\smash{(\cE^\bu,
\phi)}})\ra \cL_{\smash{\uf^*(\cE^\bu,\phi)}}$.
\end{itemize}
\label{sd4thm14}
\end{thm}

Now we can define orientations on d-manifolds.

\begin{dfn} Let $\bX$ be a d-manifold. Then the virtual cotangent
bundle $T^*\bX$ is a virtual vector bundle on $\uX$ by Proposition
\ref{sd4prop2}(b), so Theorem \ref{sd4thm14} gives a line bundle
$\cL_{T^*\bX}$ on $\uX$. We call $\cL_{T^*\bX}$ the {\it orientation
line bundle\/} of~$\bX$.

An {\it orientation\/} $\om$ on $\bX$ is an orientation on
$\cL_{T^*\bX}$. That is, $\om$ is an equivalence class $[\tau]$ of
isomorphisms $\tau:\OX\ra\cL_{T^*\bX}$ in $\qcoh(\uX)$, where
$\tau,\tau'$ are equivalent if they are proportional by a smooth
positive function on~$\uX$.

If $\om=[\tau]$ is an orientation on $\bX$, the {\it opposite
orientation\/} is $-\om=[-\tau]$, which changes the sign of the
isomorphism $\tau:\O_X\ra\cL_{T^*\bX}$. When we refer to $\bX$ as an
oriented d-manifold, $-\bX$ will mean $\bX$ with the opposite
orientation, that is, $\bX$ is short for $(\bX,\om)$ and $-\bX$ is
short for~$(\bX,-\om)$.
\label{sd4def9}
\end{dfn}

\begin{ex}{\bf(a)} Let $X$ be an $n$-manifold, and $\bX=F_\Man^\dMan(X)$
the associated d-manifold. Then $\uX=F_\Man^\CSch(X)$, $\EX=0$ and
$\FX=T^*\uX$. So $\EX,\FX$ are vector bundles of ranks $0,n$. As
$\La^0\EX\cong\O_X$, Theorem \ref{sd4thm14}(a) gives a canonical
isomorphism $\cL_{T^*\bX}\cong\La^nT^*\uX$. That is, $\cL_{T^*\bX}$
is isomorphic to the lift to $C^\iy$-schemes of the line bundle
$\La^nT^*X$ on the manifold~$X$.

As above, an orientation on $X$ is an orientation on the line bundle
$\La^nT^*X$. Hence orientations on the d-manifold
$\bX=F_\Man^\dMan(X)$ in the sense of Definition \ref{sd4def9} are
equivalent to orientations on the manifold $X$ in the usual sense.
\smallskip

\noindent{\bf(b)} Let $V$ be an $n$-manifold, $E\ra V$ a vector
bundle of rank $k$, and $s\in C^\iy(E)$. Then Example \ref{sd4ex1}
defines a `standard model' principal d-manifold $\bS=\bS_{V,E,s}$,
which has $\ES\cong\uE^*\vert_\uS$, $\FS\cong T^*\uV\vert_\uS$,
where $\uE,T^*\uV$ are the lifts of the vector bundles $E,T^*V$ on
$V$ to $\uV$. Hence $\ES,\FS$ are vector bundles on $\uS_{V,E,s}$ of
ranks $k,n$, so Theorem \ref{sd4thm14}(a) gives an
isomorphism~$\cL_{T^*\bS_{V,E,s}}\cong(\La^k\uE\ot
\La^nT^*\uV)\vert_\uS$.

Thus $\cL_{T^*\bS_{V,E,s}}$ is the lift to $\uS_{V,E,s}$ of the line
bundle $\La^kE\ot\La^nT^*V$ over the manifold $V$. Therefore we may
induce an orientation on the d-manifold $\bS_{V,E,s}$ from an
orientation on the line bundle $\La^kE\ot\La^nT^*V$ over $V$.
Equivalently, we can induce an orientation on $\bS_{V,E,s}$ from an
orientation on the total space of the vector bundle $E^*$ over $V$,
or from an orientation on the total space of~$E$.
\label{sd4ex5}
\end{ex}

We can construct orientations on d-transverse fibre products of
oriented d-manifolds. Note that \eq{sd4eq9} depends on an {\it
orientation convention\/}: a different choice would change
\eq{sd4eq9} by a sign depending on $\vdim\bX,\vdim\bY,\vdim\bZ$. Our
conventions follow those of Fukaya et al.\ \cite[\S 8.2]{FOOO} for
Kuranishi spaces.

\begin{thm} Work in the situation of Theorem\/ {\rm\ref{sd4thm7},} so
that\/ $\bW,\bX,\bY,\bZ$ are d-manifolds with\/ $\bW=\bX\t_{\bs
g,\bZ,\bs h}\bY$ for $\bs g,\bs h$ d-transverse, where $\bs
e:\bW\ra\bX,$ $\bs f:\bW\ra\bY$ are the projections. Then we have
orientation line bundles $\cL_{T^*\bW},\ldots,\cL_{T^*\bZ}$ on
$\uW,\ldots,\uZ,$ so\/ $\cL_{T^*\bW},\ue^*(\cL_{T^*\bX}),
\uf^*(\cL_{T^*\bY}),(\ug\ci\ue)^*(\cL_{T^*\bZ})$ are line bundles on
$\uW$. With a suitable choice of orientation convention, there is a
canonical isomorphism
\e
\Phi:\cL_{T^*\bW}\longra\ue^*(\cL_{T^*\bX})\ot_{\O_W}\uf^*
(\cL_{T^*\bY})\ot_{\O_W}(\ug\ci\ue)^*(\cL_{T^*\bZ})^*.
\label{sd4eq9}
\e

Hence, if\/ $\bX,\bY,\bZ$ are oriented d-manifolds, then $\bW$ also
has a natural orientation, since trivializations of\/
$\cL_{T^*\bX},\cL_{T^*\bY},\cL_{T^*\bZ}$ induce a trivialization
of\/ $\cL_{T^*\bW}$ by\/~\eq{sd4eq9}.
\label{sd4thm15}
\end{thm}

Fibre products have natural commutativity and associativity
properties. When we include orientations, the orientations differ by
some sign. Here is an analogue of results of Fukaya et al.\
\cite[Lem.~8.2.3]{FOOO} for Kuranishi spaces.

\begin{prop} Suppose $\bV,\ldots,\bZ$ are oriented d-manifolds,
$\bs e,\ldots,\bs h$ are $1$-morphisms, and all fibre products below
are d-transverse. Then the following hold, in oriented d-manifolds:
\begin{itemize}
\setlength{\itemsep}{0pt}
\setlength{\parsep}{0pt}
\item[{\bf(a)}] For\/ $\bs g:\bX\ra\bZ$ and\/
$\bs h:\bY\ra\bZ$ we have
\begin{equation*}
\bX\t_{\bs g,\bZ,\bs h}\bY\simeq(-1)^{(\vdim \bX-\vdim\bZ)(\vdim
\bY-\vdim \bZ)}\bY\t_{\bs h,\bZ,\bs g}\bX.
\end{equation*}
In particular, when $\bZ=\bs{*}$ so that\/
$\bX\t_\bZ\bY=\bX\t\bY$ we have
\begin{equation*}
\bX\t\bY\simeq(-1)^{\vdim \bX\vdim \bY}\bY\t\bX.
\end{equation*}
\item[{\bf(b)}] For\/ $\bs e:\bV\ra\bY,$ $\bs f:\bW\ra\bY,$ $\bs
g:\bW\ra\bZ,$ and\/ $\bs h:\bX\ra\bZ$ we have
\begin{equation*}
\bV\t_{\bs e,\bY,\bs f\ci\bs\pi_\bW}\bigl(\bW\t_{\bs g,\bZ,\bs
h}\bX\bigr)\simeq \bigl(\bV\t_{\bs e,\bY,\bs f}\bW\bigr)\t_{\bs
g\ci\bs\pi_\bW,\bZ,\bs h}\bX.
\end{equation*}
\item[{\bf(c)}] For\/ $\bs e:\bV\ra\bY,$ $\bs f:\bV\ra\bZ,$
$\bs g:\bW\ra\bY,$ and\/ $\bs h:\bX\ra\bZ$ we have
\begin{align*}
&\bV\t_{(\bs e,\bs f),\bY\t\bZ,\bs g\t\bs h}(\bW\t\bX)\simeq \\
&\quad(-1)^{\vdim\bZ(\vdim\bY+\vdim\bW)} (\bV\t_{\bs e,\bY,\bs
g}\bW)\t_{\bs f\ci\bs\pi_\bV,\bZ,\bs h}\bX.
\end{align*}
\end{itemize}
\label{sd4prop4}
\end{prop}

\subsection{D-manifolds with boundary and corners, d-orbifolds}
\label{sd49}

For brevity, this section will give much less detail than
\S\ref{sd41}--\S\ref{sd48}. So far we have discussed only manifolds
{\it without boundary\/} (locally modelled on $\R^n$). One can also
consider {\it manifolds with boundary\/} (locally modelled on
$[0,\iy)\t\R^{n-1}$) and {\it manifolds with corners} (locally
modelled on $[0,\iy)^k\t\R^{n-k}$). In \cite{Joyc2} the author
studied manifolds with boundary and with corners, giving a new
definition of {\it smooth map\/} $f:X\ra Y$ between manifolds with
corners $X,Y$, satisfying extra conditions over $\pd^kX,\pd^lY$.
This yields categories $\Manb,\Manc$ of manifolds with boundary and
with corners with good properties {\it as categories}.

In \cite[\S 6--\S 7]{Joyc5}, the author defined 2-categories
$\dSpab,\dSpac$ of {\it d-spaces with boundary\/} and {\it with
corners}, and 2-subcategories $\dManb,\dManc$ of {\it d-manifolds
with boundary\/} and {\it with corners}. Objects in
$\dSpab,\dSpac,\dManb,\dManc$ are quadruples $\rX=(\bX,\bpX,\bs
i_\rX,\om_\rX)$, where $\bX,\bpX$ are d-spaces, and $\bs
i_\rX:\bpX\ra\bX$ is a 1-morphism, such that $\bpX$ is locally
equivalent to a fibre product $\bX\t_{\bs{[0,\iy)}}\bs{*}$ in
$\dSpa$, in a similar way to Theorem \ref{sd4thm10}(b). This implies
that the `conormal bundle' $\cN_\rX$ of $\bpX$ in $\bX$ is a line
bundle on $\upX$. The final piece of data $\om_\rX$ is an
orientation on $\cN_\rX$, giving a notion of `outward-pointing
normal vectors' to $\pd\rX$ in $\rX$. Here are some properties of
these:

\begin{thm} In\/ {\rm\cite[\S 7]{Joyc5}} we define strict\/
$2$-categories $\dManb,\dManc$ of \begin{bfseries}d-manifolds with
boundary\end{bfseries} and \begin{bfseries}d-manifolds with
corners\end{bfseries}. These have the following properties:
\begin{itemize}
\setlength{\itemsep}{0pt}
\setlength{\parsep}{0pt}
\item[{\bf(a)}] $\dManb$ is a full\/ $2$-subcategory of\/ $\dManc$.
There is a full and faithful $2$-functor\/
$F_\dMan^\dManc:\dMan\hookra\dManc$ whose image is a full\/
$2$-subcategory $\bdMan$ in $\dManb,$ so that\/
$\bdMan\subset\dManb\subset\dManc$.
\item[{\bf(b)}] There are full and faithful\/ $2$-functors
$F_\Manb^\dManb:\Manb\ra\dManb$ and\/
$F_\Manc^\dManc:\Manc\ra\dManc$. We write $\bManb,\bManc$ for
the full\/ $2$-subcategories of objects in $\dManb,\dManc$
equivalent to objects in the images
of\/~$F_\Manb^\dManb,F_\Manc^\dManc$.
\item[{\bf(c)}] Each object\/ $\rX=(\bX,\bpX,\bs i_\rX,\om_\rX)$ in
$\dManb$ or $\dManc$ has a \begin{bfseries}virtual
dimension\end{bfseries} $\vdim\rX\in\Z$. The virtual cotangent
sheaf\/ $T^*\bX$ of the underlying d-space $\bX$ is a virtual
vector bundle on $\uX$ with rank\/~$\vdim\rX$.
\item[{\bf(d)}] Each d-manifold with corners $\rX$ has a
\begin{bfseries}boundary\end{bfseries} $\pd\rX,$ which is another
d-manifold with corners with $\vdim \pd\rX=\vdim\rX-1$. The
d-space $1$-morphism $\bs i_\rX$ in $\rX$ is also a $1$-morphism
$\bs i_\rX:\pd\rX\ra\rX$ in $\dManc$. If\/ $\rX\in\dManb$ then
$\pd\rX\in\bdMan,$ and if\/ $\rX\in\bdMan$ then\/~$\pd\rX=\es$.
\item[{\bf(e)}] Boundaries in $\dManb,\dManc$ have strong functorial
properties. For instance, there is an interesting class of
\begin{bfseries}simple\end{bfseries} $1$-morphisms $\bs f:\rX\ra\rY$
in $\dManb$ and\/ $\dManc,$ which satisfy a discrete condition
broadly saying that $\bs f$ maps $\pd^k\rX$ to $\pd^k\rY$ for
all\/ $k$. These have the property that for all simple $\bs
f:\rX\ra\rY$ there is a unique simple $1$-morphism $\bs
f_-:\pd\rX\ra\pd\rY$ with $\bs f\ci\bs i_\rX=\bs i_\rY\ci\bs
f_-,$ and the following diagram is $2$-Cartesian in $\dManc$
\begin{equation*}
\xymatrix@C=120pt@R=10pt{ \pd\rX \ar[r]_(0.2){\bs f_-}
\ar[d]_{\bs i_\rX} \drtwocell_{}\omit^{}\omit{^{\id_{\bs f\ci\bs
i_\rX}\,\,\,\,\,\,\,\,\,\,\,\,\,\,{}}}
& \pd\rY \ar[d]^{\bs i_\rY} \\
\rX \ar[r]^(0.7){\bs f} & \rY,}
\end{equation*}
so that\/ $\pd\rX\simeq\rX\t_{\bs f,\rY,\bs i_\rY}\pd\rY$ in
$\dManc$. If\/ $\bs f,\bs g:\rX\ra\rY$ are simple\/
$1$-morphisms and $\eta:\bs f\Ra\bs g$ is a $2$-morphism in
$\dManc$ then there is a natural\/ $2$-morphism $\eta_-:\bs
f_-\Ra\bs g_-$ in\/~$\dManc$.
\item[{\bf(f)}] An \begin{bfseries}orientation\end{bfseries} on a
d-manifold with corners $\rX=(\bX,\bpX,\bs i_\rX,\om_\rX)$ is an
orientation on the line bundle $\cL_{T^*\bX}$ on $\uX$. If\/
$\rX$ is an oriented d-manifold with corners, there is a natural
orientation on $\pd\rX,$ constructed using the orientation on
$\rX$ and the data $\om_\rX$ in\/~$\rX$.
\item[{\bf(g)}] Almost all the results of\/
{\rm\S\ref{sd41}--\S\ref{sd48}} on d-manifolds without boundary
extend to d-manifolds with boundary and with corners, with some
changes.
\end{itemize}
\label{sd4thm16}
\end{thm}

One moral of \cite{Joyc2} and \cite[\S 5--\S 7]{Joyc5} is that doing
`things with corners' properly is a great deal more complicated, but
also more interesting, than you would believe if you had not thought
about the issues involved.

\begin{ex}{\bf(i)} Let $\rX$ be the fibre product
$\bs{[0,\iy)}\t_{\bs i,\bR,\bs 0}\bs *$ in $\dManc$, where
$i:[0,\iy)\hookra\R$ is the inclusion. Then $\rX=(\bX,\bpX,\bs
i_\rX,\om_\rX)$ is `a point with point boundary', of virtual
dimension 0, and its boundary $\pd\rX$ is an `obstructed point', a
point with obstruction space $\R$, of virtual dimension~$-1$.

The conormal bundle $\cN_\rX$ of $\bpX$ in $\bX$ is the obstruction
space $\R$ of $\bpX$. In this case, the orientation $\om_\rX$ on
$\cN_\rX$ cannot be determined from $\bX,\bpX,\bs i_\rX$, in fact,
there is an automorphism of $\bX,\bpX,\bs i_\rX$ which reverses the
orientation of $\cN_\rX$. So $\om_\rX$ really is extra data. We
include $\om_\rX$ in the definition of d-manifolds with corners to
ensure that orientations of d-manifolds with corners are
well-behaved. If we omitted $\om_\rX$ from the definition, there
would exist oriented d-manifolds with corners $\rX$ whose boundaries
$\pd\rX$ are not orientable.
\smallskip

\noindent{\bf(ii)} The fibre product $\bs{[0,\iy)}\t_{\bs
i,\bs{[0,\iy)},\bs 0}\bs *$ is a point $\bs *$ without boundary. The
only difference with {\bf(i)} is that we have replaced the target
$\bR$ with $\bs{[0,\iy)}$, adding a boundary. So in a fibre product
$\rW=\rX\t_\rZ\rY$ in $\dManc$, the boundary of $\rZ$ affects the
boundary of $\rW$. This does not happen for fibre products
in~$\Manc$.

\smallskip

\noindent{\bf(iii)} Let $\rX'$ be the fibre product
$\bs{[0,\iy)}\t_{\bs i,\bR,\bs i}\bs{(-\iy,0]}$ in $\dManc$, that
is, the derived intersection of submanifolds $[0,\iy),(-\iy,0]$ in
$\R$. Topologically, $\rX'$ is just the point $\{0\}$, but as a
d-manifold with corners $\rX'$ has virtual dimension 1. The boundary
$\pd\rX'$ is the disjoint union of two copies of $\rX$ in {\bf(i)}.
The $C^\iy$-scheme $\uX$ in $\bX$ is the spectrum of the
$C^\iy$-ring $C^\iy\bigl([0,\iy)^2\bigr)/(x+y)$, which is
infinite-dimensional, although its topological space is a point.
\label{sd4ex6}
\end{ex}

{\it Orbifolds\/} are generalizations of manifolds locally modelled
on $\R^n/G$ for $G$ a finite group. They are related to manifolds as
Deligne--Mumford stacks are related to schemes in algebraic
geometry, and form a 2-category $\Orb$. Lerman \cite{Lerm} surveys
definitions of orbifolds, and explains why $\Orb$ is a 2-category.
As for $\Manb,\Manc$ one can also consider 2-categories of {\it
orbifolds with boundary\/} $\Orbb$ and {\it orbifolds with
corners\/} $\Orbc$, discussed in~\cite[\S 8]{Joyc5}.

In \cite[\S 9 \& \S 11]{Joyc5} we define 2-categories of {\it
d-stacks\/} $\dSta$, {\it d-stacks with boundary\/} $\dStab$ and
{\it d-stacks with corners\/} $\dStac$, which are orbifold versions
of $\dSpa,\dSpab,\dSpac$. Broadly, to go from d-spaces
$\bX=(\uX,\OXp,\EX,\im_X,\jm_X)$ to d-stacks we just replace the
$C^\iy$-scheme $\uX$ by a {\it Deligne--Mumford\/ $C^\iy$-stack\/}
$\cX$, where the theory of the 2-category of $C^\iy$-stacks $\CSta$
is developed in \cite[\S 7--\S 11]{Joyc3} and summarized in \cite[\S
4]{Joyc4}. Then in \cite[\S 10 \& \S 12]{Joyc5} we define
2-categories of {\it d-orbifolds\/} $\dOrb$, {\it d-orbifolds with
boundary\/} $\dOrbb$ and {\it d-orbifolds with corners\/} $\Orbc$,
which are orbifold versions of~$\dMan,\dManb,\dManc$.

One might expect that combining the 2-categories $\Orb$ and $\dMan$
should result in a 3-category $\dOrb$, but in fact a 2-category is
sufficient. For 1-morphisms $\bs f,\bs g:\bcX\ra\bcY$ in $\dOrb$, a
2-morphism $\bs\eta:\bs f\Ra\bs g$ in $\dOrb$ is a pair
$(\eta,\eta')$, where $\eta:f\Ra g$ is a 2-morphism in $\CSta$, and
$\eta': f^*(\cF_\cY)\ra\cE_\cX$ is as for 2-morphisms in $\dMan$.
These $\eta,\eta'$ do not interact very much.

The generalizations to d-orbifolds are mostly straightforward, with
few surprises. Almost all the results of \S\ref{sd3}--\S\ref{sd48},
and Theorem \ref{sd4thm16}, extend to d-stacks and d-orbifolds with
only cosmetic changes. One exception is that the generalizations of
Theorems \ref{sd3thm3} and \ref{sd4thm5} to d-stacks and d-orbifolds
need extra conditions on the $C^\iy$-stack 2-morphism components
$\eta_{ijk}$ of $\bs\eta_{ijk}$ on quadruple overlaps
$\bcX_i\cap\bcX_j\cap\bcX_k \cap\bcX_l$, as in Remark \ref{sd3rem2}.
This is because 2-morphisms $\eta_{ijk}$ in $\CSta$ are discrete,
and cannot be glued using partitions of unity.

\subsection{D-manifold bordism, and virtual cycles}
\label{sd410}

Classical bordism groups $MSO_k(Y)$ were defined by Atiyah
\cite{Atiy} for topological spaces $Y$, using continuous maps
$f:X\ra Y$ for $X$ a compact oriented manifold. Conner \cite[\S
I]{Conn} gives a good introduction. We define bordism $B_k(Y)$ only
for manifolds $Y$, using smooth $f:X\ra Y$, following Conner's {\it
differential bordism groups\/} \cite[\S I.9]{Conn}. By
\cite[Th.~I.9.1]{Conn}, the natural projection $B_k(Y)\ra MSO_k(Y)$
is an isomorphism, so our notion of bordism agrees with the usual
definition.

\begin{dfn} Let $Y$ be a manifold without boundary, and $k\in\Z$.
Consider pairs $(X,f)$, where $X$ is a compact, oriented manifold
without boundary with $\dim X=k$, and $f:X\ra Y$ is a smooth map.
Define an equivalence relation $\sim$ on such pairs by $(X,f)\sim
(X',f')$ if there exists a compact, oriented \ab $(k+1)$-manifold
with boundary $W$, a smooth map $e:W\ra Y$, and a diffeomorphism of
oriented manifolds $j:-X\amalg X'\ra\pd W$, such that $f\amalg
f'=e\ci i_W\ci j$, where $-X$ is $X$ with the opposite orientation.

Write $[X,f]$ for the $\sim$-equivalence class ({\it bordism
class\/}) of a pair $(X,f)$. For each $k\in\Z$, define the $k^{\it
th}$ {\it bordism group\/} $B_k(Y)$ of $Y$ to be the set of all such
bordism classes $[X,f]$ with $\dim X=k$. We give $B_k(Y)$ the
structure of an abelian group, with zero element $0_Y=[\es,\es]$,
and addition given by $[X,f]+[X',f']=[X\amalg X',f\amalg f']$, and
additive inverses $-[X,f]=[-X,f]$.

Define $\Pi_\bo^\hom:B_k(Y)\ra H_k(Y;\Z)$ by
$\Pi_\bo^\hom:[X,f]\mapsto f_*([X])$, where $H_*(-;\Z)$ is singular
homology, and $[X]\in H_k(X;\Z)$ is the fundamental class.

If $Y$ is oriented and of dimension $n$, there is a biadditive,
associative, supercommutative {\it intersection product\/}
$\bu:B_k(Y)\t B_l(Y)\ra B_{k+l-n}(Y)$, such that if $[X,f],[X',f']$
are classes in $B_*(Y)$, with $f,f'$ transverse, then the fibre
product $X\t_{f,Y,f'}X'$ exists as a compact oriented manifold, and
\begin{equation*}
[X,f]\bu[X',f']=[X\t_{f,Y,f'}X',f\ci\pi_X].
\end{equation*}
\label{sd4def10}
\end{dfn}

As in \cite[\S I.5]{Conn}, bordism is a generalized homology theory.
Results of Thom, Wall and others in \cite[\S I.2]{Conn} compute the
bordism groups $B_k(*)$ of the point $*$. This partially determines
the bordism groups of general manifolds $Y$, as there is a spectral
sequence $H_i\bigl(Y;B_j(*)\bigr)\Ra B_{i+j}(Y)$. We define {\it
d-manifold bordism\/} by replacing manifolds $X$ in $[X,f]$ by
d-manifolds~$\bX$:

\begin{dfn} Let $Y$ be a manifold without boundary, and
$k\in\Z$. Consider pairs $(\bX,\bs f)$, where $\bX\in\dMan$ is a
compact, oriented d-manifold without boundary with $\vdim\bX=k$, and
$\bs f:\bX\ra\bY$ is a 1-morphism in $\dMan$,
where~$\bY=F_\Man^\dMan(Y)$.

Define an equivalence relation $\sim$ between such pairs by
$(\bX,\bs f)\sim (\bX',\bs f')$ if there exists a compact, oriented
d-manifold with boundary $\rW$ with $\vdim\rW=k+1$, a 1-morphism
$\bs e:\rW\ra\bY$ in $\dManb$, an equivalence of oriented
d-manifolds $\bs j:-\bX\amalg\bX'\ra\pd \rW$, and a
2-morphism~$\eta:\bs f\amalg\bs f'\Ra\bs e\ci\bs i_\rW\ci \bs j$.

Write $[\bX,\bs f]$ for the $\sim$-equivalence class ({\it d-bordism
class\/}) of a pair $(\bX,\bs f)$. For each $k\in\Z$, define the
$k^{\it th}$ {\it d-manifold bordism group}, or {\it d-bordism
group}, $dB_k(Y)$ of $Y$ to be the set of all such d-bordism classes
$[\bX,\bs f]$ with $\vdim\bX=k$. As for $B_k(Y)$, we give $dB_k(Y)$
the structure of an abelian group, with zero element
$0_Y=[\bs\es,\bs\es]$, addition $[\bX,\bs f]+[\bX',\bs
f']=[\bX\amalg\bX',\bs f\amalg\bs f']$, and additive
inverses~$-[\bX,\bs f]=[-\bX,\bs f]$.

If $Y$ is oriented and of dimension $n$, we define a biadditive,
associative, supercommutative {\it intersection product\/}
$\bu:dB_k(Y)\t dB_l(Y)\ra dB_{k+l-n}(Y)$ by
\begin{equation*}
[\bX,\bs f]\bu[\bX',\bs f']=[\bX\t_{\bs f,\bY,\bs f'}\bX',
\bs f\ci\bs\pi_\bX].
\end{equation*}
Here $\bX\t_{\bs f,\bY,\bs f'}\bX'$ exists as a d-manifold by
Theorem \ref{sd4thm8}(a), and is oriented by Theorem \ref{sd4thm15}.
Note that we do not need to restrict to $[X,f],[X',f']$ with $f,f'$
transverse as in Definition \ref{sd4def10}. Define a morphism
$\Pi_\bo^\dbo:B_k(Y)\ra dB_k(Y)$ for $k\ge 0$ by~$\Pi_\bo^\dbo:
[X,f]\mapsto\bigl[F_\Man^\dMan(X),F_\Man^\dMan(f)\bigr]$.
\label{sd4def11}
\end{dfn}

In \cite[\S 13.2]{Joyc5} we prove that $B_*(Y)$ and $dB_*(Y)$ are
isomorphic. See \cite[Th.~2.6]{Spiv} for the analogous result for
Spivak's derived manifolds.

\begin{thm} For any manifold\/ $Y,$ we have $dB_k(Y)=0$ for $k<0,$
and\/ $\Pi_\bo^\dbo:B_k(Y)\ra dB_k(Y)$ is an isomorphism for $k\ge
0$. When $Y$ is oriented, $\Pi_\bo^\dbo$ identifies the intersection
products $\bu$ on\/ $B_*(Y)$ and\/~$dB_*(Y)$.
\label{sd4thm17}
\end{thm}

Here is the main idea in the proof of Theorem \ref{sd4thm17}. Let
$[\bX,\bs f]\in dB_k(Y)$. By Corollary \ref{sd4cor1} there exists an
embedding $\bs g:\bX\ra\bR^{\bs n}$ for $n\gg 0$. Then the direct
product $(\bs f,\bs g):\bX\ra\bY\t\bR^{\bs n}$ is also an embedding.
Theorem \ref{sd4thm13} shows that there exist an open set
$V\subseteq Y\t\R^n$, a vector bundle $E\ra V$ and $s\in C^\iy(E)$
such that $\bX\simeq\bS_{V,E,s}$. Let $\ti s\in C^\iy(E)$ be a
small, generic perturbation of $s$. As $\ti s$ is generic, the graph
of $\ti s$ in $E$ intersects the zero section transversely. Hence
$\ti X=\ti s^{-1}(0)$ is a $k$-manifold for $k\ge 0$, which is
compact and oriented for $\ti s-s$ small, and $\ti X=\es$ for $k<0$.
Set $\ti f=\pi_Y\vert_{\smash{\ti X}}:\ti X\ra Y$. Then
$\Pi_\bo^\dbo\bigl([\ti X,\ti f]\bigr)=[\bX,\bs f]$, so that
$\Pi_\bo^\dbo$ is surjective. A similar argument for $\rW,\bs e$ in
Definition \ref{sd4def11} shows that $\Pi_\bo^\dbo$ is injective.

By Theorem \ref{sd4thm17}, we may define a projection
$\Pi_\dbo^\hom:dB_k(Y)\ra H_k(Y;\Z)$ for $k\ge 0$ by
$\Pi_\dbo^\hom=\Pi_\bo^\hom\ci(\Pi_\bo^\dbo)^{-1}$. We think of
$\Pi_\dbo^\hom$ as a {\it virtual class map}. Virtual classes (or
virtual cycles, or virtual chains) are used in several areas of
geometry to construct enumerative invariants using moduli spaces. In
algebraic geometry, Behrend and Fantechi \cite{BeFa} construct
virtual classes for schemes with obstruction theories. In symplectic
geometry, there are many versions --- see for example Fukaya et al.\
\cite[\S 6]{FuOn}, \cite[\S A1]{FOOO}, Hofer et al.\ \cite{HWZ2},
and McDuff~\cite{McDu}.

The main message we want to draw from this is that {\it oriented
d-manifolds and d-orbifolds admit virtual classes\/} (or virtual
cycles, or virtual chains, as appropriate). Thus, we can use
d-manifolds and d-orbifolds as the geometric structure on moduli
spaces in enumerative invariants problems such as Gromov--Witten
invariants, Lagrangian Floer cohomology, Donaldson--Thomas
invariants, \ldots, as this structure is strong enough to contain
all the `counting' information.

In future work the author intends to define a virtual chain
construction for d-manifolds and d-orbifolds, expressed in terms of
new (co)homology theories whose (co)chains are built from
d-manifolds or d-orbifolds, as for the `Kuranishi (co)homology'
described in~\cite{Joyc1}.

\subsection{Relation to other classes of spaces in mathematics}
\label{sd411}

In \cite[\S 14]{Joyc5} the author studied the relationships between
d-manifolds and d-orbifolds and other classes of geometric spaces in
the literature. The next theorem summarizes our results:

\begin{thm} We may construct `truncation functors' from various
classes of geometric spaces to d-manifolds and d-orbifolds, as
follows:
\begin{itemize}
\setlength{\itemsep}{0pt}
\setlength{\parsep}{0pt}
\item[{\bf(a)}] There is a functor\/ $\Pi_{\bf BManFS}^{\bf
dMan}:\mathop{\bf BManFS}\ra\Ho(\dMan),$ where $\mathop{\bf
BManFS}$ is a category whose objects are triples $(V,E,s)$ of a
Banach manifold $V,$ Banach vector bundle $E\ra V,$ and smooth
section $s:V\ra E$ whose linearization $\d s\vert_x:T_xV\ra
E\vert_x$ is Fredholm with index $n\in\Z$ for each $x\in V$
with\/ $s\vert_x=0,$ and\/ $\Ho(\dMan)$ is the homotopy category
of the $2$-category of d-manifolds\/~$\dMan$.

There is also an orbifold version $\Pi_{\bf BOrbFS}^{\bf
dOrb}:\Ho(\mathop{\bf BOrbFS})\ra\Ho(\dOrb)$ of this using
Banach orbifolds $V,$ and `corners' versions of both.
\item[{\bf(b)}] There is a functor\/ $\Pi_{\bf MPolFS}^{\bf
dMan}:\mathop{\bf MPolFS}\ra\Ho(\dMan),$ where $\mathop{\bf
MPolFS}$ is a category whose objects are triples $(V,E,s)$ of an
\begin{bfseries}M-polyfold\end{bfseries} without boundary\/ $V$
as in Hofer, Wysocki and Zehnder\/ {\rm \cite[\S 3.3]{HWZ1},} a
fillable strong M-polyfold bundle $E$ over\/ $V$ {\rm \cite[\S
4.3]{HWZ1},} and an sc-smooth Fredholm section $s$ of\/ $E$ {\rm
\cite[\S 4.4]{HWZ1}} whose linearization $\d s\vert_x:T_xV\ra
E\vert_x$ {\rm \cite[\S 4.4]{HWZ1}} has Fredholm index\/
$n\in\Z$ for all\/ $x\in V$ with\/~$s\vert_x=0$.

There is also an orbifold version $\Pi_{\bf PolFS}^{\bf
dOrb}:\Ho(\mathop{\bf PolFS})\ra\Ho(\dOrb)$ of this using
\begin{bfseries}polyfolds\end{bfseries} $V,$ and `corners'
versions of both.
\item[{\bf(c)}] Given a d-orbifold with corners $\bcX,$ we can
construct a \begin{bfseries}Kuranishi space\end{bfseries}
$(X,\ka)$ in the sense of Fukaya, Oh, Ohta and Ono\/
{\rm\cite[\S A]{FOOO},} with the same underlying topological
space $X$. Conversely, given a Kuranishi space $(X,\ka),$ we can
construct a d-orbifold with corners $\bcX'$. Composing the two
constructions, $\bcX$ and\/ $\bcX'$ are equivalent
in\/~$\dOrbc$.

Very roughly speaking, this means that the `categories' of
d-orbifolds with corners, and Kuranishi spaces, are equivalent.
However, Fukaya et al.\ {\rm\cite{FOOO}} do not define morphisms
of Kuranishi spaces, nor even when two Kuranishi spaces are `the
same', so we have no category of Kuranishi spaces.
\item[{\bf(d)}] There is a functor\/ $\Pi_{\bf SchObs}^{\bf
dMan}:\mathop{\bf Sch_\C Obs}\ra\Ho(\dMan),$ where $\mathop{\bf
Sch_\C Obs}$ is a category whose objects are triples
$(X,E^\bu,\phi),$ for $X$ a separated, second countable
$\C$-scheme and\/ $\phi:E^\bu\ra \tau_{\ge -1}(L_X)$ a perfect
obstruction theory on $X$ with constant virtual dimension, in
the sense of Behrend and Fantechi\/ {\rm\cite{BeFa}}. We may
define a natural orientation on $\Pi_{\bf SchObs}^{\bf
dMan}(X,E^\bu,\phi)$ for each\/~$(X,E^\bu,\phi)$.

There is also an orbifold version\/ $\Pi_{\bf StaObs}^{\bf
dOrb}:\Ho(\mathop{\bf Sta_\C Obs})\ra\Ho(\dOrb),$ taking $X$ to
be a Deligne--Mumford\/ $\C$-stack.
\item[{\bf(e)}] There is a functor\/ $\Pi_{\bf QsDSch}^{\bf
dMan}:\Ho(\mathop{\bf QsDSch_\C})\longra\Ho(\dMan),$ where
$\mathop{\bf QsDSch_\C}$ is the $\iy$-category of separated,
second countable, quasi-smooth derived\/ $\C$-schemes $X$ of
constant dimension, as in To\"en and Vezzosi\/
{\rm\cite{Toen,ToVe}}. We may define a natural orientation on
$\Pi_{\bf QsDSch}^{\bf dMan}(X)$ for each\/~$X$.

There is also an orbifold version\/ $\Pi_{\bf QsDSta}^{\bf
dOrb}:\Ho(\mathop{\bf QsDSta_\C})\ra\Ho(\dOrb),$ taking $X$ to
be a derived Deligne--Mumford\/~$\C$-stack.
\item[{\bf(f)}]{\bf (Borisov \cite{Bori})} There is a natural
functor\/ $\Pi_{\bf DerMan}^{\bf dMan}:\Ho(\mathop{\bf
DerMan}^{\bf pd}_{\smash{\bf ft}})\ra\Ho(\dMan_{\bf pr})$ from
the homotopy category of the $\iy$-category $\smash{\mathop {\bf
DerMan}^{\bf pd}_{\smash{\bf ft}}}$ of \begin{bfseries}derived
manifolds\end{bfseries} of finite type with pure dimension, in
the sense of Spivak\/ {\rm\cite{Spiv},} to the homotopy category
of the full\/ $2$-subcategory $\dMan_{\smash{\bf pr}}$ of
principal d-manifolds in\/ $\dMan$. This functor induces a
bijection between isomorphism classes of objects in
$\Ho(\mathop{\bf DerMan}^{\bf pd}_{\smash{\bf ft}})$ and\/
$\Ho(\dMan_{\bf pr})$. It is full, but not faithful. If\/ $[\bs
f]$ is a morphism in $\Ho(\mathop{\bf DerMan}^{\bf
pd}_{\smash{\bf ft}}),$ then $[\bs f]$ is an isomorphism if and
only if\/ $\Pi_{\bf DerMan}^{\bf dMan}([\bs f])$ is an
isomorphism.
\end{itemize}
\label{sd4thm18}
\end{thm}

One moral of Theorem \ref{sd4thm18} is that essentially every
geometric structure on moduli spaces which is used to define
enumerative invariants, either in differential geometry, or in
algebraic geometry over $\C$, has a truncation functor to
d-manifolds or d-orbifolds. Combining Theorem \ref{sd4thm18} with
proofs from the literature of the existence on moduli spaces of the
geometric structures listed in Theorem \ref{sd4thm18}, in \cite[\S
14]{Joyc5} we deduce:

\begin{thm}{\bf(i)} Any solution set of a smooth nonlinear elliptic
equation with fixed topological invariants on a compact manifold
naturally has the structure of a d-manifold, uniquely up to
equivalence in\/~$\dMan$.

For example, let\/ $(M,g),(N,h)$ be Riemannian manifolds, with\/ $M$
compact. Then the family of \begin{bfseries}harmonic
maps\end{bfseries} $f:M\ra N$ is a d-manifold\/ $\bs\cH_{M,N}$
with\/ $\vdim\bs\cH_{M,N}=0$. If\/ $M={\cal S}^1,$ then
$\bs\cH_{M,N}$ is the moduli space of \begin{bfseries}parametrized
closed geodesics\end{bfseries} in\/~$(N,h)$.
\smallskip

\noindent{\bf(ii)} Let\/ $(M,\om)$ be a compact symplectic manifold
of dimension $2n,$ and\/ $J$ an almost complex structure on\/ $M$
compatible with\/ $\om$. For\/ $\be\in H_2(M,\Z)$ and\/ $g,m\ge 0$,
write $\oM_{g,m}(M,J,\be)$ for the moduli space of stable triples
$(\Si,\vec z,u)$ for $\Si$ a genus $g$ prestable Riemann surface
with\/ $m$ marked points $\vec z=(z_1,\ldots,z_m)$ and\/ $u:\Si\ra
M$ a $J$-holomorphic map with\/ $[u(\Si)]=\be$ in\/ $H_2(M,\Z)$.
Using results of Hofer, Wysocki and Zehnder\/ {\rm\cite{HWZ3}}
involving their theory of polyfolds, we can make
$\oM_{g,m}(M,J,\be)$ into a compact, oriented
d-orbifold\/~$\bs\oM_{g,m}(M,J,\be)$.

\smallskip

\noindent{\bf(iii)} Let\/ $(M,\om)$ be a compact symplectic
manifold, $J$ an almost complex structure on\/ $M$ compatible with\/
$\om,$ and\/ $L$ a compact, embedded Lagrangian submanifold in\/
$M$. For\/ $\be\in H_2(M,L;\Z)$ and\/ $k\ge 0,$ write
$\oM_k(M,L,J,\be)$ for the moduli space of
$J$-\begin{bfseries}holomorphic stable maps\end{bfseries} $(\Si,\vec
z,u)$ to $M$ from a prestable holomorphic disc $\Si$ with\/ $k$
boundary marked points $\vec z=(z_1,\ldots,z_k),$ with\/
$u(\pd\Si)\subseteq L$ and\/ $[u(\Si)]=\be$ in $H_2(M,L;\Z)$. Using
results of Fukaya, Oh, Ohta and Ono\/ {\rm\cite[\S 7--\S 8]{FOOO}}
involving their theory of Kuranishi spaces, we can make
$\oM_k(M,L,J,\be)$ into a compact d-orbifold with corners\/
$\bs\oM_k(M,L,J,\be)$. Given a relative spin structure for $(M,L),$
we may define an orientation on\/~$\bs\oM_k(M,L,J,\be)$.
\smallskip

\noindent{\bf(iv)} Let\/ $X$ be a complex projective manifold, and\/
$\oM_{g,m}(X,\be)$ the Deligne--Mumford moduli\/ $\C$-stack of
stable triples $(\Si,\vec z,u)$ for $\Si$ a genus $g$ prestable
Riemann surface with\/ $m$ marked points $\vec z=(z_1,\ldots,z_m)$
and\/ $u:\Si\ra X$ a morphism with\/ $u_*([\Si])=\be\in H_2(X;\Z)$.
Then Behrend\/ {\rm\cite{Behr}} defines a perfect obstruction theory
on $\oM_{g,m}(X,\be),$ so we can make $\oM_{g,m}(X,\be)$ into a
compact, oriented d-orbifold\/~$\bs\oM_{g,m}(X,\be)$.
\smallskip

\noindent{\bf(v)} Let $X$ be a complex algebraic surface, and\/
$\cM$ a stable moduli\/ $\C$-scheme of vector bundles or coherent
sheaves $E$ on $X$ with fixed Chern character. Then Mochizuki
{\rm\cite{Moch}} defines a perfect obstruction theory on $\cM,$ so
we can make $\cM$ into an oriented d-manifold\/~$\bs\cM$.
\smallskip

\noindent{\bf(vi)} Let\/ $X$ be a complex Calabi--Yau $3$-fold or
smooth Fano $3$-fold, and\/ $\cM$ a stable moduli\/ $\C$-scheme of
coherent sheaves $E$ on $X$ with fixed Hilbert polynomial. Then
Thomas\/ {\rm\cite{Thom}} defines a perfect obstruction theory on
$\cM,$ so we can make $\cM$ into an oriented d-manifold\/~$\bs\cM$.
\smallskip

\noindent{\bf(vii)} Let\/ $X$ be a smooth complex projective
$3$-fold, and\/ $\cM$ a moduli\/ $\C$-scheme of `stable PT pairs'
$(C,D)$ in $X,$ where $C\subset X$ is a curve and $D\subset C$ is a
divisor. Then Pandharipande and Thomas\/ {\rm\cite{PaTh}} define a
perfect obstruction theory on $\cM,$ so we can make $\cM$ into a
compact, oriented d-manifold\/~$\bs\cM$.
\smallskip

\noindent{\bf(ix)} Let\/ $X$ be a complex Calabi--Yau $3$-fold,
and\/ $\cM$ a separated moduli\/ $\C$-scheme of simple perfect
complexes in the derived category $D^b\coh(X)$. Then Huybrechts and
Thomas\/ {\rm\cite{HuTh}} define a perfect obstruction theory on
$\cM,$ so we can make $\cM$ into an oriented d-manifold\/~$\bs\cM$.
\label{sd4thm19}
\end{thm}

We can use d-manifolds and d-orbifolds to construct {\it virtual
classes\/} or {\it virtual chains\/} for all these moduli spaces.

\begin{rem} D-manifolds should not be confused with {\it
differential graded manifolds}, or {\it dg-manifolds}. This term is
used in two senses, in algebraic geometry to mean a special kind of
dg-scheme, as in Ciocan-Fontanine and Kapranov
\cite[Def.~2.5.1]{CiKa}, and in differential geometry to mean a
supermanifold with extra structure, as in Cattaneo and Sch\"atz
\cite[Def.~3.6]{CaSc}. In both cases, a dg-manifold $\mathfrak E$ is
roughly the total space of a graded vector bundle $E^\bu$ over a
manifold $V$, with a vector field $Q$ of degree 1
satisfying~$[Q,Q]=0$.

For example, if $E$ is a vector bundle over $V$ and $s\in C^\iy(E)$,
we can make $E$ into a dg-manifold $\mathfrak E$ by giving $E$ the
grading $-1$, and taking $Q$ to be the vector field on $E$
corresponding to $s$. To this $\mathfrak E$ we can associate the
d-manifold $\bS_{V,E,s}$ from Example \ref{sd4ex1}. Note that
$\bS_{V,E,s}$ only knows about an infinitesimal neighbourhood of
$s^{-1}(0)$ in $V$, but $\mathfrak E$ remembers all of~$V,E,s$.
\label{sd4rem2}
\end{rem}

\appendix

\section{Basics of 2-categories}
\label{sdA}

Finally we discuss 2-categories. A good reference is Behrend et
al.~\cite[App.~B]{BEFF}.

\begin{dfn} A ({\it strict\/}) 2-{\it category\/} $\fC$ consists of
a proper class of {\it objects\/} $\Obj(\fC)$, for all
$X,Y\in\Obj(\fC)$ a category $\Hom(X,Y)$, for all $X$ in $\Obj(\fC)$
an object $\id_X$ in $\Hom(X,X)$ called the {\it identity
$1$-morphism}, and for all $X,Y,Z$ in $\Obj(\fC)$ a functor
$\mu_{X,Y,Z}:\Hom(X,Y)\t\Hom(Y,Z)\ra\Hom(X,Z)$. These must satisfy
the {\it identity property}, that
\begin{equation*}
\mu_{X,X,Y}(\id_X,-)\!=\! \mu_{X,Y,Y}(-,\id_Y)\!=\!\id_{\Hom(X,Y)}
\end{equation*}
as functors $\Hom(X,Y)\ra\Hom(X,Y)$, and the {\it associativity
property}, that
\begin{equation*}
\mu_{W,Y,Z}\ci(\mu_{W,X,Y}\t\id_{\Hom(Y,Z)})
=\mu_{W,X,Z}\ab\ci\ab(\id_{\Hom(W,X)}\ab\t\mu_{X,Y,Z})
\end{equation*}
as functors $\Hom(W,X)\t\Hom(X,Y)\t\Hom(Y,Z)\ra\Hom(W,X)$.

Objects $f$ of $\Hom(X,Y)$ are called 1-{\it morphisms}, written
$f:X\ra Y$. For 1-morphisms $f,g:X\ra Y$, morphisms $\eta\in
\Hom_{\Hom(X,Y)}(f,g)$ are called 2-{\it morphisms}, written
$\eta:f\Ra g$. Thus, a 2-category has objects $X$, and two kinds of
morphisms, 1-morphisms $f:X\ra Y$ between objects, and 2-morphisms
$\eta:f\Ra g$ between 1-morphisms.
\label{sdAdef1}
\end{dfn}

There are three kinds of composition in a 2-category, satisfying
various associativity relations. If $f:X\ra Y$ and $g:Y\ra Z$ are
1-morphisms then $\mu_{X,Y,Z}(f,g)$ is the {\it horizontal
composition of\/ $1$-morphisms}, written $g\ci f:X\ra Z$. If
$f,g,h:X\ra Y$ are 1-morphisms and $\eta:f\Ra g$, $\ze:g\Ra h$ are
2-morphisms then composition of $\eta,\ze$ in $\Hom(X,Y)$ gives the
{\it vertical composition of\/ $2$-morphisms} of $\eta,\ze$, written
$\ze\od\eta:f\Ra h$, as a diagram
\e
\begin{gathered}
\xymatrix@C=25pt{ X \rruppertwocell^f{\eta} \rrlowertwocell_h{\ze}
\ar[rr]_(0.35)g && Y & \ar@{~>}[r] && X
\rrtwocell^f_h{{}\,\,\,\,\ze\od\eta\!\!\!\!\!} && Y.}
\end{gathered}
\label{sdAeq1}
\e
And if $f,\ti f:X\ra Y$ and $g,\ti g:Y\ra Z$ are 1-morphisms and
$\eta:f\Ra\ti f$, $\ze:g\Ra\ti g$ are 2-morphisms then
$\mu_{X,Y,Z}(\eta,\ze)$ is the {\it horizontal composition of\/
$2$-morphisms}, written $\ze*\eta:g\ci f\Ra\ti g\ci\ti f$, as a
diagram
\e
\begin{gathered}
\xymatrix@C=20pt{ X \rrtwocell^f_{\ti f}{\eta} && Y
\rrtwocell^g_{\ti g}{\ze} && Z & \ar@{~>}[r] && X \rrtwocell^{g\ci
f}_{\ti g\ci\ti f}{{}\,\,\,\ze*\eta\!\!\!\!\!} && Z. }
\end{gathered}
\label{sdAeq2}
\e
There are also two kinds of identity: {\it identity\/
$1$-morphisms\/} $\id_X:X\ra X$ and {\it identity\/
$2$-morphisms\/}~$\id_f:f\Ra f$.

A basic example is the {\it $2$-category of categories\/}
$\mathfrak{Cat}$, with objects categories $\cC$, 1-morphisms
functors $F:\cC\ra\cD$, and 2-morphisms natural transformations
$\eta:F\Ra G$ for functors $F,G:\cC\ra\cD$. Orbifolds naturally form
a 2-category, as do stacks in algebraic geometry.

In a 2-category $\fC$, there are three notions of when objects $X,Y$
in $\fC$ are `the same': {\it equality\/} $X=Y$, and {\it
isomorphism}, that is we have 1-morphisms $f:X\ra Y$, $g:Y\ra X$
with $g\ci f=\id_X$ and $f\ci g=\id_Y$, and {\it equivalence}, that
is we have 1-morphisms $f:X\ra Y$, $g:Y\ra X$ and 2-isomorphisms
$\eta:g\ci f\Ra\id_X$ and $\ze:f\ci g\Ra\id_Y$. Usually equivalence
is the correct notion.

{\it Commutative diagrams\/} in 2-categories should in general only
commute {\it up to (specified)\/ $2$-isomorphisms}, rather than
strictly. A simple example of a commutative diagram in a 2-category
$\fC$ is
\begin{equation*}
\xymatrix@C=50pt@R=8pt{ & Y \ar[dr]^g \ar@{=>}[d]^\eta \\
X \ar[ur]^f \ar[rr]_h && Z, }
\end{equation*}
which means that $X,Y,Z$ are objects of $\fC$, $f:X\ra Y$, $g:Y\ra
Z$ and $h:X\ra Z$ are 1-morphisms in $\fC$, and $\eta:g\ci f\Ra h$
is a 2-isomorphism.

We define fibre products in 2-categories, following
\cite[Def.~B.13]{BEFF}.

\begin{dfn} Let $\fC$ be a 2-category and $g:X\ra Z$, $h:Y\ra Z$ be
1-morphisms in $\fC$. A {\it fibre product\/} $X\t_ZY$ in $\fC$
consists of an object $W$, 1-morphisms $\pi_X:W\ra X$ and
$\pi_Y:W\ra Y$ (we usually write $e=\pi_X$ and $f=\pi_Y$) and a
2-isomorphism $\eta:g\ci\pi_X\Ra h\ci\pi_Y$ in $\fC$ with the
following universal property: suppose $\pi_X':W'\ra X$ and
$\pi_Y':W'\ra Y$ are 1-morphisms and $\eta':g\ci\pi_X'\Ra
h\ci\pi_Y'$ is a 2-isomorphism in $\fC$. Then there should exist a
1-morphism $b:W'\ra W$ and 2-isomorphisms $\ze_X:\pi_X\ci
b\Ra\pi_X'$, $\ze_Y:\pi_Y\ci b\Ra\pi_Y'$ such that the following
diagram of 2-isomorphisms commutes:
\begin{equation*}
\xymatrix@C=50pt@R=11pt{ g\ci\pi_X\ci b \ar@{=>}[r]_{\eta*\id_b}
\ar@{=>}[d]_{\id_g*\ze_X} & h\ci\pi_Y\ci b
\ar@{=>}[d]^{{}\,\id_h*\ze_Y} \\ g\ci\pi_X'
\ar@{=>}[r]^{\eta'} & h\ci\pi_Y'.}
\end{equation*}
Furthermore, if $\ti b,\ti\ze_X,\ti\ze_Y$ are alternative choices of
$b,\ze_X,\ze_Y$ then there should exist a unique 2-isomorphism
$\th:\ti b\Ra b$ with
\begin{equation*}
\ti\ze_X=\ze_X\od(\id_{\pi_X}*\th)\quad\text{and}\quad
\ti\ze_Y=\ze_Y\od(\id_{\pi_Y}*\th).
\end{equation*}
If a fibre product $X\t_ZY$ in $\fC$ exists then it is unique up to
equivalence.
\label{sdAdef2}
\end{dfn}

Orbifolds, and stacks in algebraic geometry, form 2-categories, and
Definition \ref{sdAdef2} is the right way to define fibre products
of orbifolds or stacks.

\medskip

\noindent{\small\sc The Mathematical Institute, 24-29 St. Giles,
Oxford, OX1 3LB, U.K.}

\noindent{\small\sc E-mail: \tt joyce@maths.ox.ac.uk}

\end{document}